\documentclass[11pt]{article}

\usepackage{amsmath}
\usepackage{amsthm}
\usepackage{_a4}
\usepackage{amsfonts}
\usepackage{amssymb}
\usepackage{amscd}

\begin{document}

%
%

\newcommand{\lb}[1]{\label{#1}}
\newcommand{\G}{\ensuremath{{\cal G}}}
\newcommand{\F}{\ensuremath{{\cal F}}}
\newcommand{\Fh}{\ensuremath{\hat{{\cal F}}}}
\newcommand{\Ah}{\ensuremath{\hat{{\cal A}}}}
\newcommand{\K}{\ensuremath{{\cal K}}}
\newcommand{\Bc}{\ensuremath{{\cal B}}}
\newcommand{\U}{\ensuremath{{\cal U}}}
\newcommand{\M}{\ensuremath{{\cal M}}}
\newcommand{\N}{\ensuremath{{\cal N}}}
\newcommand{\A}{\ensuremath{{\cal A}}}
\newcommand{\dG}{\ensuremath{\hat{{\cal G}}}}
\newcommand{\D}{\ensuremath{{\cal D}}}
\newcommand{\MR}{{\M}_R}
\newcommand{\ML}{{\M}_L}
\newcommand{\B}{\tilde{B}}

\newcommand{\LL}[1][]{{\rm \bf L}^{#1}}
\newcommand{\TT}[1][]{{\rm \bf T}^{#1}}
\newcommand{\RR}[1][]{{\rm \bf R}^{#1}}
\newcommand{\DD}[1][]{{\rm \bf D}^{#1}}
\newcommand{\DDt}[1][]{{\rm \bf \tilde{D}}^{#1}}
\newcommand{\LLt}[1][]{{\rm \bf \tilde{L}}^{#1}}
\newcommand{\Lt}{\tilde{L}}
\newcommand{\Tt}{\tilde{T}}
\newcommand{\MM}[1][]{{\rm \bf M}^{#1}}
%
\newcommand{\Fi}[2][]{\phi^{#1}_{#2}}
\newcommand{\Fii}[2][]{(\phi_{#2}^{-1})^{#1}}
\newcommand{\vi}{\varphi}
\newcommand{\Fib}[2][]{\bar{\phi}^{#1}_{#2}}
\newcommand{\la}{\lambda}
\newcommand{\ep}{\epsilon}
\newcommand{\om}{\omega}
\newcommand{\Om}{\Omega}
\newcommand{\de}{\delta}
\newcommand{\laR}{{\la}_R}
\newcommand{\rhoL}{{\rho}_L}
\newcommand{\Psh}{\hat{\Psi}}

\newcommand{\cop}{\Delta}
\newcommand{\copd}{\Delta_D}
\newcommand{\copf}{\Delta_f}

\newcommand{\me}{m_{(1)}}
\newcommand{\mz}{m_{(2)}}
\newcommand{\meb}{m_{(\bar{1})}}
\newcommand{\mzb}{m_{(\bar{2})}}
\newcommand{\vie}{\vi_{(1)}}
\newcommand{\viz}{\vi_{(2)}}
\newcommand{\vid}{\vi_{(3)}}
\newcommand{\wb}{\bar{w}}
\newcommand{\we}{w_{(1)}}
\renewcommand{\ae}{a_{(1)}}
\newcommand{\wz}{w_{(2)}}
\newcommand{\az}{a_{(2)}}
\newcommand{\web}{w_{(\bar{1})}}
\newcommand{\wzb}{w_{(\bar{2})}}
\newcommand{\wez}{w_{(1,2)}}
\newcommand{\webz}{w_{(\bar{1},2)}}
\newcommand{\wezb}{w_{(1,\bar{2})}}
\newcommand{\weeb}{w_{(1,\bar{1})}}
\newcommand{\wzbz}{w_{(\bar{2},2)}}
\newcommand{\wzbe}{w_{(\bar{2},1)}}
\newcommand{\wzbzb}{w_{(\bar{},\bar{2})}}
\newcommand{\eb}{{(\bar{1})}}
\newcommand{\zb}{{(\bar{2})}}
\newcommand{\Ye}{Y_{(1)}}
\newcommand{\Yz}{Y_{(2)}}
\newcommand{\psie}{\psi_{(1)}}
\newcommand{\psiz}{\psi_{(2)}}
\newcommand{\psid}{\psi_{(3)}}
\newcommand{\Si}{S^{-1}}
\newcommand{\tp}{\otimes}
\newcommand{\bo}{\boxtimes}

\newcommand{\coa}[1][]{\delta_{#1}}
\newcommand{\pfeil}{\longrightarrow}
\newcommand{\e}{{\bf 1}}
\newcommand{\eG}{{\e}_{\G}}
\newcommand{\edG}{{\e}_{\dG}}
\newcommand{\eM}{{\e}_{\M}}
\newcommand{\eA}{{\e}_{\A}}
\newcommand{\eN}{{\e}_{\N}}
\newcommand{\eMR}{{\e}_{\MR}}
\newcommand{\eML}{{\e}_{\ML}}
\newcommand{\id}{{\rm id}}
\newcommand{\idM}{{\id}_{\M}}
\newcommand{\idMR}{{\id}_{\MR}}
\newcommand{\idML}{{\id}_{\ML}}
\newcommand{\idG}{{\id}_{\G}}
\newcommand{\lpa}{\langle}
\newcommand{\rpa}{\rangle}
%
\newcommand{\wir}{\triangleright}
\def\lef{{\,\hbox{$\textstyle\triangleright$}\,}}

\def\re{\lef}
\def\li{{\,\hbox{$\textstyle\triangleleft$}\,}}
\def\reli{\bowtie}
%
\newcommand{\rep}{\makebox{\rm Rep} \,}
\newcommand{\End}{\makebox{\rm End}}
\newcommand{\Co}{\mathbb{C}}

\newcommand{\arr}{\rightharpoonup}
\newcommand{\arl}{\leftharpoonup}

%
\newtheorem{df}{Definition}[section]
\newtheorem{thm}[df]{Theorem}
\newtheorem{lem}[df]{Lemma}
\newtheorem{cor}[df]{Corollary}
\newtheorem{prop}[df]{Proposition}

\numberwithin{equation}{section}

\title{ \bf Doubles of Quasi--Quantum Groups}

\renewcommand{\thefootnote}{\fnsymbol{footnote}}
\author{{\sc Frank Hau{\ss}er\footnotemark[1]}\hspace{1ex} and {\sc Florian
    Nill\footnotemark[2]} 
\\
{\small Freie Universit\"at Berlin,
Institut f\"ur theoretische Physik,}\\
{\small Arnimallee 14, D-14195 Berlin, Germany}}

\maketitle

\footnotetext[1]{
supported by DFG, SFB 288 {\em Differentialgeometrie und
  Quantenphysik},\quad
e-mail: hausser@physik.fu-berlin.de}
\footnotetext[2]{
supported by DFG, SFB 288 {\em Differentialgeometrie und
  Quantenphysik},\quad
e-mail:  nill@physik.fu-berlin.de}

\renewcommand{\thefootnote}{\arabic{footnote}}
\setcounter{footnote}{0}

\begin{abstract}
In [Dr1] Drinfeld showed that any finite dimensional Hopf algebra
$\G$ extends to 
a quasitriangular Hopf algebra $\D(\G)$, the 
quantum double of $\G$. Based on the construction of a so--called {\em
 diagonal crossed product} developed by the 
authors in [HN], we generalize this result to the
case of quasi--Hopf algebras $\G$. As for ordinary Hopf algebras, as a
vector space the
``quasi--quantum double''  $\D(\G)$ is isomorphic to
$\dG\tp\G$, where $\dG$ denotes the dual of $\G$.  We give explicit
formulas for the product, the 
coproduct, the $R$--matrix and the antipode on $\D(\G)$ and prove that
they fulfill Drinfeld's axioms of a quasitriangular quasi--Hopf
algebra. In particular $\D(\G)$ becomes an associative algebra
containing $\G \equiv \edG\tp\G$ as a quasi--Hopf subalgebra. On the
other hand, $\dG\equiv \dG\tp\eG$ is not a subalgebra of $\D(\G)$
unless the coproduct on $\G$ is strictly coassociative. 
It is shown that the category $\rep\D(\G)$ of finite dimensional
representations of $\D(\G)$ coincides with what has been called the double
category of $\G$--modules by S. Majid [M2]. Thus our construction
gives a concrete realization of Majid's abstract definition of
quasi--quantum doubles in terms of a Tannaka--Krein--like reconstruction
procedure. The whole construction is shown to generalize to weak
quasi--Hopf algebras with $\D(\G)$ now being linearly isomorphic to a
subspace of $\dG\tp\G$. 
\end{abstract}

\footnotesize
\tableofcontents

\normalsize
\newpage

%
%
\section{Introduction}
Given a finite dimensional Hopf algebra $\G$ and its dual $\dG$
Drinfeld [Dr1] has introduced the quantum double $\D(\G) \supset \G$
as the universal Hopf algebra extension of $\G$ satisfying 
\begin{enumerate}
\item There exists a unital algebra embedding $D : \dG \pfeil \D(\G)$
  such that $\D(\G)$ is algebraically generated by $\G$ and $D(\dG)$.
\item Let $e_\mu \in \G$ be a basis with dual basis $e^\mu \in
  \dG$. Then $R_D : = e_\mu \tp D(e^\mu) \in \D(\G) \tp \D(\G)$ is
  quasitriangular. 
\end{enumerate}
It follows that as a coalgebra $\D(\G) = \dG^{cop} \tp \G$, where
``cop'' refers to the opposite coproduct. However, when realized on
$\dG\tp\G$, the algebraic structure of $\D(\G)$ becomes more
involved. It has been analyzed in detail by S.\ Majid as a particular
example of his notion of double crossed products, see [M3,M4] and
references therein. The dual version of the quantum double has been
introduced for infinite dimensional compact quantum groups in [PW] as
the mathematical structure underlying the quantum Lorentz group.

During the 90's the quantum double has become of increasing
importance as a quantum symmetry in two--dimensional lattice and
continuum QFT. In continuum theories the quantum double $\D(\G)$ of a
finite group $G$ (i.e.\ $\G = \mathbb{C} G$) has first been applied
(mostly in a twisted version, which we will come back to below) to
describe the symmetry underlying the sector structure of orbifold
models in [DPR]. Quite interestingly, the same structure appears as a
residual generalized ``dyon--symmetry'' in spontanously broken
(2+1)-dimensional Higgs models with a finite unbroken subgroup $G$
[BaWi]. For the role of quantum doubles in integrable field theories see,
e.g.\ [BL].

More recently, in the framework of algebraic QFT, M.\ M\"uger [M\"u]
has also found the double of a finite group $G$ acting as a global
symmetry on a ``disorder--field extension'' $\Fh$ of a massive
2--dimensional field algebra $\F$ with global gauge symmetry $G$. As
opposed to the above cases, in this type of models the
``disorder--part'' $\dG$ of the double is also spontaneously broken,
corresponding to a violation of Haag duality (for double cones) for
the $\D(\G)$--invariant observable algebra $\A\subset \Fh$. The Haag
dual extension $\Ah \supset \A$ is then recovered as the invariant
subalgebra of $\Fh$ under the unbroken symmetry $G$.

On the lattice, related but prior to M\"uger's work, the double of a
finite group $G$ has been realized by  
K.\ Szlach\'anyi and P.\ Vecserny\'es as a symmetry realized on the
order$\times$disorder field 
algebra of a $G$--spin quantum chain [SzV]. Since for $G =
\mathbb{Z}_N$ the double coincides with (the group algebra of) $\mathbb{Z}_N
\times \mathbb{Z}_N$, this generalizes the well known
order$\times$disorder symmetry of abelian $G$--spin models. This
investigation has 
been substantially extended to arbitrary finite dimensional
$C^*$--Hopf algebras $\G$ in [NSz], where the authors show that such
``Hopf spin models'' always have $\D(\G)$ as a universal localized
cosymmetry. This means that under the assumption of a Haag
dual vacuum representation (i.e.\ absence of spontanous symmetry
breaking) the full superselection structure
of these models is precisely created by the irreducible
representations of $\D(\G)$.
The formulation of [NSz] also allowed for a generalization of duality
transformations to the non--commutative and non--cocommutative setting.

As it has turned out meanwhile, very much related results have been
obtained independently for lattice current algebras on finite periodic
lattice chains by A.\ Alekseev et al.\ [AFFS]. For these models the
authors have completely determined the representation category,
showing that it is in one-to-one correspondence with $\rep \K_1$,
where $\K_1$ is the algebra living on a minimal loop consisting of one
site and one link biting into its own tail. Using the braided--group
theory of [M5] (see also [M4]), it has been realized by one of us [N1],
that $\K_1$ is in fact again isomorphic to a quantum double $\D(\G)$. Also,
requiring $\G$ to be a modular Hopf algebra as in [AFFS], the Hopf
spin model of [NSz] has been shown in [N1] to be isomorphic to the
lattice current algebra of [AFFS] by a local transformation of the
generators\footnote{More generally, even without these assumptions on $\G$,
periodic Hopf spin chains are meanwhile known to be isomorphic to
$\D(\G)\tp \text{Mat}(N_L)$, where $N_L \in \mathbb{N}$ depends on the
length of the loop $L$ [Sz], thus explaining the representation theory
of [AFFS].}.\\

As a common feature of all these models we emphasize that under
the quantum physical requirement of positivity they only give rise to
quantum symmetries with integer $q$--dimensions [N2]\footnote{
This result is frequently ignored in the literature and relates to the
finite dimensionality of $\G$. By a Perron--Frobenius argument it also
applies to the twisted double of [DPR].}. Thus, to construct ``rational''
models with a finite sector 
theory and non--integer dimensions one is inevitably forced to depart
from ordinary Hopf algebras $\G$. Here, the most fashionable
candidates are the truncated semisimple versions of the
$q$--deformations $\U_q(\mathfrak{g})$, $\mathfrak{g}$ a simple Lie
algebra, at roots of unity, $q^N = 1$.
Also, since lattice current algebras have been invented as regularized
verions of WZNW--models [AFFS, AFSV, AFS, ByS, Fa, FG], they should
eventually be studied at roots of unity.

 Following G.\ Mack and V.\
Schomerus [MS], truncated quantum groups at $q^N=1$ have to be
described as {\em weak quasi--Hopf algebras} in 
the sense of Drinfeld [Dr2], with the additional feature $\cop(\e)
\neq \e\tp\e$, where $\cop: \G \pfeil \G\tp\G$ denotes the
coproduct\footnote{The twisted double of [DPR] is also a quasi--Hopf
  algebra but it still satisfies $\cop(\e) = \e\tp\e$.}.

To formulate lattice current algebras at roots of unity one may now
combine the methods of [AFFS] with those developped by [AGS,AS] for
lattice Chern--Simons theories. However, it remains unclear whether and
how for $q=$root  of unity the structural 
results of [AFFS] survive the truncation to the semi-simple 
(``physical'') quotients. Similarly, the generalizations of the 
model, the methods and the results of [NSz] to weak quasi 
quantum groups are by no means obvious. In particular one would 
like to know whether and in what sense in such models universal 
localized cosymmetries $\rho:\A\pfeil\A\tp\G$ still provide 
coactions and whether $\G$ would still be (an analogue of) a 
quantum double of a quasi-Hopf algebra.

In fact, a definition of a quantum double $\D(\G)$ for quasi--Hopf
algebras $\G$ has recently been proposed by S.\ Majid
[M2]. Unfortunately this has only been done in form of an implicit
Tannaka--Krein reconstruction procedure, which makes it hard to
identify this algebra in terms of generators and relations in concrete
models.\\

In [HN] we have started a program where we generalize standard notions
of Hopf algebra theory (like coactions and crossed products) to (weak)
quasi--Hopf algebras and apply them to quantum chains based on weak
quasi--quantum groups in the spirit of [NSz, AFFS]. As a central
mathematical structure underlying these constructions we have
developed the concept of a {\em diagonal crossed product} by the
dual $\dG$ of a (weak) quasi--quantum group $\G$. In
this way we have obtained as one of our main nontrivial examples an
explicit algebraic definition of the double $\D(\G)$. We have shown
that, as for 
ordinary Hopf algebras, $\D(\G)$ may be realized as a new
quasi--bialgebra structure on $\dG\tp\G$ \footnote{This has also been
announced in [M2].} (or, in the weak
  case, a certain subspace thereof) containing $\G\equiv\edG\tp\G$ as a
sub-bialgebra. Generalizing the results of [NSz, AFFS] we have also
constructed the above lattice models for weak quasi--Hopf algebras
$\G$ and established that they always admit localized coactions of
$\D(\G)$ in the sense of [NSz].

In this work we extend our analysis of $\D(\G)$ by proving that it is
always a quasitriangular quasi--Hopf algebra, which is weak if and
only if $\cop(\eG) \neq \eG\tp\eG$. Our main results are summarized by
the following
%
%
\vspace{0.25cm}
  \begin{sloppypar} \noindent{\bf Theorem A} {\it 
 Let $(\G,\cop,\phi,S)$ be a finite dimensional quasi--Hopf algebra with
 coproduct\\
  $\cop: \G \pfeil \G\tp\G$, reassociator $\phi \in
 \G^{\tp^3}$ and invertible antipode $S$. Assume $\D(\G)$ to be 
  a quasi--Hopf algebra extension $\D(\G) \supset \G$ satisfying
\begin{itemize}
\item[(i)] There exists a linear map $D: \dG \pfeil \D(\G)$ such that
  $\D(\G)$ is algebraically generated by $\G$ and $D(\dG)$ 
\item[(ii)] $R_D : = \sum_\mu e_\mu \tp D(e^\mu) \in \D(\G)\tp \D(\G)$
  is quasitriangular.
\item[(iii)] If $\tilde{\D}\supset \G$ and $\tilde{D}: \dG \pfeil
  \tilde{\D}$ have the same properties, then their exists a 
  bialgebra homomorphism $f: \D(\G) \pfeil \tilde{\D}$ restricting to
  the identity on $\G$ and satisfying $f \circ D = \tilde{D}$.
\end{itemize}
Then $\D(\G)$ exists uniquely up to equivalence and 
the map $\mu: \dG \tp \G \pfeil \D(\G)$ given by\footnote{
Here $\alpha \in \G$ is one of the two structural elements appearing
in Drinfeld's antipode axioms, see Sect.\ 2.1.}
\begin{equation}
  \label{thA}
  \mu(\vi\tp a) = (\id \tp \vi_{(1)})(q_\rho)\, D(\vi_{(2)}),\quad
  \text{where}\quad q_\rho =  \phi^{1} \tp \Si(\alpha\phi^{3})\phi^{2}
  \in \G\tp\G 
\end{equation}
provides a linear bijection.
}
\end{sloppypar} \vspace{0.25cm}  
Theorem A will be proven at the end of Section 3.4. We will also have
a generalization to weak quasi--Hopf algebras, which is stated as
Theorem~B in Section~4. 

The major achievement of Theorem~A in
comparison with [HN] consists in the construction of the antipode on
$\D(\G)$ To this end, as a central technical result we establish 
a formula for $(S\tp S)(R)$ and the
relations between $R^{-1}$, $(S\tp\id)(R)$ and $(\id\tp\Si)(R)$ for a
quasitriangular $R \in \G\tp\G$ in any quasi--Hopf algebra
$\G$. Recall, that in ordinary Hopf algebras the last three quantities
coincide and therefore $(S\tp S)(R) = R$. 

To prove these results we combine the methods of [HN] with the very efficient
graphical calculus developed by [RT,T,AC]. This will also allow to
give nice intuitive interpretations of many of our almost untraceable
identities derived in [HN]. In fact, without this graphical machinery
we would have been lost in proving or even only trying to guess these
formulas. In particular, a purely algebraic proof of the formulas for
$R^{-1}$ and $(S\tp S)(R)$ in Theorem~2.1 would most likely be
unreadable and therefore also untrustworthy. This is why we think it
worthwhile to put more emphasis on this graphical technique in the
present paper. \\

We start in Section~2.1 with shortly reviewing Drinfeld's theory of
quasi--Hopf algebras and introduce our graphical conventions in
Section~2.2. In Section~2.3 we derive our main formulas for $R^{-1}$
and $(S\tp S)(R)$ for any quasi--triangular $R \in \G\tp\G$. In
Section~3.1 we review our construction [HN] of the double $\D(\G)$ as
an associative algebra on the vector space $\dG\tp\G$. In Section~3.2 we
reformulate this construction in the spirit of [N1] in terms of the
universal $\cop$--flip operator $\DD \in \G\tp\D(\G)$. Section~3.3
roughly sketches, how the double may also be realized on the vector
space $\G\tp\dG$. In Section~3.4 we establish the quasi--triangular
quasi--Hopf structure of $\D(\G)$ and prove Theorem~A. Finally, in
Section~3.5 we identify the category $\rep\D(\G)$ as the double of the
category $\rep \G$ in the sense of [M2], thus proving that our
construction of $\D(\G)$ provides a concrete realization of Majid's
Tannaka--Krein like reconstruction procedure. In Section~4 we
generalize our results to weak quasi--Hopf algebras $\G$. As an
application we discuss the twisted double $\D^{\om}(G)$ of [DPR] in
Appendix~A and generalize the results of [N1] on the relation with the
monodromy algebras of [AGS,AS] in Appendix~B.\\

Throughout, all linear spaces are assumed finite dimensional over the
field $\mathbb{C}$. We will use standard Hopf algebra notations, see
e.g. [A,Sw,K,M3]. By an extension $\Bc\supset \A$ of algebras we always
mean a unital injective algebra morphism $\A\pfeil \Bc$. Two
extensions $\Bc_1\supset\A$ and $\Bc_2\supset \A$ are called equivalent,
if there exists an isomorphism of algebras $\Bc_1 \cong\Bc_2$ restricting
to the identity on $\A$.

\section{Quasitriangular quasi--Hopf algebras}

\subsection{Basic definitions and properties}
In this subsection we review the basic 
definitions and properties of quasitriangular
quasi--Hopf algebras as introduced by
Drinfeld [Dr2], where the interested reader will find a
more detailed discussion.

A {\it quasi-bialgebra} $(\G,\cop,\ep,\phi)$ is an associative algebra
$\G$ with 
unity, algebra morphisms $\cop:\, \G\pfeil \G\tp\G$ and $\ep: \,
\G\pfeil\mathbb{C}$, and an invertible element $\phi \in
\G\tp\G\tp\G$, such that
\begin{gather}
  \lb{eq11}
   (\id\tp\cop)(\cop(a)) \phi = \phi(\cop\tp\id)(\cop(a)), \quad
   a\in\G \\
   \lb{eq12}
   (\id\tp\id\tp\cop)(\phi)(\cop\tp\id\tp\id)(\phi) =
   (\e\tp\phi)(\id\tp\cop\tp\id)(\phi)(\phi\tp\e), \\
   \lb{eq13}
   (\ep\tp\id)\circ\cop = \id = (\id\tp\ep)\circ\cop , \\
   \lb{eq14}
   (\id\tp\ep\tp\id)(\phi)=\e\tp\e
\end{gather}
The map $\cop$ is called the coproduct and $\ep$ the counit. A
coproduct with the above properties is called {\em quasi-coassociative} and
the element $\phi$ will be called the {\em reassociator}. The identities
(\ref{eq12}) and (\ref{eq14}) together imply 
\begin{equation}
  \lb{eq15}
  (\ep\tp\id\tp\id)(\phi) = (\id\tp\id\tp\ep)(\phi)= \e\tp\e.
\end{equation}

Let us briefly recall some of the main consequences of these
definitions for the representation theory of $\G$. Let $\rep\G$ be the
category of finite dimensional representations of $\G$, i.e. of pairs
$(\pi_V,V)$, where $V$ is a finite dimensional vector space and $\pi_V
:\, \G\pfeil \End_\Co (V)$ is a unital algebra morphism. We will also use the
equivalent notion of a $\G$--module $V$ with multiplication 
$g \cdot v \equiv \pi_V(g) v$. Given two pairs 
$(\pi_V,V),(\pi_U,U)$, the coproduct allows for the definition of a
tensor product $(\pi_{V\tp U},V\tp U)$ by setting
  $\pi_{V\tp U} = (\pi_V\tp\pi_U) \circ \cop$.
The counit defines a one dimensional representation. Equation
(\ref{eq13}) says, that this representation is a left and right unit
with respect to the tensor product, and (\ref{eq11}) says that given three
representations $(\pi_U,\pi_V,\pi_W)$, then $\pi_{(U\tp V)\tp W} \cong 
\pi_{U\tp (V\tp W)}$ with intertwiner $\phi_{UVW} = (\pi_U \tp\pi_V
\tp\pi_W)(\phi)$. 

The meaning of (\ref{eq12}) is the commutativity of
the pentagon
\begin{equation}  
\label{pent}
  \unitlength0.5cm
  \begin{picture}(30,5)
  \put(3.5,4){\makebox(0,0){$((U\tp V)\tp W)\tp X$}}
  \put(13.5,4){\makebox(0,0){$(U\tp V)\tp(W\tp X)$}}
  \put(23.5,4){\makebox(0,0){$U\tp (V\tp(W\tp X))$}}
  \put(8.5,0){\makebox(0,0){$(U \tp (V \tp W)) \tp X$}}
  \put(18.5,0){\makebox(0,0){$U\tp ((V \tp W) \tp X)\quad ,$}}
  \put(12.5,0){\vector(1,0){2}}
    \put(7.5,4){\vector(1,0){2}}
   \put(17.5,4){\vector(1,0){2}}
    \put(3.5,3){\vector(2,-1){4}}
    \put(18.5,1){\vector(2,1){4}}
  \end{picture}
\end{equation}
where the arrows stand for the corresponding rebracketing intertwiners.
For example the first one is given by
$ (\pi_{U\tp V}\tp \pi_W \tp \pi_X)(\phi) =
(\pi_U\tp\pi_V\tp\pi_W\tp\pi_X)\big((\cop\tp\id\tp\id)(\phi)\big)$. 
The diagram (\ref{pent}) explains the name pentagon identity for equation 
(\ref{eq12}). The importance of axiom (\ref{eq12}) lies in the fact, 
that in any tensor product representation the intertwiner
connecting two different bracket conventions is given by a suitable
product of $\phi$'s, as in (\ref{pent}). The pentagon identity then
guarantees, that this intertwiner is
independent of the chosen sequence of intermediate rebracketings.
This is known as Mac Lanes coherence theorem [ML].

A quasi--bialgebra $\G$ is called {\it quasi-Hopf algebra}, if there
is a linear 
antimorphism $S: \, \G \rightarrow \G$ and elements $\alpha, \beta \in
\G$ satisfying (for all $a\in\G$)
\begin{align}
  \label{eq16}
  &\sum_i S(a^i_{(1)})\alpha a^i_{(2)} = \alpha \ep (a), \,\,\quad
  \sum_i a^i_{(1)} \beta S(a^i_{(2)}) = \beta \ep(a)\\
  \label{eq17}
 & \sum_j X^j \beta S(Y^j)\alpha Z^j = 1 = \sum_j S(P^j)\alpha Q^j
  \beta S(R^j).
\end{align}
Here and throughout we use the notation 
$\sum_i a^i_{(1)}\tp a^i_{(2)} = \cop(a)$ and
\begin{equation}
  \label{eq17a}
\phi = \sum_j X^j\tp Y^j \tp Z^j; \quad\phi^{-1}= \sum_j P^j\tp Q^j\tp R^j.
\end{equation}
To simplify the notation, we will in the following also
frequently suppress the summation symbol and write $\phi = X^i\tp Y^i
\tp Z^i$, $\cop(a) = a_{(1)}\tp a_{(2)}$,
etc. The map $S$ is called an antipode. 
We will also always suppose that $S$ is invertible. Note that as opposed to
ordinary Hopf algebras, an antipode is not uniquely determined,
provided it exists. The antipode allows to define 
the (left) dual representation $({}^*\pi,{}^*V)$ of $(\pi,V)$, where
${}^*V$ is the 
dual space of $V$, by ${}^*\pi(a) = \pi(S(a))^t$, the superscript $t$
denoting the transposed map. Analogously one defines
a right dual representation $(\pi^*,V^*)$, where $V^*\equiv {}^*V$ and
$\pi^*(a) = \pi(\Si(a))^t$.  

A quasi-Hopf algebra $\G$ is called {\it quasitriangular}, if there exists
an invertible element $R\in\G\tp\G$, such that
\begin{align}
   \label{eq18}
  \cop^{op}(a) R &= R \cop(a),\quad a\in\G \\
   \label{eq19}
  (\cop\tp\id)(R) &= \Fi[312]{} \,R^{13}\, \Fii[132]{}\, R^{23}\, \Fi{} \\
   \label{eq110}
  (\id\tp\cop)(R) &= \Fii[231]{} \,R^{13}\,\Fi[213]{}\, R^{12}\, \phi^{-1},
\end{align}
where we use the following notation: If $\psi = \sum_i \psi_i^1
\tp\dots \psi_i^n \,\in\G^{\tp^m}$, then, for $m\le n$, $\psi^{n_1 n_2
  \dots n_m} \in 
\G^{\tp^n}$ denotes the element of $\G^{\tp^n}$ having $\psi_i^k$
in the ${n_k}^{\rm th}$ slot and $\e$ in the remaining ones. The element
$R$ is called the R-matrix. The above
relations imply the quasi-Yang-Baxter equation
\begin{equation}
  \label{eq111}
   R^{12}\,\Fi[312]{}\, R^{13}\, \Fii[132]{}\, R^{23} \,\Fi{} =
   \Fi[321]{}\, R^{23} \,\Fii[231]{} \,R^{13}\,\Fi[213]{}\, R^{12}
\end{equation}
and the property
\begin{equation}
  \label{eq112}
  (\ep\tp\id)(R) = (\id\tp\ep)(R) = \e.
\end{equation}
Eq.\ (\ref{eq18}) implies, that for any pair $\pi_U,\pi_V$ the two
representations $(\pi_{U\tp 
  V},U\tp V)$ and $(\pi_{V\tp U},V\tp U)$ are equivalent with
intertwiner $B^{UV}: = \tau^{12}\circ (\pi_U\tp\pi_V)(R)$,
where $\tau^{12}$ denotes the permutation of 
tensor factors in $U\tp V$. Eqs.~\eqref{eq19},\eqref{eq110} imply the
commutativity of two hexagon diagrams obtained by taking
$\pi_U\tp\pi_V\tp\pi_W$ on both sides.

$\G$ being a quasitriangular quasi--Hopf algebra
implies that $\rep
\G$ is a rigid monoidal category with braiding, where the
associativity and 
commutativity constraints for the tensor product functor $\tp : \rep
\G \times \rep 
\G \pfeil \rep \G$ are given by the natural families $\phi_{UVW}$ and
$\tau^{12}\circ R_{UV}$ and the (left) duality is defined with the help of the
antipode $S$ and the elements $\alpha,\beta$, see
(\ref{graf1}--\ref{graf3}) below.  

Together with a quasi--Hopf algebra $\G\equiv
(\G,\cop,\ep,\phi,S,\alpha,\beta)$ we also have $\G_{op}, \G^{cop}$
and $\G_{op}^{cop}$ as quasi--Hopf algebras, where ``$op$'' means
opposite multiplication and ``$cop$'' means opposite
comultiplication. The quasi--Hopf structures are obtained by putting
$\phi_{op}:=\phi^{-1}$, $\phi^{cop}:=\Fii[321]{}$,
$\phi_{op}^{cop}:=\Fi[321]{}$, $S_{op}=S^{cop} =
(S_{op}^{cop})^{-1}:=\Si$, $\alpha_{op}:= \Si(\beta)$, $\beta_{op}
:=\Si(\alpha)$, $\alpha^{cop} := \Si(\alpha)$, $\beta^{cop} :=
\Si(\beta)$, $\alpha_{op}^{cop} :=\beta$ and $\beta_{op}^{cop} :=
\alpha$. Also if $R \in \G\tp\G$ is quasitriangular in $\G$, then
$R^{-1}$ is quasitriangular in $\G_{op}$, $R^{21}$ is quasitriangular
in $\G^{cop}$ and $(R^{-1})^{21}$ is quasitriangular in $\G_{op}^{cop}$

Next we recall that the definition of a quasitriangular quasi-Hopf algebra is
`twist covariant' in the following sense:
An element $F\in \G\tp\G$ which is invertible and satisfies
$(\ep\tp\id)(F) = (\id\tp\ep)(F) = \e$, induces a so--called {\it
twist transformation}
\begin{align}
  \label{eq113}
  \cop_F(a) : &= F \cop(a) F, \\
   \label{eq114}
  \phi_F : &= (\e\tp F)\, (\id\tp\cop)(F) \, \phi \,
  (\cop\tp\id)(F^{-1})\, (F^{-1}\tp\e)
\end{align}
It has been noticed by Drinfel'd [Dr2] that
$(\G,\cop_F,\ep,\phi_F)$ is again a quasi--bialgebra.
Setting 
\begin{equation*}
  \alpha_F := S(h^i)\alpha k^i, \quad 
  \beta_F := f^i \beta S(g^i),
\end{equation*}
where $h^i\tp k^i = F^{-1}$ and $f^i\tp g^i  = F$,
$(\G,\cop_F,\ep,\phi_F, S, \alpha_F, \beta_F)$ is also a quasi-Hopf
algebra. Moreover, if $R$ is quasitriangular with respect to $(\cop,
\phi)$, then
\begin{equation}
  \label{eq115}
  R_F : = F^{21} R F^{-1}
\end{equation}
is quasitriangular w.r.t. $(\cop_F,\phi_F)$.
This means that a twist preserves the class of quasitriangular
quasi-Hopf algebras [Dr2]. 

For Hopf algebras, one knows, that the antipode is an anti coalgebra
morphism, i.e. $\cop(a) = (S\tp S)\big(\cop^{op}(\Si(a))\big)$. For
quasi-Hopf algebras this is true only up to a twist: Following
Drinfeld we define the elements $\gamma, \delta \in \G\tp\G$ by
setting\footnote{suppressing summation symbols} 
\begin{align}
  \label{eq116}
   \gamma & := (S(U^i)\tp S(T^i))\cdot (\alpha\tp\alpha)\cdot
   (V^i\tp W^i) 
   \\
   \label{eq117}
   \delta & := (K^j\tp L^j)\cdot (\beta\tp\beta)\cdot (S(N^j)\tp S(M^j))
\end{align}
where
\begin{align}
  T^i\tp U^i \tp V^i \tp W^i & =
               (\e\tp\phi^{-1})\cdot (\id\tp\id\tp\cop)(\phi), \\
   K^j \tp L^j \tp M^j \tp N^j & =
              (\cop\tp\id\tp\id)(\phi)\cdot (\phi^{-1}\tp\e).
\end{align}
With these definitions Drinfel'd has shown in [Dr2], that $f\in \G\tp\G$ given by
\begin{equation}
  \label{eq118}
  f:= (S\tp S)(\cop^{op}(P^i)) \cdot \gamma \cdot \cop(Q^i\beta R^i).
\end{equation}
defines a twist with inverse given by 
\begin{equation}
  \label{eq119}
  f^{-1} = \cop(S(P^j)\alpha Q^j) \cdot \delta \cdot (S\tp
  S)(\cop^{op}(R^i)), 
\end{equation}
such that for all $a\in\G$
\begin{equation}
  \label{eq120} f\cop(a)f^{-1} = (S\tp S)\big(\cop^{op}(\Si(a))\big).
\end{equation}
The elements $\gamma, \delta$ and the twist $f$ fulfill the relation
\begin{equation}
  \label{eq119b}
  f\, \cop(\alpha) = \gamma,\quad \cop(\beta) \, f^{-1} = \delta
\end{equation}
Furthermore, the corresponding twisted reassociator (\ref{eq114}) is
given by
\begin{equation}
  \label{eq121}
  \phi_f = (S\tp S\tp S)(\phi^{321}).
\end{equation}
Setting $h : = (\Si\tp\Si)(f^{21})$, the above relations imply
\begin{align}
  \label{eq122}
  h\cop(a) h^{-1} &= (\Si\tp\Si)\big(\cop^{op}(S(a))\big)  \\
  \label{eq123}
   \phi_h &= (\Si\tp\Si\tp\Si)(\phi^{321}).
\end{align}
The importance of the twist $f$ for the representation theory of $\G$ is
the existence of an intertwiner $U\tp V  \pfeil
 ({}^*V\tp {}^*U)^*$ given by $\tau^{12}\circ (\pi_U\tp\pi_V)(f)$.\\

Finally we introduce $\dG$ as the dual space of $\G$ with its natural
coassociative coalgebra structure $(\hat{\cop},\hat{\ep})$ given by
$\lpa \hat{\cop}(\vi) \mid a \tp b\rpa:=\lpa \vi \mid ab\rpa$ and
$\hat{\ep}(\vi) : = \lpa \vi \mid \eG\rpa$, where $\vi \in \dG,\, a,b
\in \G$ and where $\lpa \cdot \mid \cdot \rpa : \dG \tp \G \pfeil
\mathbb{C}$ denotes the dual pairing. On $\dG$ we have the natural
left and right $\G$--actions 
\begin{equation*}
  a \arr \vi : = \vi_{(1)} \lpa \vi_{(2)}\mid a \rpa, \quad 
  \vi \arl a := \vi_{(2)} \lpa \vi_{(1)}\mid a\rpa,
\end{equation*}
where $a \in \G, \, \vi \in \dG$. By transposing the coproduct on $\G$
we also get a multiplication $\dG\tp\dG \pfeil \dG$, which however is
no longer associative
\begin{equation*}
  \lpa \vi\psi\mid a\rpa : = \lpa \vi\tp \psi\mid \cop(a)\rpa, \quad
  \lpa \edG \mid a \rpa : = \ep(a).
\end{equation*}
Yet, we have the identities $\edG \vi = \vi \edG = \vi$,
$\hat{\cop}(\vi\psi) = \hat{\cop}(\vi) \hat{\cop}(\psi)$, $a \arr(\vi\psi)=
(a_{(1)}\arr \vi)$ $(a_{(2)} \arr \psi)$ and $(\vi\psi)\arl a = (\vi \arl
a_{(1)}) (\psi \arl a_{(2)})$ for all $\vi, \psi \in \dG$ and $a \in
\G$. We also introduce $\hat{S}: \dG \pfeil \dG$ as the coalgebra
anti-mophism dual to $S$, i.e.\ $\lpa \hat{S}(\vi)\mid a \rpa := \lpa \vi
\mid S(a)\rpa$.

%

\unitlength0.5cm
\newcommand{\Rmatrix}{
\begin{picture}(2,2)(0.25,0) 
\put(0,0){\line(1,1){2}}
\put(0,2){\line(1,-1){0.8}}
\put(1.2,0.8){\line(1,-1){0.8}}
\end{picture} }
\newcommand{\Dmatrix}{
\begin{picture}(2,2)(0.25,0) 
\put(0,0){\line(1,1){2}}
\multiput(0.1,1.9)(0.2,-0.2){10}{\circle*{0.15}}
\end{picture} }
\newcommand{\Dmatrixgross}{
\begin{picture}(3,3)(0.25,0) 
\put(0,0){\line(1,1){3}}
\multiput(0.1,2.9)(0.2,-0.2){15}{\circle*{0.15}}
\end{picture} }

\newcommand{\iRmatrix}{
\begin{picture}(2,2)(0.25,0) 
\put(0,2){\line(1,-1){2}}
\put(0,0){\line(1,1){0.8}}
\put(1.2,1.2){\line(1,1){0.8}}

\end{picture} }
\newcommand{\iDmatrix}{
\begin{picture}(2,2)(0.25,0)
\multiput(0.1,0.1)(0.2,0.2){10}{\circle*{0.15}}
\put(0,2){\line(1,-1){2}}
\end{picture} }

\newcommand{\boxza}{
  \begin{picture}(10,3)(0.25,0)
  \multiput(1,0)(8,0){2}{\line(0,1){1}}
  \multiput(1,2)(8,0){2}{\line(0,1){1}}
  \multiput(2,0)(0.5,0){13}{\circle*{0.15}}
  \put(0,1){\framebox(10,1){}}
  \put(0.5,0.25){(} \put(0.5,2.5){(}
  \put(8.5,0.25){)} \put(9.25,2.5){)}
  \end{picture}
}

\newcommand{\boxzb}{
  \begin{picture}(6,3)(0.25,0)
  \multiput(1,0)(4,0){2}{\line(0,1){1}}
  \multiput(1,2)(4,0){2}{\line(0,1){1}}
  \multiput(2,3)(0.5,0){5}{\circle*{0.15}}
  \put(0,1){\framebox(6,1){}}
  \put(0.5,0.25){(} \put(0.5,2.5){(}
  \put(4.5,2.5){)} \put(5.25,0.25){)}
  \end{picture}
}
 
\newcommand{\boxd}{
  \begin{picture}(6,3)(0.25,0)
  \multiput(1,2)(2,0){3}{\line(0,1){1}}
  \multiput(1,0)(2,0){3}{\line(0,1){1}}
  \put(0,1){\framebox(6,1){}}
  \end{picture}
}

\newcommand{\boxdi}{
  \begin{picture}(6,3)(0.25,0)
  \multiput(1,2)(2,0){3}{\line(0,1){1}}
  \multiput(1,0)(2,0){3}{\line(0,1){1}}
  \put(0,1){\framebox(6,1){}}
  \put(0.5,0.25){(}
  \put(2.5,2.5){(}
  \put(5.25,2.5){)}
  \put(3.25,0.25){)}
  \end{picture}
}

\newcommand{\boxdii}{
  \begin{picture}(6,3)(0.25,0)
  \multiput(1,2)(2,0){3}{\line(0,1){1}}
  \multiput(1,0)(2,0){3}{\line(0,1){1}}
  \put(0,1){\framebox(6,1){}}
  \put(0.5,2.5){(}
  \put(2.5,0.25){(}
  \put(5.25,0.25){)}
  \put(3.25,2.5){)}
  \end{picture}
}

\newcommand{\boxv}{
  \begin{picture}(8,3)(0.25,0)
  \multiput(1,2)(2,0){4}{\line(0,1){1}}
  \multiput(1,0)(2,0){4}{\line(0,1){1}}
  \put(0,1){\framebox(8,1){}}
  \end{picture}
}
\newcommand{\boxvi}{
  \begin{picture}(8,3)(0.25,0)
  \multiput(1,2)(2,0){4}{\line(0,1){1}}
  \multiput(1,0)(2,0){4}{\line(0,1){1}}
  \multiput(0.5,2.5)(4,0){2}{(}
  \multiput(3.25,2.5)(4,0){2}{)}
  \put(2.25,0.25){[(} \put(5.25,0.25){)} \put(7.25,0.25){]}
  \put(0,1){\framebox(8,1){}}
  \end{picture}
}
\newcommand{\boxvii}{
  \begin{picture}(8,3)(0.25,0)
  \multiput(1,2)(2,0){4}{\line(0,1){1}}
  \multiput(1,0)(2,0){4}{\line(0,1){1}}
  \multiput(0.5,0.25)(4,0){2}{(}
  \multiput(3.25,0.25)(4,0){2}{)}
  \put(2.25,2.5){[(} \put(5.25,2.5){)} \put(7.25,2.5){]}
  \put(0,1){\framebox(8,1){}}
  \end{picture}
}
\newcommand{\boxviii}{
  \begin{picture}(8,3)(0.25,0)
  \multiput(1,2)(2,0){4}{\line(0,1){1}}
  \multiput(1,0)(2,0){4}{\line(0,1){1}}
  \multiput(0.5,0.25)(4,0){2}{(}
  \multiput(3.25,0.25)(4,0){2}{)}
  \put(0.25,2.5){[(} \put(3.25,2.5){)} \put(5.25,2.5){]}
  \put(0,1){\framebox(8,1){}}
  \end{picture}
}
\newcommand{\boxviv}{
  \begin{picture}(8,3)(0.25,0)
  \multiput(1,2)(2,0){4}{\line(0,1){1}}
  \multiput(1,0)(2,0){4}{\line(0,1){1}}
  \multiput(2.5,0.25)(0,2.25){2}{(}
  \multiput(5.25,0.25)(0,2.25){2}{)}
  \put(0.5,2.5){[} \put(5.5,2.5){]} \put(7.5,0.25){]}\put(2.25,0.25){[}
  \put(0,1){\framebox(8,1){}}
  \end{picture}
}
\newcommand{\boxvv}{
  \begin{picture}(8,3)(0.25,0)
  \multiput(1,2)(2,0){4}{\line(0,1){1}}
  \multiput(1,0)(2,0){4}{\line(0,1){1}}
  \multiput(0.5,0.25)(4,0){2}{(}
  \multiput(3.25,0.25)(4,0){2}{)}
  \put(4.5,2.5){(} \put(7.25,2.5){)]} \put(2.5,2.5){[}
  \put(0,1){\framebox(8,1){}}
  \end{picture}
}
\newcommand{\boxvvi}{
  \begin{picture}(8,3)(0.25,0)
  \multiput(1,2)(2,0){4}{\line(0,1){1}}
  \multiput(1,0)(2,0){4}{\line(0,1){1}}
  \multiput(0.5,0.25)(4,0){2}{(}
  \multiput(3.25,0.25)(4,0){2}{)}
  \put(2.5,2.5){(} \put(5.25,2.5){)]} \put(0.5,2.5){[}
  \put(0,1){\framebox(8,1){}}
  \end{picture}
}
\newcommand{\boxvvii}{
  \begin{picture}(8,3)(0.25,0)
  \multiput(1,2)(2,0){4}{\line(0,1){1}}
  \multiput(1,0)(2,0){4}{\line(0,1){1}}
  \multiput(0.5,2.5)(4,0){2}{(}
  \multiput(3.25,2.5)(4,0){2}{)}
  \put(2.5,0.25){(} \put(5.25,0.25){)]} \put(0.5,0.25){[}
  \put(0,1){\framebox(8,1){}}
  \end{picture}
}
\newcommand{\boxvviii}{
  \begin{picture}(8,3)(0.25,0)
  \multiput(1,2)(2,0){4}{\line(0,1){1}}
  \multiput(1,0)(2,0){4}{\line(0,1){1}}
  \multiput(2.5,0.25)(0,2.25){2}{(}
  \multiput(5.25,0.25)(0,2.25){2}{)}
  \put(2.25,2.5){[} \put(7.5,2.5){]} \put(5.5,0.25){]}\put(0.5,0.25){[}
  \put(0,1){\framebox(8,1){}}
  \end{picture}
}
\newcommand{\boxvviv}{
  \begin{picture}(8,3)(0.25,0)
  \multiput(1,2)(2,0){4}{\line(0,1){1}}
  \multiput(1,0)(2,0){4}{\line(0,1){1}}
  \multiput(2.5,0.25)(2,2.25){2}{(}
  \multiput(5.25,0.25)(2,2.25){2}{)}
  \put(2.5,2.5){[} \put(7.5,2.5){]} \put(5.5,0.25){]}\put(0.5,0.25){[}
  \put(0,1){\framebox(8,1){}}
  \end{picture}
}

\newcommand{\boxva}{
  \begin{picture}(10,3)(0.25,0)
  \put(0,1){\framebox(10,1){}}
  \put(1,0){\line(0,1){1}}
  \put(1,2){\line(0,1){1}}
  \multiput(5,0)(2,0){3}{\line(0,1){1}}
  \multiput(5,2)(2,0){3}{\line(0,1){1}}
    \multiput(2,0)(0.5,0){5}{\circle*{0.15}}
  \put(0.5,0.25){(}\put(0.25,2.5){\big[}
    \put(0.5,2.5){(}
  \put(5.25,0.25){)} \put(6.5,0.25){[}\put(9.25,0.25){]}
  \put(7.25,2.5){\big]} \put(5.25,2.5){)}
  \end{picture}
} 
\newcommand{\boxvb}{
  \begin{picture}(10,3)(0.25,0)
  \put(0,1){\framebox(10,1){}}
  \put(1,0){\line(0,1){1}}
  \put(1,2){\line(0,1){1}}
  \multiput(5,0)(2,0){3}{\line(0,1){1}}
  \multiput(5,2)(2,0){3}{\line(0,1){1}}
    \multiput(2,3)(0.5,0){5}{\circle*{0.15}}
  \put(0.5,0.25){(}\put(0.25,2.5){\big[}
    \put(0.5,2.5){(}
  \put(5.25,0.25){)} \put(6.5,0.25){(}\put(9.25,0.25){)}
  \put(7.25,2.5){)\big]} \put(4.25,2.5){)(}
  \end{picture}
} 
\newcommand{\boxs}{
  \begin{picture}(12,3)(0.25,0)
  \multiput(0.5,2.25)(4,0){3}{(}
  \multiput(3.25,2.25)(4,0){3}{)}
  \multiput(0.5,0.25)(4,0){3}{(}
  \multiput(3.25,0.25)(4,0){3}{)}
  \multiput(1,2)(2,0){6}{\line(0,1){1}}
  \multiput(1,0)(2,0){6}{\line(0,1){1}}
  \put(0,1){\framebox(12,1){}}
  \end{picture}
}
\newcommand{\boxda}{
  \begin{picture}(8,3)(0.25,0)
  \put(0,1){\framebox(8,1){}}
  \multiput(5,0)(2,0){2}{\line(0,1){1}}
  \multiput(5,2)(2,0){2}{\line(0,1){1}}
  \put(1,0){\line(0,1){1}}
  \put(1,2){\line(0,1){1}}
  \put(0.5,0.25){(}
  \put(4.5,2.5){(}
  \put(5.25,0.25){)}
  \put(7.25,2.5){)}
  \multiput(2,3)(0.5,0){5}{\circle*{0.15}}
  \end{picture}
}
\newcommand{\boxdb}{
  \begin{picture}(8,3)(0.25,0)
  \put(0,1){\framebox(8,1){}}
  \multiput(5,0)(2,0){2}{\line(0,1){1}}
  \multiput(5,2)(2,0){2}{\line(0,1){1}}
  \put(1,0){\line(0,1){1}}
  \put(1,2){\line(0,1){1}}
  \put(0.5,2.5){(}
  \put(4.5,0.25){(}
  \put(5.25,2.5){)}
  \put(7.25,0.25){)}
  \multiput(2,0)(0.5,0){5}{\circle*{0.15}}
  \end{picture}
}

\newcommand{\boxdc}{
  \begin{picture}(8,3)(0.25,0)
  \put(0,1){\framebox(8,1){}}
  \multiput(1,0)(2,0){2}{\line(0,1){1}}
  \multiput(1,2)(2,0){2}{\line(0,1){1}}
  \put(7,0){\line(0,1){1}}
  \put(7,2){\line(0,1){1}}
  \put(0.5,0.25){(}
  \put(2.5,2.5){(}
  \put(3.25,0.25){)}
  \put(7.25,2.5){)}
  \multiput(4,0)(0.5,0){5}{\circle*{0.15}}
  \end{picture}
}
\newcommand{\boxdd}{
  \begin{picture}(8,3)(0.25,0)
  \put(0,1){\framebox(8,1){}}
  \multiput(1,0)(2,0){2}{\line(0,1){1}}
  \multiput(1,2)(2,0){2}{\line(0,1){1}}
  \put(7,0){\line(0,1){1}}
  \put(7,2){\line(0,1){1}}
  \put(2.5,0.25){(}
  \put(0.5,2.5){(}
  \put(3.25,2.5){)}
  \put(7.25,0.25){)}
  \multiput(4,3)(0.5,0){5}{\circle*{0.15}}
  \end{picture}
}

\newcommand{\okreis}{
  \begin{picture}(2,2)(0.25,0)
  \put(1,0){\oval(2,2)[t]}
  \end{picture}
}
\newcommand{\gokreis}{
  \begin{picture}(6,4)(0.25,0)
  \put(3,0){\oval(6,4)[t]}
  \end{picture}
}
\newcommand{\ukreis}{
  \begin{picture}(2,2)(0.25,0)
  \put(1,0){\oval(2,2)[b]}
  \end{picture}
}
\newcommand{\gukreis}{
  \begin{picture}(6,4)(0.25,0)
  \put(3,0){\oval(6,4)[b]}
  \end{picture}
}

\newcommand{\bilda}{ 
  \begin{picture}(15,6)(0.25,0)
   \put(0,2){\framebox(6,2){\it g}}
   \multiput(1,0)(4,0){2}{\line(0,1){2}}
    \multiput(1,4)(4,0){2}{\line(0,1){2}}
    \multiput(2,1)(0.5,0){5}{\circle*{0.15}}
    \multiput(2,5)(0.5,0){5}{\circle*{0.15}}
    \put(8,5){target}
    \put(8,0){source}
  \end{picture}
}
\newcommand{\bildbprime}{ 
  \begin{picture}(6,4)(0.25,0)
   \put(0,1){\framebox(6,2){$g'$}}
   \multiput(1,0)(4,0){2}{\line(0,1){1}}
    \multiput(1,3)(4,0){2}{\line(0,1){1}}
    \multiput(2,0)(0.5,0){5}{\circle*{0.15}}
    \multiput(2,4)(0.5,0){5}{\circle*{0.15}}
  \end{picture}
}
\newcommand{\bildb}{ 
  \begin{picture}(6,4)(0.25,0)
   \put(0,1){\framebox(6,2){\it g}}
   \multiput(1,0)(4,0){2}{\line(0,1){1}}
    \multiput(1,3)(4,0){2}{\line(0,1){1}}
    \multiput(2,0)(0.5,0){5}{\circle*{0.15}}
    \multiput(2,4)(0.5,0){5}{\circle*{0.15}}
  \end{picture}
}

\newcommand{\ogabel}{
  \begin{picture}(1,1)(0.25,0)
   \put(0,0){\line(0,1){1}}
   \put(0,0){\line(1,0){1}}
   \put(1,0){\line(0,1){1}}
  \end{picture}
}
\newcommand{\ugabel}{
  \begin{picture}(1,1)(0.25,0)
   \put(1,0){\line(0,1){1}}
   \put(0,1){\line(1,0){1}}
   \put(0,0){\line(0,0){1}}
  \end{picture}
}
\newcommand{\liniez}{
  \begin{picture}(2,2)(0.25,0)
  \line(0,1){2}
  \end{picture}
}

\newcommand{\bildc}{
\begin{picture}(8,10)(0.25,0)
\put(5,1){\ukreis}
\put(0,1){\boxvii}
\put(1,4){\liniez}
\put(3,4){\Rmatrix}
\put(7,4){\liniez}
\put(0,6){\boxvi}
\put(1,9){\okreis}
\multiput(5,9)(2,0){2}{\line(0,1){1}}
\multiput(1,0)(2,0){2}{\line(0,1){1}}
\end{picture}
}

\newcommand{\bildcD}{
\begin{picture}(8,10)(0.25,0)
\put(5,1){\ukreis}
\put(0,1){\boxvii}
\put(1,4){\liniez}
\put(3,4){\Dmatrix}
\put(7,4){\liniez}
\put(0,6){\boxvi}
\put(1,9){\okreis}
\multiput(5,9)(2,0){2}{\line(0,1){1}}
\multiput(1,0)(2,0){2}{\line(0,1){1}}
\end{picture}
}

\newcommand{\bildd}{
\begin{picture}(8,10)(0.25,0)
\put(1,1){\ukreis}
\put(0,1){\boxvii}
\put(1,4){\liniez}
\put(3,4){\Rmatrix}
\put(7,4){\liniez}
\put(0,6){\boxvi}
\put(5,9){\okreis}
\multiput(1,9)(2,0){2}{\line(0,1){1}}
\multiput(5,0)(2,0){2}{\line(0,1){1}}
\end{picture}
}
\newcommand{\bilddD}{
\begin{picture}(8,10)(0.25,0)
\put(1,1){\ukreis}
\put(0,1){\boxvii}
\put(1,4){\liniez}
\put(3,4){\Dmatrix}
\put(7,4){\liniez}
\put(0,6){\boxvi}
\put(5,9){\okreis}
\multiput(1,9)(2,0){2}{\line(0,1){1}}
\multiput(5,0)(2,0){2}{\line(0,1){1}}
\end{picture}
}
\newcommand{\bilde}{
\begin{picture}(8,5)(0.25,0)
\put(0,0){\boxvii}
\put(3,3){\okreis}
\put(1,3){\gokreis}
\end{picture}
}
\newcommand{\bildf}{
\begin{picture}(8,5)(0.25,0)
\put(0,0){\boxvi}
\put(3,0){\ukreis}
\put(1,0){\gukreis}
\end{picture}
}

%
%
%
\subsection{Graphical calculus}
In the following it will be useful to have a graphical notation for the
identities and definitions given so far. The graphical calculus
introduced below has been developed and used in many papers, 
e.g. [RT,AC,T], mainly in the setting of ribbon--Hopf algebras.
Formally speaking, it consists of a functor from the braided monoidal
category $\rep \G$ into a category of colored graphs.
 For an introduction into the category terminology see [K], [T].
We will use the  graphical notation to have a pictorial way to understand - and
deduce - certain relations and identities between morphisms (intertwiners) in
$\rep \G$, which -
written out algebraically - would look very complicated. By 
morphisms in $\rep_\G$ we mean elements $t\in \makebox{Hom}_\G(U,V)$,
i.e. linear maps $t:\, U \rightarrow V$ satisfying $t\, \pi_U (a) = \pi_V
(a)\, t, \, \forall a \in \G$. As discussed in Section 2.1, the
$n$--fold tensor product of  $\G$--modules is again a $\G$--module (where one
has to take care of the bracketing of the tensor factors). A morphism
$t$ from an $n$-fold to an $m$-fold tensor product of $\G$--modules is
represented by a graph consisting of a ``coupon'' with $n$ lower legs
and $m$ upper legs ``coloured'' with the 
source and target modules respectively. The upper and lower legs are
always equipped with a definite bracketing corresponding to the
bracketing defining the associated tensor module. For example the picture

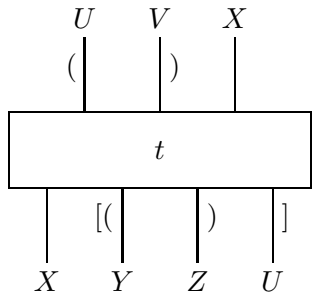
\begin{equation*}
\begin{picture}(8,8)
 \put(1,0.5){\makebox(0,0){$X$}}
 \put(3,0.5){\makebox(0,0){$Y$}}
 \put(5,0.5){\makebox(0,0){$Z$}}
 \put(7,0.5){\makebox(0,0){$U$}}
 \put(2.25,2){[(} \put(5.25,2){)} \put(7.25,2){]}
 \multiput(1,1)(2,0){4}{\line(0,1){2}}
 \put(0,3){\framebox(8,2){$t$}}
 \multiput(2,5)(2,0){3}{\line(0,1){2}}
 \put(2,7.5){\makebox(0,0){$U$}}
 \put(4,7.5){\makebox(0,0){$V$}}
 \put(6,7.5){\makebox(0,0){$X$}}
 \put(1.5,6){(} \put(4.25,6){)}
\end{picture}
\end{equation*}
corresponds to the morphism
\begin{equation*}
  t: \, X\tp\big[ (Y\tp Z) \tp U \big] \pfeil
    (U \tp V) \tp X
\end{equation*}
The tensor product of two morphisms corresponds to the juxtaposition of
diagrams and the composition of morphisms is depicted by gluing the
corresponding graphs together. Here one has to take care that the
gluing $t \circ k$ is only admissible if source($t$) = target($k$),
which in particular implies that the bracketing conventions of the
associated tensor factors have to coincide. We also use the convention
that the lower legs always represent the source, i.e. the graph
$t\circ k$ is obtained by gluing $t$ on top of $k$. 

Following the conventions of [AC] we now give a list of some
special morphisms depicted by the following graphs:
\begin{equation}
  \label{graf1}
  \begin{picture}(9,4)
  \put(0,0){\line(0,1){3}}
  \put(0,3.5){\makebox(0,0){V}}
  \put(1,1.75){$:= \,\, \id_V$,}
   \end{picture}
\end{equation}
\begin{equation}
      \label{graf2}
       \begin{picture}(13,5)
   \put(0,2){\ukreis}
   \put(0,2.5){\makebox(0,0){$V$}}
   \put(2,2.5){\makebox(0,0){${}^*V$}}
   \put(3,1.75){:=}
   \put(5,1.25){\framebox(4,1.5){$b_V$}}
   \put(7,0){\line(0,1){1.25}}
   \multiput(6,2.75)(2,0){2}{\line(0,1){1.25}}
   \put(6,4.5){\makebox(0,0){$V$}}
   \put(8,4.5){\makebox(0,0){${}^*V$}}
   \put(7,-0.5){\makebox(0,0){$\mathbb{C}$}}
   \put(9.5,1.75){,}
   \end{picture}
    \begin{picture}(11,5)
   \put(0,2){\okreis}
   \put(2,1.5){\makebox(0,0){$V$}}
   \put(0,1.5){\makebox(0,0){${}^*V$}}
   \put(3,1.75){:=}
   \put(5,1.25){\framebox(4,1.5){$a_V$}}
   \put(7,2.75){\line(0,1){1.25}}
   \multiput(6,0)(2,0){2}{\line(0,1){1.25}}
   \put(8,-0.5){\makebox(0,0){$V$}}
   \put(6,-0.5){\makebox(0,0){${}^*V$}}
   \put(7,4.5){\makebox(0,0){$\mathbb{C}$}}
   \end{picture}   
\end{equation}
%
%
\begin{equation}
  \label{graf3}
  \begin{picture}(13,4)
  \put(0,1){\Rmatrix}
  \put(0,0.5){\makebox(0,0){$V$}}
  \put(2,0.5){\makebox(0,0){$W$}}
  \put(0,3.5){\makebox(0,0){$W$}}
  \put(2,3.5){\makebox(0,0){$V$}}
  \put(3,1.75){$:=\,\,B_{VW}$,}
  \end{picture}
 \begin{picture}(11,4)
  \put(0,1){\iRmatrix}
  \put(0,0.5){\makebox(0,0){$W$}}
  \put(2,0.5){\makebox(0,0){$V$}}
  \put(0,3.5){\makebox(0,0){$V$}}
  \put(2,3.5){\makebox(0,0){$W$}}
  \put(3,1.75){$:=\,\,B_{VW}^{-1}$}
  \end{picture}
\end{equation}
where $\mathbb{C}$ stands for the one dimensional representation given by the
counit and where
\begin{align*}
  & b_V :\, \mathbb{C} \pfeil V\tp{}^*V, \,\, 
     1 \longmapsto \sum_i \beta_V v_i \tp v^i \\
  & a_V :\, {}^*V\tp V \pfeil \mathbb{C},\,\, 
   \hat{v} \tp w \longmapsto \lpa \hat{v} \mid \alpha_V w\rpa \\
 &  B_{VW} = \tau_{VW} \circ R_{VW}, \,\quad
  B_{VW}^{-1} = R_{VW}^{-1} \circ \tau_{WV} .
\end{align*}
Here $\{v_i\}$ is a basis of $V$ with dual basis $\{v^i\}$ and
$\tau _{VW}$ denotes the permutation of tensor factors in $V\tp W$. We
also use the shortcut notation $\alpha_V \equiv \pi_V(\alpha)$,
$R_{VW} \equiv (\pi_V\tp\pi_W)(R)$, etc.
The properties of $\G$ being a quasitriangular quasi-Hopf algebra
ensure, that the above defined maps are in fact intertwiners
(morphisms of $\G$-modules). Note that within higher tensor products
the graphs (\ref{graf2}) and (\ref{graf3}) are only admissible if
their legs are ``bracketed together''. In order to change the bracket
convention one has to use  rebracketing morphisms. These are given as
products of the basic elements
\begin{equation}
  \label{graf4}
  \begin{picture}(26,5)
  \put(1,0.5){\makebox(0,0){$U$}}
  \put(3,0.5){\makebox(0,0){$V$}}
  \put(5,0.5){\makebox(0,0){$W$}}
  \put(0,1){\boxdi}
  \put(7,2.25){=} 
  \put(8,2.25){$\phi_{UVW}$,}  
  \put(15,0.5){\makebox(0,0){$U$}}
  \put(17,0.5){\makebox(0,0){$V$}}
  \put(19,0.5){\makebox(0,0){$W$}}
  \put(14,1){\boxdii}
    \put(21,2.25){=} 
  \put(22,2.25){$\phi_{UVW}^{-1}$}              
\end{picture}
\end{equation}
where each of the three individual legs in (\ref{graf4}) may again
represent a tensor product of $\G$--modules. In this way we adopt the
convention that any empty coupon with the same number of upper and
lower legs - where the colouring only differs by the bracket convention
- always represents the associated unique rebracketing morphism in
$\rep\G$ given in terms of suitable products of $\phi$'s. We have
already remarked that the uniqueness of this rebracketing morphism
(i.e. the independence of the chosen sequence of intermediate
rebracketings) is guaranteed by McLanes coherence theorem and the
``pentagon axiom'' (\ref{eq12}). This is why it is often not even
necessary to spell out one of the possible formulas for such an
intertwiner. Explicitly, the pentagon identity (\ref{eq12}) may be
expressed  as 
\begin{equation}
  \label{graf5d}
  \begin{picture}(20,9)
  \put(0,0){\boxviii}
  \put(0,3){\boxd}
  \put(7,3){\line(0,1){3}}
  \put(0,6){\boxvviii}
  \put(10,4.5){\makebox(0,0){=}}
  \put(12,1){\boxvv}
  \put(13,4){\line(0,1){3}}
  \put(14,4){\boxd}
  \put(14.25,6.5){[(} \put(17.25,6.5){)} \put(19.25,6.5){]}
  \end{picture}
\end{equation}
which is the graphical notation for 
\begin{equation*}
    (\id\tp\cop\tp\id)(\phi)\,  (\phi \tp\e)\,
    (\cop\tp\id\tp\id)(\phi^{-1}) =
   (\e\tp\phi^{-1}) \, (\id\tp\id\tp\cop)(\phi) 
\end{equation*}
In the same philosophy one may rewrite a simple rebracketing as a
product of more complicated ones, as long as the overall source and
target brackets coincide, for example 
\begin{equation}
  \label{graf4a}
  \begin{picture}(20,7)
  \multiput(0,0)(2,0){4}{\line(0,1){2}}
  \put(0,2){\line(0,1){3}}
  \put(1,2){\boxdi}
  \multiput(1.25,2.25)(0,2.25){2}{[}
   \multiput(6.75,2.25)(0,2.25){2}{]}
  \multiput(0,5)(2,0){4}{\line(0,1){2}}
  \put(9.5,3.5){\makebox(0,0){=}}
  \put(11,0){\boxvi}
  \multiput(12,3)(2,0){4}{\line(0,1){1}}
  \put(11,4){\boxvv}
  \end{picture}
\end{equation}
As done in the above pictures we will
frequently not specify the modules sitting at the source and 
target legs. 
Also note that by Eqs.\. (\ref{eq14}) and (\ref{eq15}) the rebracketing
of the (invisible) ``white'' leg corresponding to the trivial
$\G$--module $\mathbb{C}$ is always given by the trivial
identification.

If $\G$ is finite dimensional, it may itself be
viewed as a $\G$-module under left multiplication and algebraic
identities may directly be 
translated into identities of the corresponding graphs and vice versa.
So e.g. Eq.\ (\ref{eq17}) is equivalent to
\begin{equation}
  \label{graf5}
  \begin{picture}(20,7)
  \put(0,2){\boxdi}
  \put(1,5){\line(0,1){1}}
  \put(3,5){\okreis}
  \put(1,2){\ukreis}
  \put(5,1){\line(0,1){1}}
  \put(7,3.25){=} 
  \put(8.5,1){\line(0,1){5}}
  \put(9.5,3.25){=} 
  \put(11,2){\boxdii}
  \put(12,5){\okreis}
  \put(16,5){\line(0,1){1}}
  \put(14,2){\ukreis}
  \put(12,1){\line(0,1){1}}
  \end{picture}
\end{equation}
and Eqs.\ (\ref{eq14}) and (\ref{eq15}) together with \eqref{eq16} imply
\begin{equation}
  \label{graf5b}
\begin{picture}(20,4)
 \multiput(1,0)(2,0){4}{\line(0,1){1.5}}
  \multiput(1,2.5)(2,0){4}{\line(0,1){1}}
   \multiput(1,3.5)(6,0){2}{\line(0,1){1}}
   \put(0.5,0.75){[}
   \put(2.5,0.75){(}
   \put(5.25,0.75){)]}
   \put(2.25,3){[(}
   \put(7.25,3){]}
   \put(5.25,3){)}
   \put(0,1.5){\framebox(8,1){}}
   \put(3,3.5){\okreis}
   \put(10,2){=}
    \multiput(12,1)(2,0){4}{\line(0,1){1}}
   \multiput(12,2)(6,0){2}{\line(0,1){1}}
    \put(14,2){\okreis}
     \put(11.5,1.25){[}
   \put(13.5,1.25){(}
   \put(16.25,1.25){)]} 
\end{picture}
\end{equation}
and
\begin{equation}
  \label{graf5c}
\begin{picture}(20,4)
  \multiput(1,0)(2,0){4}{\line(0,1){1.5}}
  \multiput(1,2.5)(2,0){4}{\line(0,1){1}}
   \multiput(1,3.5)(2,0){2}{\line(0,1){1}}
   \put(2.5,0.75){[}
   \put(4.5,0.75){(}
   \put(7.25,0.75){)]}
   \put(0.5,3){(}
   \put(4.5,3){(}
   \put(7.25,3){)}
   \put(3.25,3){)}
   \put(0,1.5){\framebox(8,1){}}
   \put(5,3.5){\okreis}
   \put(10,2){=}
    \multiput(12,1)(2,0){4}{\line(0,1){1}}
   \multiput(12,2)(2,0){2}{\line(0,1){1}}
    \put(16,2){\okreis}
     \put(13.5,1.25){[}
   \put(15.5,1.25){(}
   \put(18.25,1.25){)]} 
\end{picture}
\end{equation}
as well as the upside--down and left--right mirror images of
(\ref{graf5b}) and (\ref{graf5c}) and the graphs obtained by rotating
by $180^\circ$ in the drawing plane. In general, with every graphical
rule, where the graph is build from elementary graphs of the above
list, the rotated as well as the upside--down and left--right mirror
images are also valid and are proven analogously. This induces a
$\mathbb{Z}_2\times\mathbb{Z}_2$ - symmetry action on all graphical
identities given below, which in fact is already apparent in the axioms of a
quasi--triangular quasi--Hopf algebra given in Section 2.1 by taking
$\G^{op}$, $\G_{cop}$ or $\G_{cop}^{op}$ instead of $\G$.
 
Finally we  point out the important {\it ``pull through''} rule saying that
morphisms built from representation matrices of special elements in
$\G$ (like the braiding (\ref{graf3}) or the reassociator
(\ref{graf4})) always ``commute'' with all other intertwiners in the
appropriate sense i.e. by changing colours and orderings accordingly.
For example one has
\newcommand{\biega}[3]{
  \begin{picture}(4,6)(0.25,0)
     \put(0,1.25){\framebox(4,1.5){$h'$}}
   \put(2,2.75){\line(0,1){1.25}}
   \multiput(1,0)(2,0){2}{\line(0,1){1.25}}
   \put(3,-0.5){\makebox(0,0){$#3$}}
   \put(1,-0.5){\makebox(0,0){$#2$}}
   \put(2,4.5){\makebox(0,0){$#1$}}  
\end{picture}
  }
\newcommand{\biegb}[3]{
\begin{picture}(4,6)(0.25,0)
    \put(0,1.25){\framebox(4,1.5){$h$}}
   \put(2,0){\line(0,1){1.25}}
   \multiput(1,2.75)(2,0){2}{\line(0,1){1.25}}
   \put(1,4.5){\makebox(0,0){$#1$}}
   \put(3,4.5){\makebox(0,0){$#2$}}
   \put(2,-0.5){\makebox(0,0){$#3$}} 
\end{picture}
  }
\begin{equation}
  \label{pull}
  \begin{picture}(32,8)
  \put(0,0.5){\makebox(0,0){$X$}}
  \put(2,0.5){\makebox(0,0){$U$}}
  \put(0,1){\Rmatrix}
  \put(0,3){\line(0,1){4}}
  \multiput(2,3)(0,2.75){2}{\line(0,1){1.25}}
  \put(1,4.25){\framebox(2,1.5){$h$}}
  \put(0,7.5){\makebox(0,0){$U$}}
  \put(2,7.5){\makebox(0,0){$V$}}
  \put(4.5,4){\makebox(0,0){=}}
  \put(7,0.5){\makebox(0,0){$X$}}
  \put(9,0.5){\makebox(0,0){$U$}}
  \put(7,5){\Rmatrix}
  \put(9,1){\line(0,1){4}}
  \multiput(7,1)(0,2.75){2}{\line(0,1){1.25}}
  \put(6,2.25){\framebox(2,1.5){$h$}}
  \put(7,7.5){\makebox(0,0){$U$}}
  \put(9,7.5){\makebox(0,0){$V$}}
  \put(9.5,4){\makebox(0,0){,}}
\put(12,0){\boxvviii}
\multiput(13,3)(6,0){2}{\line(0,1){4}}
\put(14,3){\biega{}{}{}}
\put(15.5,6.25){(}
\put(19.25,6.25){)}
\put(21.5,4){=}
\multiput(24,0)(6,0){2}{\line(0,1){4}}
\put(25,0){\biega{}{}{}}
\put(23.5,0.5){[}
\put(25.5,0.5){(}
\put(28.25,0.5){)]}
  \put(23,5){\framebox(8,1){}}
  \multiput(24,4)(3,0){2}{\line(0,1){1}}
  \multiput(24,6)(3,0){2}{\line(0,1){1}}
  \put(30,4){\line(0,1){1}}
  \put(30,6){\line(0,1){1}}
  \put(23.5,4.25){(}
  \put(26.5,6.5){(}
  \put(30.25,6.5){)}
  \put(27.25,4.25){)}
  \end{picture}
\end{equation}
In the language of categories this means that the braidings and the
reassociators provide {\it natural transformations} [ML].
%
%
%
%
\subsection{The antipode image of the R--matrix}
In this subsection we exploit the full power of  our graphical
machinery by proving various important identities involving the action
of the antipode on a quasitriangular R--matrix. We recall that for
ordinary Hopf algebras (i.e. where $\phi, \alpha$ and $\beta$ are
trivial) one has $(S\tp\id)(R) = R^{-1} = (\id\tp\Si)(R)$ and
$(S\tp S)(R) = R$. To generalize these identities to the  quasi--Hopf
case we introduce the following four elements in $\G\tp\G$
(using the notation \eqref{eq17a})
\begin{equation}
 \label{defpq}
\begin{split}
  p_\la := Y^i\,\Si(X^i\,\beta)\tp Z^i \quad \quad & p_\rho : = P^i\tp
  Q^i\,\beta\, S(R^i)  \\
  q_\la : =S(P^i)\,\alpha\, Q^i \tp R^i \quad\quad & q_\rho : = X^i\tp
  \Si(\alpha\,Z^i)\,Y^i
\end{split}
\end{equation}
These elements have already been considered by [Dr2,S], see also
Eqs. (9.20)-(9.23) of [HN]. They obey the commutation relations (for
all $a \in \G$)
\begin{equation}
 \label{eqpq}  
\begin{split}
  \cop(\az)\,p_\la\, [\Si(\ae)\tp \e] = p_\la\, [\e\tp a], \quad\quad &
  \cop(\ae)\, p_\rho\, [\e\tp S(\az)] = p_\rho\, [a \tp \e] \\
  [S(\ae)\tp\e]\, q_\la \, \cop(\az) = [\e\tp a]\, q_\la, \quad\quad &
  [\e\tp \Si(\az)]\,q_\rho\, \cop(\ae) = [a\tp\e]\, q_\rho
\end{split}
\end{equation}
see e.g. [HN,Lem.\ 9.1]. A graphical interpretation of these
identities will be given in Eqs.\ (\ref{grafpq}) below. With these
definitions we now have
\begin{thm}
  \label{SS}
  Let $(\G,\cop,\phi,S,\alpha,\beta)$ be a finite dimensional
  quasi--Hopf algebra, let $\gamma \in \G\tp\G$ be as in
  \eqref{eq116} and let $R\in\G\tp\G$ be
  quasitriangular. Then
\begin{align}
   R^{-1}= & [X^j\beta S(P^iY^j)\tp\e ] \cdot [(S\tp\id)(q^{op}_\rho\,R)] \cdot
  [(R^i\tp Q^i) \,\cop^{op}(Z^j)] \notag\\ 
      \label{graf7a}
   = &[\cop(R^j)\, (Y^i \tp Z^i)] \cdot [(\id\tp\Si)(R \,
  p_\rho)]\cdot [\e\tp \Si(\alpha Q^j X^i)P^j] \\  
  \label{graf8a}
     (S\tp S)(R)\, \gamma  = & \gamma^{21} \,R 
\end{align}
\end{thm}
A direct implication of equation \eqref{graf8a} is the following
formula, which has already been stated without proof in [AC].
\begin{cor}
\label{corfrf}
Under the conditions of Theorem \ref{SS}
let $f \in \G\tp\G$ be the twist defined in \eqref{eq118}, then
  \begin{equation}
    \label{cor01}
     f^{op} \, R \, f^{-1} = (S\tp S)(R)
  \end{equation}
\end{cor}
\begin{proof}[Proof of Corollary \ref{corfrf}.]
  Using the formula (\ref{eq118}) for $f$  and (\ref{graf8a}) one
  computes
  \begin{align*}
     f^{op} \, R &=  (S\tp S)(\cop(P^i))\,\gamma^{op} \,
     \cop^{op}(Q^i\beta R^i) \, R \\ 
       &=  (S\tp S)(\cop(P^i)) \, (S\tp S)(R) \,  \gamma \,
       \cop(Q^i\beta R^i)\\
       & = (S\tp S)(R)\, (S\tp S)(\cop^{op}(P^i)) \,\gamma\,
          \cop(Q^i\beta R^i) = (S\tp S)(R) \, f.
  \end{align*}
\end{proof}
To prepare the proof of Theorem \ref{SS} we need the following 3
Lemmata. First we have
\begin{lem} {\rm
  \label{cora}
  \begin{equation}
    \label{graf5e}
    \begin{picture}(15,8)
    \put(1,1){\ukreis}
    \put(0,1){\boxvv}
    \multiput(5,0)(2,0){2}{\line(0,1){1}}
    \put(1,4){\line(0,1){4}}
    \put(2,4){\boxd}
    \put(2.25,6.5){[(}
    \put(5.25,6.5){)}\put(7.25,6.5){]}
    \put(3,7){\okreis}
    \put(7,7){\line(0,1){1}}
    \put(10,4){\makebox(0,0){=}}
    \put(12,0){\line(0,1){8}}
    \put(14,0){\line(0,1){8}}
    \end{picture}
  \end{equation}  }
  and three mirror images.
\end{lem}
\begin{proof} This is straightforward and left to the reader.
  (Use first \eqref{graf5d}, then \eqref{graf5b} and
  \eqref{graf5c} (or the suitable mirror image) and finally \eqref{graf5}). 
\end{proof}
To give an algebraic formulation of the four identities of Lemma
\ref{cora} let us introduce the notation 
\begin{equation}
\label{grafpq}
  \begin{picture}(36,12)
   \put(2,8.5){\makebox(0,0){$P^\la_{VW}:=$}}
   \put(10,6){\line(0,1){1}}
   \put(5,7){\boxdi}
   \put(6,7){\ukreis}
   \multiput(6,10)(2,0){3}{\line(0,1){0.5}}
   \put(6,11){\makebox(0,0){$V^*$}}
\put(8,11){\makebox(0,0){$V$}}
\put(10,11){\makebox(0,0){$W$}}
\put(16,8.5){\makebox(0,0){$P^\rho_{WV}:=$}}
   \put(20,6){\line(0,1){1}}
   \put(19,7){\boxdii}
   \put(22,7){\ukreis}
   \multiput(20,10)(2,0){3}{\line(0,1){0.5}}
   \put(20,11){\makebox(0,0){$W$}}
\put(22,11){\makebox(0,0){$V$}}
\put(24,11){\makebox(0,0){${}^*V$}}
\put(2,2.5){\makebox(0,0){$Q^\la_{VW}:=$}}
   \put(10,4){\line(0,1){1}}
   \put(5,1){\boxdii}
   \put(6,4){\okreis}
   \multiput(6,0.5)(2,0){3}{\line(0,1){0.5}}
   \put(6,0){\makebox(0,0){${}^*V$}}
\put(8,0){\makebox(0,0){$V$}}
\put(10,0){\makebox(0,0){$W$}}
\put(16,2.5){\makebox(0,0){$Q^\rho_{WV}:=$}}
   \put(20,4){\line(0,1){1}}
   \put(19,1){\boxdi}
   \put(22,4){\okreis}
   \multiput(20,0.5)(2,0){3}{\line(0,1){0.5}}
   \put(20,0){\makebox(0,0){$W$}}
\put(22,0){\makebox(0,0){$V$}}
\put(24,0){\makebox(0,0){$V^*$}}
  \end{picture}
\end{equation}
In a more general scenario these morphisms have already been
introduced in [HN]. Algebraically, they are given by (using the module
notation)
\begin{align}
\label{idpq}
  P^\la_{VW} &: w \mapsto v^i\tp p_\la\cdot (v_i\tp w), &
 P^\rho_{WV}&: w \mapsto p_\rho \cdot (w\tp v_i) \tp v^i, \\
 Q^\la_{VW}&: \hat{v}\tp v\tp w \mapsto (\hat{v} \tp\id)\big(q_\la
 \cdot (v\tp w)\big), &
 Q^\rho_{WV}&: w \tp v \tp \hat{v} \mapsto (\id\tp \hat{v})\big(q_\rho
 \cdot (w\tp v)\big) \notag
\end{align}
see Eqs.\ (9.36),(9.37),(9.42) and (9.43) of [HN].
Note that the identities \eqref{eqpq} precisely reflect the fact that
these maps are morphisms in $\rep \G$. 
In this way Lemma \ref{cora} is also contained in [HN, Lem 9.1], since
it is equivalent to the four identities, respectively
\begin{align}
\label{eqvw}
  [S(p^1_\la) \tp \e]\, q_\la \, \cop(p^2_\la) = \e\tp\e, \quad \quad
  & [\e\tp\Si(p_\rho^2)]\, q_\rho \, \cop(p^1_\rho) = \e\tp\e \\
 \notag \cop (q^2_\la) \, p_\la \, [\Si (q_\la^1) \tp \e] =\e\tp\e,
  \quad\quad
  & \cop(q_\rho^1) \, p_\rho\, [\e\tp S(q_\rho^2)] =\e\tp\e.
\end{align}
Next we define intertwiners 
$g_{VW}: ({}^*W\tp{}^*V)\tp (V\tp W)\rightarrow \mathbb{C}$ and 
$d_{VW}: \mathbb{C} \rightarrow  (V\tp W)\tp ({}^*W\tp{}^*V)$ by
\begin{equation}
  \label{graf6}
  \begin{picture}(14,6)(1,0)
  \put(1.5,3){$g_{VW}:=$}
  \put(6,4){\gokreis}
  \put(8,4){\okreis}
  \put(5,1){\boxvii}
  \put(6,0.5){\makebox(0,0){${}^*W$}}
  \put(8,0.5){\makebox(0,0){${}^*V$}}
  \put(10,0.5){\makebox(0,0){$V$}}
  \put(12,0.5){\makebox(0,0){$W$}}
  \end{picture}
  \begin{picture}(14,6)(1,1)
  \put(1.5,4){$d_{VW}:=$}
  \put(6,3){\gukreis}
  \put(8,3){\ukreis}
  \put(5,3){\boxvi}
  \put(6,6.5){\makebox(0,0){$V$}}
  \put(8,6.5){\makebox(0,0){$W$}}
  \put(10,6.5){\makebox(0,0){${}^*W$}}
   \put(12,6.5){\makebox(0,0){${}^*V$}}
  \end{picture}  
\end{equation}
One directly verifies that
\begin{align}
  \label{graf6al}
  &g_{VW}(\hat{w}\tp \hat{v} \tp v\tp w) = \lpa \hat{v}\tp \hat{w} | \,
    \gamma_{VW} \,(v\tp w) \rpa\\
& d_{VW} (1) = \sum_i \big(\delta_{VW} \, (v_i\tp w_i)
 \big) \tp \big(w^i\tp v^i\big), \notag
\end{align}
where $ \hat{v} \in {}^*V,\,
     \hat{w} \in {}^* W,\, v\in V,\, w\in W$ and 
 where $\{v_i\tp w_i\}$ is a basis of $V\tp W$ with dual basis
 $\{v^i\tp w^i\} \in {}^*V\tp {}^*W$.
Here $\gamma$,$\delta\in \G\tp\G$ are 
 given in (\ref{eq116}),(\ref{eq117}) and 
$\gamma_{VW}=
     (\pi_V\tp\pi_W)(\gamma)$,
$\delta_{VW}
 =(\pi_V\tp\pi_W)(\delta)$. 
We remark that in terms of these intertwiners the identities
(\ref{eq119b}) may now be depicted as
\begin{equation}
\label{grafdelta}
  \begin{picture}(33,7)
   \put(0,1.5){\makebox(0,0){${}^*(V\tp W)$}}  
   \put(4,1.5){\makebox(0,0){$V\tp W$}}
   \put(2,2){\oval(4,5)[t]}
   \put(6,3.5){\makebox(0,0){=}} 
   \put(9,0.5){\makebox(0,0){${}^*(V\tp W)$}}  
   \put(13,0.5){\makebox(0,0){$V\tp W$}}
    \put(9,1){\line(0,1){1.25}}
    \put(13,1){\line(0,1){1.25}}
    \put(7.5,2.25){\framebox(3,1.5){$f^{-1}$}}
    \put(11.5,2.25){\framebox(3,1.5){$\id_{V\tp W}$}}
    \put(8,3.75){\line(0,1){1.25}}
    \put(10,3.75){\line(0,1){1.25}}
    \put(12,3.75){\line(0,1){1.25}}
    \put(14,3.75){\line(0,1){1.25}}
    \put(7.5,4.25){(} \put(11.5,4.25){(} 
    \put(10.25,4.25){)} \put(14.25,4.25){)}  
   \put(7.5,5){\framebox(7,1.5){$g_{VW}$}}
    \put(15,3.5){\makebox(0,0){;}}
      \put(18,4.5){\makebox(0,0){$(V\tp W)^*$}}  
   \put(22,4.5){\makebox(0,0){$V\tp W$}}
   \put(20,4){\oval(4,5)[b]}
   \put(24,3.5){\makebox(0,0){=}} 
   \put(27,6){\makebox(0,0){$(V\tp W)^*$}}  
   \put(31,6){\makebox(0,0){$V\tp W$}}
    \put(27,4.25){\line(0,1){1.25}}
    \put(31,4.25){\line(0,1){1.25}}
    \put(25.5,2.75){\framebox(3,1.5){$f$}}
    \put(29.5,2.75){\framebox(3,1.5){$\id_{V\tp W}$}}
    \put(26,1.5){\line(0,1){1.25}}
    \put(28,1.5){\line(0,1){1.25}}
    \put(30,1.5){\line(0,1){1.25}}
    \put(32,1.5){\line(0,1){1.25}}
    \put(25.5,2){(} \put(29.5,2){(} 
    \put(28.25,2){)} \put(32.25,2){)}  
   \put(25.5,0){\framebox(7,1.5){$d_{V^*W^*}$}}
  \end{picture}
\end{equation}
Moreover, we have the following
\begin{lem}
\label{lemgamma}
{\rm
  \begin{equation}
  \label{graf6a}
  \begin{picture}(14,6)(1,0)
  \put(1.5,3){$g_{VW}=$}
  \put(6,4){\gokreis}
  \put(8,4){\okreis}
  \put(5,1){\boxvvi}
  \put(6,0.5){\makebox(0,0){${}^*W$}}
  \put(8,0.5){\makebox(0,0){${}^*V$}}
  \put(10,0.5){\makebox(0,0){$V$}}
  \put(12,0.5){\makebox(0,0){$W$}}
  \end{picture}
  \begin{picture}(14,6)(1,1)
  \put(1.5,4){$d_{VW}=$}
  \put(6,3){\gukreis}
  \put(8,3){\ukreis}
  \put(5,3){\boxvvii}
  \put(6,6.5){\makebox(0,0){$V$}}
  \put(8,6.5){\makebox(0,0){$W$}}
  \put(10,6.5){\makebox(0,0){${}^*W$}}
   \put(12,6.5){\makebox(0,0){${}^*V$}}
  \end{picture}  
\end{equation} 
}
\end{lem}
\begin{proof}
  We prove the first identity: 
  \begin{equation*}
  \begin{picture}(32,11)
  \put(1,6){\gokreis}
  \put(3,6){\okreis}
  \put(0,3){\boxvii}
  \put(10,5.5){\makebox(0,0){$\equiv$}}
  \put(11,1.5){\boxvv}
  \put(13,4.5){\boxd}
  \put(13.25,7){[(}
  \put(16.25,7){)}
  \put(18.25,7){]}
  \put(12,4.5){\line(0,1){3}}
  \put(12,7.5){\gokreis}
  \put(14,7.5){\okreis}
    \put(21,5.5){\makebox(0,0){=}}
  \put(23,0){\boxviii}
  \put(23,3){\boxd}
   \put(30,3){\line(0,1){3}}
   \put(23,6){\boxvviii}
   \put(24,9){\gokreis}
  \put(26,9){\okreis}
  \end{picture}
\end{equation*}
\begin{equation*}
  \begin{picture}(22,8)
  \put(0,4){\makebox(0,0){=}}
  \put(2,0){\boxviii}
  \put(2,3){\boxd}
  \put(2.5,5.5){[}
  \put(4.5,5.5){(}
  \put(7.25,5.5){)]}
  \put(9,3){\line(0,1){3}}
  \put(3,6){\gokreis}
  \put(5,6){\okreis}
    \put(12,4){\makebox(0,0){$\equiv$}}
  \put(14,1.5){\boxvvi}
  \put(15,4.5){\gokreis}
  \put(17,4.5){\okreis}
  \end{picture}
\end{equation*}
where in the second equality we have plugged in the pentagon identity 
(\ref{graf5d}) and in the third we have used (\ref{graf5b}). The
second identity in \eqref{graf6a} is the upside--down mirror image of
the first one and is proved analogously.
\end{proof}
We invite the reader to check that algebraically Lemma \ref{lemgamma}
implies the following identities for $\gamma$ and $\delta$ defined in 
(\ref{eq116}) and (\ref{eq117}):
\begin{align*}
   &\gamma  = (S(\tilde{U}^i)\tp S(\tilde{T}^i))\cdot (\alpha\tp\alpha)\cdot
   (\tilde{V}^i\tp \tilde{W}^i) 
   \\
  & \delta  = (\tilde{K}^j\tp \tilde{L}^j)\cdot (\beta\tp\beta)\cdot
   (S(\tilde{N}^j)\tp S(\tilde{M}^j)),\quad\text{where} \\
  &  \tilde{T}^i\tp \tilde{U}^i \tp \tilde{V}^i \tp \tilde{W}^i  =
               (\phi\tp\e)\cdot (\cop\tp\id\tp\id)(\phi^{-1}), \\
 & \tilde{K}^j \tp \tilde{L}^j \tp \tilde{M}^j \tp \tilde{N}^j  =
              (\id\tp\id\tp\cop)(\phi^{-1})\cdot (\e\tp\phi).
\end{align*}
These identities have already been obtained by Drinfel'd [Dr2].

Finally we note the following linear isomorphisms of intertwiner
spaces holding in fact in any rigid monoidal category.
\begin{lem}
  \label{umbieg}
  Let $X,V,W$ be finite dimensional $\G$-modules. Then there exist
  linear bijections
  \begin{align}
    \label{bieg1}
   \Psi^V_{X,W}\,&:\,\, \text{Hom}_\G(X\tp W,V) \pfeil 
                 \text{Hom}_\G(X,V\tp {}^*W) \\
    \label{bieg1a}    
   \Phi^{V,W}_X \,&:\,\,\text{Hom}_\G(X,V\tp W) \pfeil
                \text{Hom}_\G({}^*V\tp X, W)
  \end{align}
  given by
{\rm
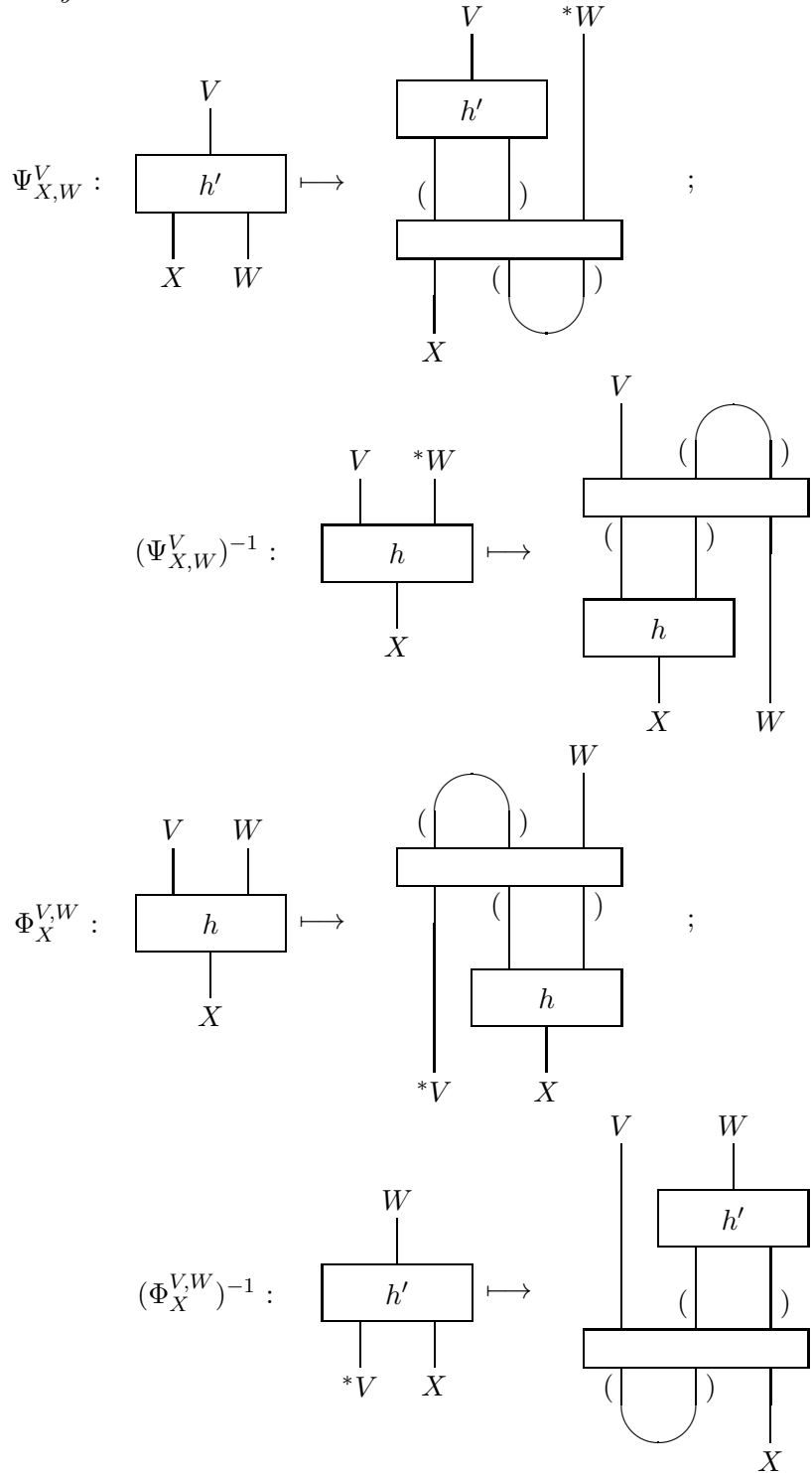
\begin{equation*}
  \label{bieg2}
  \begin{picture}(24,8)(0,0)
  \put(0,4){\makebox(0,0){$\Psi^V_{X,W}$ : }}
  \put(2,2){\biega{V}{X}{W}}
  \put(7,4){\makebox(0,0){$\longmapsto$}}
  \put(9,4){\biega{V}{}{}}
   \put(9,1){\boxdii}
   \put(12,1){\ukreis}
   \put(10,0){\line(0,1){1}}
   \put(10,-0.5){\makebox(0,0){$X$}}
   \put(14,4){\line(0,1){4}}
   \put(14,8.5){\makebox(0,0){${}^*W$}}
    \put(17,4){\makebox(0,0){; }}
  \end{picture}
\end{equation*}
\begin{equation*}
  \begin{picture}(24,8)(-4,0)
     \put(0,4){\makebox(0,0){$(\Psi^V_{X,W})^{-1}$ : }}
  \put(3,2){\biegb{V}{{}^*W}{X}}
  \put(8,4){\makebox(0,0){$\longmapsto$}}
  \put(10,0){\biegb{}{}{X}}
   \put(10,4){\boxdi}
   \put(13,7){\okreis}
   \put(11,7){\line(0,1){1}}
   \put(11,8.5){\makebox(0,0){$V$}}
   \put(15,0){\line(0,1){4}}
   \put(15,-0.5){\makebox(0,0){$W$}}
  \end{picture}
\end{equation*}
\begin{equation*}
  \label{bieg3}
  \begin{picture}(24,8)(0,0)
  \put(0,4){\makebox(0,0){$\Phi^{V,W}_X$ : }}
  \put(2,2){\biegb{V}{W}{X}}
  \put(7,4){\makebox(0,0){$\longmapsto$}}
  \put(11,0){\biegb{}{}{X}}
   \put(9,4){\boxdii}
   \put(10,7){\okreis}
   \put(14,7){\line(0,1){1}}
   \put(14,8.5){\makebox(0,0){$W$}}
   \put(10,0){\line(0,1){4}}
   \put(10,-0.5){\makebox(0,0){${}^*V$}}
  \put(17,4){\makebox(0,0){; }}
  \end{picture}
\end{equation*}
\begin{equation*}
  \begin{picture}(24,8)(-4,0)
  \put(0,4){\makebox(0,0){$(\Phi^{V,W}_X)^{-1}$ : }}
  \put(3,2){\biega{W}{{}^*V}{X}}
  \put(8,4){\makebox(0,0){$\longmapsto$}}
  \put(12,4){\biega{W}{}{}}
   \put(10,1){\boxdi}
   \put(11,1){\ukreis}
   \put(15,0){\line(0,1){1}}
   \put(15,-0.5){\makebox(0,0){$X$}}
   \put(11,4){\line(0,1){4}}
   \put(11,8.5){\makebox(0,0){$V$}}
  \end{picture}
\end{equation*}
   }
\end{lem}
\begin{proof}
  We prove $\Psi^V_{X,W} \circ (\Psi^V_{X,W})^{-1} = \id$ by determine
  its action on $h \in Hom_\G(X,V\tp {}^*W)$ as follows: 
  
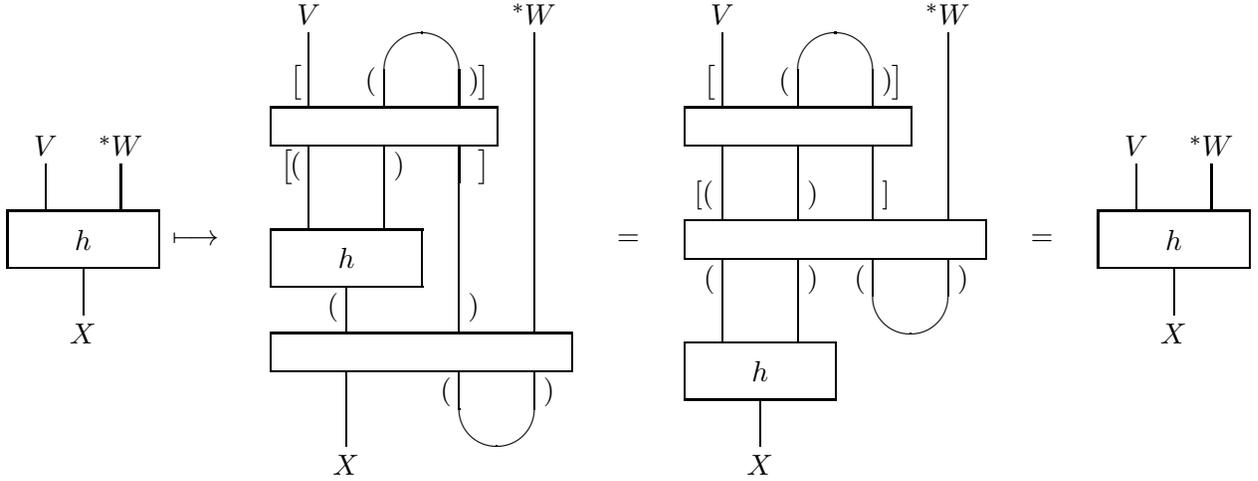
\begin{equation*}
    \begin{picture}(36,12)(0,0)
   \put(0,3.5){\biegb{V}{{}^*W}{X}}
  \put(5,5.5){\makebox(0,0){$\longmapsto$}}
  \put(7,3){\biegb{}{}{}}
   \put(7,7){\boxdi}
   \put(7.5,9.5){\big[}
   \put(12.5,9.5){\big]}
   \put(7.25,7.25){\big[}
   \put(12.5,7.25){\big]}
   \put(10,10){\okreis}
   \put(8,10){\line(0,1){1}}
   \put(8,11.5){\makebox(0,0){$V$}}
   \put(12,3){\line(0,1){4}}
   \put(14,11.5){\makebox(0,0){${}^*W$}}
   \put(9,0){\line(0,1){2}}
   \put(9,-0.5){\makebox(0,0){$X$}}
   \put(7,2){\framebox(8,1){}}
   \put(9,3){\line(0,1){1}}
   \multiput(12,1)(2,0){2}{\line(0,1){1}}
   \multiput(12,3)(2,0){2}{\line(0,1){1}}
   \put(14,4){\line(0,1){7}}
   \put(12,1){\ukreis}
   \put(11.5,1.25){(}
   \put(8.5,3.5){(}
   \put(12.25,3.5){)}
   \put(14.25,1.25){)}
\put(16.5,5.5){\makebox(0,0){=}}
  \put(18,0){\biegb{}{}{X}}
   \put(18,7){\boxd}
   \put(18.5,9.5){\big[}
   \put(23.5,9.5){\big]}
   \put(20.5,9.5){(}
   \put(23.25,9.5){)}
   \put(21,10){\okreis}
   \put(19,10){\line(0,1){1}}
   \put(19,11.5){\makebox(0,0){$V$}}
   \put(25,11.5){\makebox(0,0){${}^*W$}}
   \put(25,7){\line(0,1){4}}
   \put(18,4){\boxviii}
   \put(23,4){\ukreis}
\put(27.5,5.5){\makebox(0,0){=}}
   \put(29,3.5){\biegb{V}{{}^*W}{X}}
    \end{picture}
  \end{equation*}
where in the first equality we have used a ``pull through''
rule for $h$,
and in the second equality a left--right mirror image of
\eqref{graf5e}. Analogously 
one shows that $\Psi^{-1} \circ \Psi = \id$ and $\Phi\circ \Phi^{-1} =
\id = \Phi^{-1} \circ \Phi$.
\end{proof}

We are now in the position to prove Eqs.~(\ref{graf7a}) and
(\ref{graf8a}) of Theorem~\ref{SS} by rewriting  them as graphical
identities as follows: 
\begin{lem}
  \label{lem02}
  For all finite dimensional $\G$-modules the inverse braiding
  $B_{UV}^{-1}$
  obeys
{\rm 
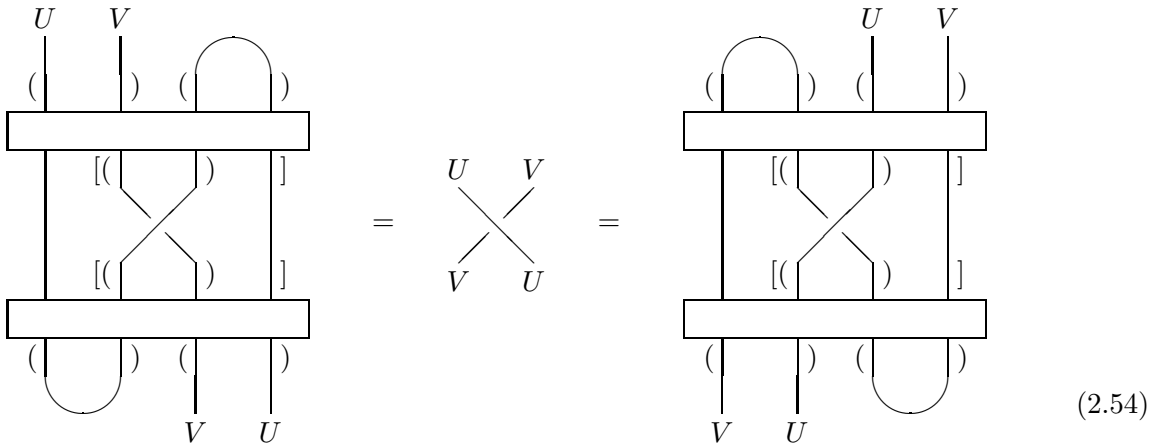
  \begin{equation}
    \label{graf7}
     \begin{picture}(28,11)
        \put(5,-0.5){\makebox(0,0){$V$}}
        \put(7,-0.5){\makebox(0,0){$U$}}
        \put(0,0){\bildd}
        \put(3,10.5){\makebox(0,0){$V$}}
        \put(1,10.5){\makebox(0,0){$U$}}
        \put(10,5){\makebox(0,0){=}}
        \put(12,3.5){\makebox(0,0){$V$}}
        \put(14,3.5){\makebox(0,0){$U$}}
        \put(12,4){\iRmatrix}
        \put(14,6.5){\makebox(0,0){$V$}}
        \put(12,6.5){\makebox(0,0){$U$}}
        \put(16,5){\makebox(0,0){=}}
        \put(19,-0.5){\makebox(0,0){$V$}}
        \put(21,-0.5){\makebox(0,0){$U$}}
        \put(25,10.5){\makebox(0,0){$V$}}
        \put(23,10.5){\makebox(0,0){$U$}}
         \put(18,0){\bildc}
    \end{picture}
 \end{equation} }
and the (left) conjugate braiding $B_{{}^*U{}^*V}$ obeys
{\rm
 \begin{equation}
    \begin{picture}(22,8)(0,-1)
       \put(1,-0.5){\makebox(0,0){${}^*U$}}
       \put(3,-0.5){\makebox(0,0){${}^*V$}}
       \put(5,-0.5){\makebox(0,0){$U$}}
       \put(7,-0.5){\makebox(0,0){$V$}}
       \put(1,0){\Rmatrix}
       \multiput(5,0)(2,0){2}{\line(0,2){2}}
       \put(0,2){\bilde}
       \put(10,4){=}
       \put(13,-0.5){\makebox(0,0){${}^*U$}}
       \put(15,-0.5){\makebox(0,0){${}^*V$}}
       \put(17,-0.5){\makebox(0,0){$U$}}
       \put(19,-0.5){\makebox(0,0){$V$}}
       \put(17,0){\Rmatrix}
       \multiput(13,0)(2,0){2}{\line(0,2){2}}
       \put(12,2){\bilde}
    \end{picture} 
    \label{graf8}
  \end{equation} }
Taking $\G$ itself as a $\G$-module yields 
$\eqref{graf7} \Leftrightarrow \eqref{graf7a}$ and
 $\eqref{graf8} \Leftrightarrow \eqref{graf8a}$ and therefore the
 above identities prove Theorem \ref{SS}.
\end{lem}
\begin{proof}
To prove the first equation in \eqref{graf7} note, that \eqref{eq19}
and \eqref{eq112} imply the identity
\begin{equation}
  \label{graf9}
  \begin{picture}(15,12)
  \put(1,0.5){\makebox(0,0){${}^*U$}}
  \put(3,0.5){\makebox(0,0){$U$}}
  \put(5,0.5){\makebox(0,0){$V$}}
  \put(3,1){\Rmatrix}
  \put(1,1){\line(0,1){2}}
  \put(0,3){\boxdii}
  \put(1,6){\Rmatrix}
  \put(5,6){\line(0,1){2}}
  \put(0,8){\boxdi}
  \put(3,11){\okreis}
  \put(1,11.5){\makebox(0,0){$V$}}
  \put(8,5){\makebox(0,0){$=$}}
  \put(9,4){\boxdii}
  \put(10,7){\okreis}
  \put(14,7){\line(0,1){1}}
    \multiput(10,3)(2,0){3}{\line(0,1){1}}
  \put(10,2.5){\makebox(0,0){${}^*U$}}
  \put(12,2.5){\makebox(0,0){$U$}}
  \put(14,2.5){\makebox(0,0){$V$}}
  \end{picture}
\end{equation}
Now we apply the isomorphism $(\Phi^{U,V}_{U\tp V})^{-1}$ of Lemma
\ref{umbieg} to both sides of \eqref{graf9} to obtain
\begin{equation}
  \label{graf9a}
  \begin{picture}(32,15)
  \put(5,0.5){\makebox(0,0){$U$}}
  \put(7,0.5){\makebox(0,0){$V$}}
  \put(1,3){\ukreis}
  \put(5,1){\Rmatrix}
  \put(0,3){\boxvii}
  \put(3,6){\Rmatrix}
  \multiput(1,6)(6,0){2}{\line(0,1){2}}
  \put(2.25,8.25){[} \put(7.5,8.25){]}
  \put(2.25,10.5){[} \put(7.5,10.5){]}
  \put(2,8){\boxdi}
  \put(1,8){\line(0,1){3}}
  \put(5,11){\okreis}
  \put(1,11.5){\makebox(0,0){$U$}}
  \put(3,11.5){\makebox(0,0){$V$}}
  \put(9.5,7){\makebox(0,0){$\equiv$}}
  \put(16,0.5){\makebox(0,0){$U$}}
  \put(18,0.5){\makebox(0,0){$V$}}
  \put(12,1){\ukreis}
  \put(16,4){\Rmatrix}
  \put(11,1){\boxvv}
  \put(13,6){\boxdii}
   \put(14,9){\Rmatrix}
   \put(18,9){\line(0,1){2}}
  \put(13,11){\boxdi}
  \put(13.25,11.25){[} \put(18.5,11.25){]}
  \put(13.25,13.5){[} \put(18.5,13.5){]}
  \put(12,4){\line(0,1){10}}
  \put(14,4){\line(0,1){2}}
  \put(16,14){\okreis}
  \put(12,14.5){\makebox(0,0){$U$}}
  \put(14,14.5){\makebox(0,0){$V$}}
  \put(20.5,7){\makebox(0,0){$=$}}
  \put(27,2.5){\makebox(0,0){$U$}}
  \put(29,2.5){\makebox(0,0){$V$}}
  \put(23,3){\ukreis}
  \put(22,3){\boxvv}
  \multiput(23,6)(2,0){4}{\line(0,1){1}}
  \put(24,7){\boxdii}
  \put(24.25,7.25){[} \put(29.5,7.25){]}
  \put(24.25,9.5){[} \put(29.5,9.5){]}
  \put(23,7){\line(0,1){4}}
  \put(29,10){\line(0,1){1}}
  \put(25,10){\okreis}
\put(23,11.5){\makebox(0,0){$U$}}
  \put(29,11.5){\makebox(0,0){$V$}}
  \end{picture}
\end{equation}
where in the left identity we have used a ``pull through'' rule for
the braiding. By 
the identity \eqref{graf5c} the top of the graph in the middle of
\eqref{graf9a} may be replaced by
\begin{equation}
  \label{graf9b}
  \begin{picture}(20,4)
  \put(0,0){\line(0,1){4}}
  \put(1,0){\boxdi}
  \put(4,3){\okreis}
  \put(2,3){\line(0,1){1}}
  \put(1.25,0.25){[}   \put(6.25,0.25){]}
    \put(1.25,2.5){[}    \put(6.5,2.5){]}
  \put(9,2){\makebox(0,0){=}}
  \put(11,0){\boxvi}
   \put(16,3){\okreis}
  \multiput(12,3)(2,0){2}{\line(0,1){1}}
  \put(19.5,1.5){\makebox(0,0){,}}
  \end{picture}
\end{equation}
and using Lemma \ref{cora} for the r.h.s. of \eqref{graf9a} we end up with 
\begin{equation}
  \label{graf9c}
  \begin{picture}(12,12)
  \put(5,0.5){\makebox(0,0){$U$}}
  \put(7,0.5){\makebox(0,0){$V$}}
  \put(1,3){\ukreis}
  \put(5,1){\Rmatrix}
  \put(0,3){\boxvii}
  \put(3,6){\Rmatrix}
  \multiput(1,6)(6,0){2}{\line(0,1){2}}
  \put(0,8){\boxvi}
  \put(5,11){\okreis}
  \put(1,11.5){\makebox(0,0){$U$}}
  \put(3,11.5){\makebox(0,0){$V$}}
  \put(10,6){\makebox(0,0){$=$}}
  \multiput(12,4)(2,0){2}{\line(0,1){7}}
   \put(12,3.5){\makebox(0,0){$U$}}
  \put(14,3.5){\makebox(0,0){$V$}}
  \end{picture}
\end{equation}
Hence we get the first equality in \eqref{graf7}\footnote{By finite
  dimensionality it suffices to prove the left inverse property.}.
Analogously one
starts with  \eqref{eq110},\eqref{eq112} to get 
\begin{equation}
  \label{graf10}
  \begin{picture}(16,12)
  \put(1,0.5){\makebox(0,0){$U$}}
  \put(3,0.5){\makebox(0,0){$V$}}
  \put(5,0.5){\makebox(0,0){$V^*$}}
  \put(1,1){\Rmatrix}
  \put(5,1){\line(0,1){2}}
  \put(0,3){\boxdi}
  \put(3,6){\Rmatrix}
  \put(1,6){\line(0,1){2}}
  \put(0,8){\boxdii}
  \put(1,11){\okreis}
  \put(5,11.5){\makebox(0,0){$U$}}
  \put(8,5){\makebox(0,0){$=$}}
 \put(10,4){\boxdi}
  \put(13,7){\okreis}
  \put(11,7){\line(0,1){1}}
   \multiput(11,3)(2,0){3}{\line(0,1){1}} 
  \put(11,2.5){\makebox(0,0){$U$}}
  \put(13,2.5){\makebox(0,0){$V$}}
  \put(15,2.5){\makebox(0,0){$V^*$}}
  \end{picture}
\end{equation}
where we have used that the trivial  identification ${}^*(V^*) =
V$.  Taking the mirror image of the above proof yields the second
equality in (\ref{graf7}).  

Eq.~(\ref{graf8}) follows from Lemma \ref{corantipode} below by
putting $h' = B_{UV}$ and $h=B_{{}^*U{}^*V}$.

The identifications 
\eqref{graf7} $\Leftrightarrow$ \eqref{graf7a} and
 \eqref{graf8} $\Leftrightarrow$ \eqref{graf8a} are straight forward
 and are left to the reader. This concludes the proof of Lemma
 \ref{lem02} and therefore also of Theorem \ref{SS}.
\end{proof}
We end this section with a Lemma used in the above proof, which will
also be used in the next section.
\newcommand{\boxh}[1]{
  \begin{picture}(4,2)(0.25,0)
    \put(0,0.25){\framebox(4,1.5){$#1$}}
    \multiput(1,0)(2,0){2}{\line(0,1){0.25}}
    \multiput(1,1.75)(2,0){2}{\line(0,1){0.25}}
  \end{picture}
} 
\newcommand{\bildci}{
\begin{picture}(8,10)(0.25,0)
\put(5,1){\ukreis}
\put(0,1){\boxvii}
\put(1,4){\liniez}
\put(2,4){\boxh{h'}}
\put(7,4){\liniez}
\put(0,6){\boxvi}
\put(1,9){\okreis}
\multiput(5,9)(2,0){2}{\line(0,1){1}}
\multiput(1,0)(2,0){2}{\line(0,1){1}}
\end{picture}
}
\newcommand{\bilddi}{
\begin{picture}(8,10)(0.25,0)
\put(1,1){\ukreis}
\put(0,1){\boxvii}
\put(1,4){\liniez}
\put(2,4){\boxh{h}}
\put(7,4){\liniez}
\put(0,6){\boxvi}
\put(5,9){\okreis}
\multiput(1,9)(2,0){2}{\line(0,1){1}}
\multiput(5,0)(2,0){2}{\line(0,1){1}}
\end{picture}
}
\begin{lem}
  \label{corantipode}
  Let $V,W,X,Y$ be finite dimensional $\G$--modules with intertwiners
  \begin{equation*}
     h:\, {}^*X\tp{}^*Y \pfeil {}^*W\tp{}^*V, \quad
     h': \, V\tp W \pfeil Y\tp X 
  \end{equation*}
 then the following two identities are equivalent
{\rm
\begin{equation}
  \label{graf11}
  \begin{picture}(22,8)
       \put(1,0.5){\makebox(0,0){${}^*X$}}
       \put(3,0.5){\makebox(0,0){${}^*Y$}}
       \put(5,0.5){\makebox(0,0){$V$}}
       \put(7,0.5){\makebox(0,0){$W$}}
       \put(0,1){\boxh{h}}
       \multiput(5,1)(2,0){2}{\line(0,2){2}}
       \put(0,3){\bilde}
       \put(10,4){=}
       \put(13,0.5){\makebox(0,0){${}^*X$}}
       \put(15,0.5){\makebox(0,0){${}^*Y$}}
       \put(17,0.5){\makebox(0,0){$V$}}
       \put(19,0.5){\makebox(0,0){$W$}}
       \put(16,1){\boxh{h'}}
       \multiput(13,1)(2,0){2}{\line(0,2){2}}
       \put(12,3){\bilde}
    \end{picture} 
\end{equation} 
  \begin{equation}
    \label{graf12}
     \begin{picture}(34,12)
        \put(7,0.5){\makebox(0,0){$V$}}
        \put(5,0.5){\makebox(0,0){${}^*Y$}}
        \put(0,1){\bilddi}
        \put(3,11.5){\makebox(0,0){${}^*W$}}
        \put(1,11.5){\makebox(0,0){$X$}}
        \put(10,6){\makebox(0,0){=}}
        \put(13,0.5){\makebox(0,0){${}^*Y$}}
        \put(15,0.5){\makebox(0,0){$V$}}
        \put(17,11.5){\makebox(0,0){$X$}}
        \put(19,11.5){\makebox(0,0){${}^*W$}}
         \put(12,1){\bildci}
      \put(21.5,6){\makebox(0,0){=:}}
       \put(23,5){\boxh{H'}}
       \multiput(24,4)(2,0){2}{\line(0,1){1}}
       \multiput(24,7)(2,0){2}{\line(0,1){1}}
        \put(24,3.5){\makebox(0,0){${}^*Y$}}
        \put(26,3.5){\makebox(0,0){$V$}}
        \put(24,8.5){\makebox(0,0){$X$}}
        \put(26,8.5){\makebox(0,0){${}^*W$}}
    \end{picture}
 \end{equation}
}
\end{lem}
\begin{proof}
  Using \eqref{graf6}, (\ref{graf6a}), a ``pull through'' rule for $h$
  and $h'$ and 
  rebracketing the source legs, Eq.\ \eqref{graf11} is equivalent to
  \begin{equation}
    \label{graf13}
    \begin{picture}(22,11)
   \put(1,0.5){\makebox(0,0){${}^*X$}}
       \put(3,0.5){\makebox(0,0){${}^*Y$}}
       \put(5,0.5){\makebox(0,0){$V$}}
       \put(7,0.5){\makebox(0,0){$W$}}
       \put(0,1){\boxdii}
       \put(0.25,1.25){[} \put(5.5,1.25){]}
       \put(0.25,3.5){[} \put(5.5,3.5){]}
       \put(0,4){\boxh{h}}
       \multiput(5,4)(2,0){2}{\line(0,2){2}}
       \put(0,6){\boxdi}
        \put(0.25,8.5){[} \put(5.5,8.5){]}
       \put(1,9){\gokreis}
       \put(3,9){\okreis}
       \put(7,1){\line(0,1){8}}
       \put(10,5){=}
       \put(13,0.5){\makebox(0,0){${}^*X$}}
       \put(15,0.5){\makebox(0,0){${}^*Y$}}
       \put(17,0.5){\makebox(0,0){$V$}}
       \put(19,0.5){\makebox(0,0){$W$}}
        \put(12,1){\boxvviv}
       \put(16,4){\boxh{h'}}
       \put(13,6){\line(0,1){3}}
       \put(14,6){\boxdii}
       \put(13,9){\gokreis}
       \put(15,9){\okreis}
        \put(14.25,8.5){[} \put(19.5,8.5){]}
       \multiput(13,4)(2,0){2}{\line(0,2){2}}    
    \end{picture}
  \end{equation}
Note that
the l.h.s. of (\ref{graf13}) is of the form $\Psi^{-1}(\cdots)$
(with a ``white'' target leg).
Now we apply the isomorphism $\Psi^{\mathbb{C}}_{[{}^*X\tp({}^*Y\tp V)],W}$
of Lemma \ref{umbieg} to both sides of \eqref{graf13}. Then using for the
bottom part of the r.h.s. the identity 
\begin{equation}
  \label{graf14a}
  \begin{picture}(24,8)
  \put(0.5,1.25){[}  \put(2.5,1.25){(}  \put(6.5,1.25){\{}
   \put(5.25,1.25){)]}  \put(9.25,1.25){\}}
  \multiput(1,0)(2,0){3}{\line(0,1){2}}
  \multiput(7,1)(2,0){2}{\line(0,1){1}}
  \put(7,1){\ukreis}
  \put(0,2){\framebox(10,1){}}
  \multiput(1,3)(2,0){5}{\line(0,1){1}}
   \put(9,4){\line(0,1){3}}
   \put(0,4){\boxvviv}
    \put(0.5,6.5){\big\{}
    \put(7.75,6.5){\big\}}
    \put(0,4.25){\big\{}
     \put(7.5,4.25){\big\}}
       \put(11.5,3.5){\makebox(0,0){=}}
  \put(14,0){\line(0,1){4}}
  \multiput(16,0)(2,0){2}{\line(0,1){1}}
  \put(20,1){\ukreis}
  \put(15,1){\boxvvi}
   \put(15,3.5){\big\{}
    \put(22.25,3.5){\big\}}
  \put(13,5){\framebox(10,1){}}
    \put(13.5,6.5){\{}
    \put(20.25,6.5){)]\big\}}
    \put(15.5,6.5){[} \put(17.5,6.5){(}
  \multiput(14,4)(2,0){5}{\line(0,1){1}}
    \multiput(14,6)(2,0){5}{\line(0,1){1}}
  \end{picture}
\end{equation}
and a pull through rule to raise the upper box of the r.h.s. of
\eqref{graf14a} to the top we conclude that (\ref{graf13}) is
equivalent to 
\begin{equation}
  \label{graf14}
  \begin{picture}(32,14)
 \put(1,0.5){\makebox(0,0){${}^*X$}}
       \put(3,0.5){\makebox(0,0){${}^*Y$}}
       \put(5,0.5){\makebox(0,0){$V$}}
       \put(1,9.5){\makebox(0,0){${}^*W$}}
       \put(0,1){\boxdii}
       \put(0,4){\boxh{h}}
       \put(5,4){\line(0,1){2}}
       \put(0,6){\boxdi}
       \put(3,9){\okreis}
       \put(7.5,5){=}
      \put(10,0.5){\makebox(0,0){${}^*X$}}
       \put(12,0.5){\makebox(0,0){${}^*Y$}}
       \put(14,0.5){\makebox(0,0){$V$}}
       \put(18,13.5){\makebox(0,0){${}^*W$}}
        \put(11,1){\boxvvi}
        \put(16,1){\ukreis}
       \put(10,1){\line(0,1){8}}
       \put(18,4){\line(0,1){5}}
       \put(13,4){\boxh{h'}}
       \put(12,4){\line(0,1){2}}
       \put(11,6){\boxdii}
       \put(12,9){\okreis}        
      \put(10,12){\gokreis}
      \put(10,9){\line(0,1){2}}
         \put(16,9){\line(0,1){2}}   \put(18,9){\line(0,1){2}}
       \put(9,11){\framebox(10,1){}}
      \put(10,12){\line(0,1){1}}
         \put(16,12){\line(0,1){1}}   \put(18,12){\line(0,1){1}}
       \put(10.5,8.5){\big\{} \put(11,8.5){[} \put(16.5,8.5){]}
       \put(18.25,8.5){\big\}}
       \put(9.5,12.5){\{}    \put(16.25,12.5){\}}
              \put(15.5,10.25){\{}       \put(18.25,10.25){\}}
     \put(20.5,5){$\equiv$}
      \put(23,2.5){\makebox(0,0){${}^*X$}}
       \put(25,2.5){\makebox(0,0){${}^*Y$}}
       \put(27,2.5){\makebox(0,0){$V$}}
       \put(27,10.5){\makebox(0,0){${}^*W$}}
      \put(23,3){\line(0,1){3}}
      \multiput(25,3)(2,0){2}{\line(0,1){1}}
     \put(24,4){\boxh{H'}}
     \put(24.5,3.25){(} \put(27.25,3.25){)}    
     \put(22,6){\boxdii}
     \put(27,9){\line(0,1){1}}
     \put(23,9){\okreis}
  \end{picture}
\end{equation}
The proof is finished by applying the isomorphism
$(\Phi^{-1})^{X,{}^*W}_{{}^*Y\tp V}$.
\end{proof}

%

\section{Doubles of quasi-Hopf algebras}

\unitlength0.5cm
\newcommand{\bildL}{
\begin{picture}(6,4)(0.25,0) 
\multiput(1,0)(2,0){3}{\line(0,1){1.25}}
\put(0,1.25){\framebox(6,1.5){$\LL$}}
\put(3,2.75){\line(0,1){1.25}}
\end{picture} }
\newcommand{\bildR}{
\begin{picture}(6,4)(0.25,0) 
\multiput(1,2.75)(2,0){3}{\line(0,1){1.25}}
\put(0,1.25){\framebox(6,1.5){$\RR$}}
\put(3,0){\line(0,1){1.25}}
\end{picture} }
\newcommand{\bildgrossL}{
\begin{picture}(10,4)(0.25,0) 
\multiput(1,0)(4,0){3}{\line(0,1){1.25}}
\put(0,1.25){\framebox(10,1.5){$\LL$}}
\put(5,2.75){\line(0,1){1.25}}
\end{picture} }
\newcommand{\bildgrossLb}{
\begin{picture}(10,4)(0.25,0) 
\multiput(2,0)(3,0){3}{\line(0,1){1.25}}
\put(0,1.25){\framebox(10,1.5){$\LL$}}
\put(5,2.75){\line(0,1){1.25}}
\end{picture} }
%
\newcommand{\bildfuenf}{
\begin{picture}(10,3)(0.25,0) 
\multiput(1,0)(2,0){5}{\line(0,1){1}}
\multiput(1,2)(2,0){5}{\line(0,1){1}}
\put(0,1){\framebox(10,1){}}
\end{picture} }
%
%
\newcommand{\bildid}[1][id]{
  \begin{picture}(4,3)(0.25,0)
 \multiput(1,0)(2,0){2}{\line(0,1){0.75}}
 \put(0,0.75){\framebox(4,1.5){#1}}
 \put(2,2.25){\line(0,1){0.75}}
  \end{picture}
}
\newcommand{\bildidu}[1][id]{
  \begin{picture}(4,3)(0.25,0)
 \multiput(1,2.25)(2,0){2}{\line(0,1){0.75}}
 \put(0,0.75){\framebox(4,1.5){#1}}
 \put(2,0){\line(0,1){0.75}}
  \end{picture}
}

\subsection{$\D(\G)$ as an associative algebra}
In this section we review the definition of the double $\D(\G)$ of  a
quasi--Hopf algebra $\G$ as a diagonal crossed product as introduced in
[HN]. We also give a graphical description of this construction.
Consider 
\begin{equation}
  \label{eq3.0}
 \delta := (\cop\tp\id)\circ
 \cop :\, \G \pfeil \G\tp\G\tp\G 
\end{equation}
and let $\Phi \in \G^{\tp^5}$ be
 given by  
\begin{equation}
  \label{Phi}
  \Phi := [(\id\tp\cop\tp\id)(\phi) \tp\e]\, [\phi\tp\e\tp\e]\, 
           [(\delta \tp\id\tp\id)(\phi^{-1})].
\end{equation}
Then the pair $(\delta,\Phi)$ provides a {\it two--sided coaction} of $\G$
on itself as defined in [HN], 
i.e the following axioms are satisfied:
\begin{itemize}
\item[(i)] The map $\delta$ is a unital algebra morphism satisfying
               $(\ep\tp\id\tp\ep)\circ \delta = \id$.
\item[(ii)]The element $\Phi \in \G^{\tp^5}$ is invertible and fulfills
  \begin{align} 
 (\id\tp\delta\tp\id)(\delta(a))\,\Phi &= \Phi \, (\cop\tp\id\tp\cop)(
    \delta(a)),\quad \forall a \in \G    \\ \notag
    (\e\tp\Phi\tp\e)\,(\id\tp\cop\tp\id\tp\cop\tp\id)&(\Phi) \,
    (\phi\tp\e\tp\phi^{-1})
    \\ &=
    (\id^{\tp^2}\tp\delta\tp\id^{\tp^2})(\Phi)\, 
    (\cop\tp\id^{\tp^3}\tp\cop)(\Phi) \label{Phi1} \\ 
    (\id\tp\ep\tp\id\tp\ep\tp\id)(\Phi) &=
    (\ep\tp\id^{\tp^3}\tp\ep)(\Phi) = \e\tp\e\tp\e \label{Phi2} 
  \end{align}
\end{itemize}
Next we define $\Om \equiv \Omega^1 \tp \Omega^2 \tp
  \Omega^3 \tp \Omega^4 \tp \Omega^5 \, \in \G^{\tp^5}$ by\footnote{
where we have again suppressed summation symbol and indices}
\begin{equation}
\label{omega}
  \Omega := 
         (\id^{\tp^3}\tp\Si\tp\Si)(f^{45}\cdot \Phi^{-1}) = 
         (\id^{\tp^3}\tp\Si\tp\Si)(\Phi^{-1}) \cdot h^{54} 
\end{equation} 
where $f,h \in \G\tp\G$ are the twists defined in \eqref{eq118} and
\eqref{eq122}. 

Let now $\dG$ be the coalgebra dual to $\G$ with
its natural left and right $\G$--action and the nonassociative
multiplication given by $\lpa \vi \psi \mid a\rpa :=\lpa
\vi\tp\psi\mid \cop(a)\rpa$. With $\delta : \G \pfeil \G^{\tp ^3}$
being a two--sided coaction we then write 
$\vi\re a \li \psi : = (\psi\tp\id\tp\vi)(\delta (a))$.
Note that for the two--sided coaction $\delta$ in Eq.\ \eqref{eq3.0}
we have the identity
 $\vi\re a \li \psi = (\vi \arr a)\arl \psi$.
Considered as an element
of $\dG$ we also write $\hat{\e} \equiv \edG \equiv \ep$.
 The following proposition has been proven in [HN]:
\begin{prop}
  \label{prop311}
  Let $(\delta,\Phi)$ be a two--sided coaction of $\G$ on $\G$ and 
 define the diagonal crossed product  $\dG\reli \G$ to be
 the vector space $\dG\tp\G$  with multiplication rule 
   \begin{equation}
   \label{e312}
  (\vi\reli a)(\psi\reli b) := 
   \Big[(\Om^1 \arr\vi\arl\Om^5)(\Om^2\arr\psi_{(2)}\arl\Om^4)\Big] \reli
  \Big[\Om^3 (\hat{S}^{-1}(\psi_{(1)})\re a \li \psi_{(3)})\, b \Big],
\end{equation}
where we write $(\vi\reli a)$ in place of $(\vi\tp a)$ to distinguish
the new algebraic structure.

Then $\dG\reli \G$ is an associative algebra with
 unit $\hat{\e}\reli\e$, containing $\G \equiv \hat{\e}\reli \G$ as a
unital subalgebra.
\end{prop}
The associativity of the above  product follows  from
the axioms for two--sided coactions as has been proven in Theorem 10.2
of [HN]. Note that in general the subspace $\dG\reli\e$ is {\bf not} a
subalgebra of $\dG\reli\G$.
On the other hand if $\G$ is an ordinary Hopf algebra 
with $\phi \equiv \e\tp\e\tp\e$ then Eq.~\eqref{e312} becomes 
\begin{equation}
\label{e312old}
  (\vi \reli a)(\psi \reli b) = \big(\vi\psi_{(2)} \reli
  (\hat{S}^{-1}(\psi_{(1)}) \re a  \li \psi_{(3)}) b \big) 
\end{equation}
which is the standard multiplication rule in the quantum double
$\D(\G)$ [Dr1,M3]. This motivates the 
\begin{df}
\label{double}
 {\rm [HN]} The diagonal crossed product $\dG\reli\G$ defined in
  Proposition \ref{prop311}, with $(\delta,\Phi)$ given by
  \eqref{eq3.0},\eqref{Phi} 
  is called the {\bf quantum double} of $\G$, 
  denoted by $\D(\G)\equiv \dG\reli\G$.
\end{df}
We will now rewrite the multiplication \eqref{e312} given in Proposition
\ref{prop311}  using the ``generating
matrix'' formalism of the St. Petersburg school. 
In this way we will be able to give a graphical (i.e. a categorical)
description of the algebraic relations in $\dG \reli \G$ which should
convince the reader, that the multiplication given in \eqref{e312} is
indeed associative.
The following Corollary
is a generalization of [N, Lemma 5.2] and coincides with [HN,Cor.\
10.4] applied to the present scenario. 
\begin{cor} 
\label{lem3.3}
  {\rm [HN]} Let $\A$ be some unital algebra and $\gamma : \G\pfeil\A$ a unital algebra
  map. Then the relation 
  \begin{equation}
    \label{e313}
    \gamma_L (\vi \reli a) = (\vi\tp\id)(\LL) \cdot \gamma (a) 
  \end{equation}
provides a one to one correspondence between unital algebra morphisms 
$\gamma_L : \dG\reli \G \pfeil \A$ extending $\gamma$ and elements
$\LL \in \G\tp\A$ satisfying $(\ep\tp\id)(\LL) = \eA$ and 
\begin{align}
     \label{e314}
     [\eG\tp \gamma (a)]\,\LL &= [\Si(a_{(1)})\tp\eA]\,\LL\,[a_{(-1)}\tp
     \gamma (a_{(0)})], \quad \forall a \in \G\\
   \label{e315}
\LL[13]\LL[23] &= 
[\Om^5  \tp \Om^4  \tp\eA]
      [(\cop\tp\id )(\LL)] 
                   [\Om^1 \tp \Om^2 \tp \gamma (\Om^3)],
\end{align}
where $\Om$ has been defined in \eqref{omega} and where 
$\delta (a) \equiv a_{(-1)} \tp a_{(0)} \tp a_{(1)}$. 
\end{cor}
An element $\LL \in \G\tp \A$ satisfying $(\ep\tp\id)(\LL) = \eA$ and
(\ref{e314}/\ref{e315}) is called a {\it normal coherent (left
  diagonal) $\delta$--implementer} (with respect to $\gamma$), see [HN]. 

Note that by choosing $\A = \dG\reli\G$ and $\gamma = \id$, Cor.\ \ref{lem3.3}
implies that the multiplication on
$\dG\reli \G$ may uniquely be described by the relations of
one ``generating matrix'' $\LL = \sum e_\mu \tp (e^\mu \reli \e)$,
where $\{ e_\mu\}$ is a basis in $\G$ with dual basis $\{e^\mu\}$. In
fact this formulation is used in [HN] to prove the associativity of
the multiplication \eqref{e312}.

We now give a graphical interpretation of the identities
\eqref{e314},\eqref{e315} by using that any unital algebra map $\gamma
: \G\pfeil \A$ defines a $\G$--module structure on $\A$  
via $b \cdot A \equiv \pi_\A (b) A := \gamma (b) A$. Moreover,
given two $\G$--modules $V,W$ then $(V\tp\A)\tp W$ becomes a
$\G$--module by setting $\pi_{V\tp\A\tp W}(a) = (\pi_V\tp\pi_\A \tp
\pi_W)( \delta(a))$. Considering $\G$ as a $\G$--module by left
multiplication we now define the map 
\begin{equation*}
L^\A_{(\G\tp\A)\tp\G^*} : (\G\tp\A)\tp\G^* \pfeil \A,\quad 
 (b\tp A\tp \vi) \mapsto (\vi \tp \id)\big(\LL \cdot (b\tp A) \big)  
\end{equation*}
Then
Eq.~\eqref{e314} is equivalent to 
$L^\A_{(\G\tp\A)\tp\G^*}$ being an intertwiner of $\G$--modules,
i.e.\  to $L^\A_{(\G\tp\A)\tp\G^*} \in \makebox{Hom}_\G
\big((\G\tp\A)\tp \G^*,  
  \G\big)$. We depict this intertwiner as 
  \begin{equation*}
\begin{picture}(12,6)
\put(2,3){\makebox(0,0){$L^\A_{(\G\tp\A)\tp\G^*} : =$}}
\put(5,1){\bildL}
\put(6,0.5){\makebox(0,0){$\G$}}
\put(8,0.5){\makebox(0,0){$\A$}}
\put(10,0.5){\makebox(0,0){$\G^*$}}
\put(8,5.5){\makebox(0,0){$\A$}}
\put(5.5,1.25){(}
\put(8.25,1.25){)}
\put(11.5,3){\makebox(0,0){,}}
\end{picture}
\end{equation*}
and call this a {\it $d$--fork ($\equiv$ down fork) graph}.
The ``coherence condition'' \eqref{e315} is now equivalent to the
graphical identity 
\begin{equation}
  \label{graf31a}
\begin{picture}(25,11.5)
 \put(1,1.5){\makebox(0,0){$\G$}}
 \put(3,1.5){\makebox(0,0){$\G$}}
 \put(5,1.5){\makebox(0,0){$\A$}}
 \put(7,1.5){\makebox(0,0){$\G^*$}}
 \put(9,1.5){\makebox(0,0){$\G^*$}}
 \put(2,2){\bildL}
\put(0.5,2.25){\{}
 \put(2.25,2.25){[(} \put(5.25,2.25){)} \put(7.25,2.25){]\}} 
 \put(0,6){\bildgrossL}
 \put(0.5,6.25){(} \put(5.25,6.25){)}
 \multiput(1,2)(8,0){2}{\line(0,1){4}}
 \put(5,10.5){\makebox(0,0){$\A$}}
  \put(12,5){\makebox(0,0){$=$}}
\put(15,0.5){\makebox(0,0){$\G$}}
 \put(17,0.5){\makebox(0,0){$\G$}}
 \put(19,0.5){\makebox(0,0){$\A$}}
 \put(21,0.5){\makebox(0,0){$\G^*$}}
 \put(23,0.5){\makebox(0,0){$\G^*$}}
 \put(14,1){\bildfuenf}
 \put(14.25,3.5){[(} \put(17.25,3.5){)} \put(19.25,3.5){]}
 \put(20.5,3.5){(} \put(23.25,3.5){)}
  \put(14.5,1.25){\{} \put(16.25,1.25){[(} \put(19.25,1.25){)}
  \put(21.25,1.25){]\}}
 \put(14,4){\bildid[$\id_{\G\tp\G}$]}
 \put(20,4){\bildid[$f$]}
 \put(19,4){\line(0,1){3}}
 \put(14,6.25){\bildgrossLb}
  \put(15.5,6.5){(} \put(19.25,6.5){)}
  \put(19,10.75){\makebox(0,0){$\A$}}
\end{picture}
\end{equation}
Note that the lowest box on the r.h.s.\ represents the rebracketing
morphism $\Phi^{-1}$ defined in \eqref{Phi}.
This explains why one has to chose the complicated multiplication rule
\eqref{e312} instead of \eqref{e312old} if $\phi$ and therefore $f$
and $\Phi$ are non--trivial.

\subsection{Coherent $\cop$--flip operators}
We are now going to provide another set of generators in $\D(\G)$ which
later will be more appropriate for defining
the (quasitriangular) quasi--Hopf structure. 
Associated with any coherent (left diagonal) $\delta$--implementer
$\LL$ we define the element $\TT \in \G\tp\A$ by
\begin{equation}
  \label{L}
  \TT := [\Si(p_\rho^2) \tp\e]\cdot \LL \cdot (\id\tp
  \gamma)(\cop(p_\rho^1)), 
\end{equation}
where $p_\rho \equiv p_\rho^1\tp p_\rho^2$ has been given in
\eqref{defpq}. 
\begin{prop}
\label{prop31.1} {\rm [HN]} 
  The relation \eqref{L} defines a one-to-one correspondence between
  elements $\LL \in 
  \G\tp\A$ satisfying \eqref{e314} and \eqref{e315} and elements
  $\TT \in \G\tp\A$ satisfying
  \begin{align}
    \label{e316}
    (\id\tp \gamma)(\cop^{op}(a)) \, \TT  &= \TT (\id\tp \gamma)(\cop(a)), \quad
    \forall a \in \G \\ 
   \label{e317}
    \phi_{\A}^{312} \, \TT[13] (\phi^{-1})^{132}_{\A} \, \TT[23] \phi_{\A} &=
    (\cop\tp\id)(\TT), 
  \end{align} 
 where $\phi_\A = (\id\tp\id\tp \gamma)(\phi)$. 
$\LL$ is recovered from $\TT$ by
\begin{equation}
  \label{T}
  \LL = (\id\tp \gamma)(q^{op}_\rho)\, \TT
\end{equation}
Moreover $(\ep\tp\id)(\TT)
 = \eA$ if and only if $(\ep\tp\id)(\LL) = \eA$. 
\end{prop}
Following [HN,Sect.\ 11] we call the elements $\TT \in
\G\tp\A$ satisfying \eqref{e316} and \eqref{e317} {\it coherent
  $\cop$--flip operators}. They are special versions of {\it
  coherent $\la\rho$--intertwiners} [HN, Def.\ 10.8] associated with
quasi--commuting pairs $(\la,\rho)$ of left $\G$--coactions $\la$ and
right $\G$--coactions $\rho$ on an algebra $\M$.
Proposition \ref{prop31.1} has been proven algebraically in
[HN,Prop.\ 10.10]. Before giving an alternative proof below, using the
graphical calculus 
developed in Section 2.2 and 2.3, let us state the following central
consequence 
\begin{thm}
  \label{thm31}
Define the element $\DD \in \G\tp (\dG\tp\G)$ by\footnote{as before
  $\{e_\mu\}$ denotes a basis of $\G$ with dual basis $\{ e^\mu\}$}. 
\begin{equation}
  \label{D1}
   \DD : = \sum_\mu \Si(p_\rho^2)\, e_\mu \, p_{\rho (1)}^1 \tp (e^\mu \tp
     p_{\rho (2)}^1)
\end{equation}
and denote $i_D: \G \hookrightarrow \dG\tp\G$ the canonical
embedding $i_D (a) : = \edG \tp a$.
Then there is a unique algebra structure on the vector space $\dG\tp \G$ 
  satisfying
   \begin{align}
     & i_D(a) i_D(b) = i_D(ab) \qquad \forall a,b \in \G \label{eq31}\\
     & \DD \cdot (\id\tp i_D)(\cop(a)) = (\id\tp i_D)(\cop^{op}(a))
     \cdot \DD \qquad
     \forall a \in \G \label{eq32}  \\
     & \phi^{312} \, \DD[13] (\phi^{-1})^{132} \, \DD[23]
     \phi =     (\cop\tp\id)(\DD) \label{eq33}.
   \end{align}
 where we have identified $\phi \equiv (\id\tp\id\tp i_D)(\phi) \in
 \G\tp\G\tp\D(\G)$. This algebra is precisely the quantum double
 $\D(\G)$ defined in Prop.~\ref{prop311} and we have 
 \begin{equation}
   \label{eq:3.15}
   \vi\reli a = (i_D \tp \vi_{(1)})(q_\rho)\, (\vi_{(2)}\tp
   \id)(\DD)\, i_D (a)
 \end{equation}
\end{thm} 
\begin{proof}
 Follows from Proposition \ref{prop31.1} and Cor.~\ref{lem3.3} by choosing $\A
  = \D(\G)$ which means that $\TT = \DD$.
\end{proof}
We call $\DD$ the {\it universal $\cop$--flip operator} in $\D(\G)$.
The description of the quantum double $\D(\G)$ as given in Theorem
\ref{thm31} will be 
used in the next section to derive the quasi--Hopf structure of
$\D(\G)$. 
We now prove Proposition \ref{prop31.1}.
\begin{proof}[Proof of Proposition \ref{prop31.1}]
The equivalence $(\ep\tp\id)(\TT) = \eA \Leftrightarrow
(\ep\tp\id)(\LL) = \eA$ follows from property \eqref{eq15} of $\phi$. 
To show that the relations \eqref{e314} and \eqref{e315} for $\LL$ are
equivalent to the relations \eqref{e316} and \eqref{e317} for $\TT$,
respectively, we use the graphical calculus.
First we use the isomorphism
$\Psi^{\A}_{(\G\tp\A),\G^*}$ of Lemma \ref{umbieg}, to define the
intertwiner $T_{\G\A} \in 
Hom_\G (\G\tp\A,\A\tp\G)$ by 
\begin{equation}
  \label{graf31b}
  \begin{picture}(20,8.5)
   \put(0,4){\makebox(0,0){$T_{\G\A}\,\equiv$}}
   \put(3,2.5){\makebox(0,0){$\G$}}
   \put(5,2.5){\makebox(0,0){$\A$}}
   \put(3,3){\Dmatrix}
   \put(5,5.5){\makebox(0,0){$\G$}}
   \put(3,5.5){\makebox(0,0){$\A$}}
 \put(6.5,4){\makebox(0,0){$:=$}}
   \put(9,0.5){\makebox(0,0){$\G$}}
   \put(11,0.5){\makebox(0,0){$\A$}}
   \put(8,1){\boxviii}
   \put(13,1){\ukreis}
   \put(8,3.25){\bildL}
   \put(15,3.25){\line(0,1){4}}  
   \put(15,7.75){\makebox(0,0){$\G$}}
   \put(11,7.75){\makebox(0,0){$\A$}}
  \end{picture}
\end{equation}
Algebraically  Definition \eqref{graf31b} translate into  
$T_{\G\A} (b\tp A)  :=  \TT[21] \cdot (A\tp b)$, where $\TT \in
\G\tp\A$ is expressed in terms of $\LL$ by \eqref{L}. 
Now note that the property of $T_{\G\A}$ being an intertwiner of
$\G$--modules is equivalent to $\TT$ satisfying \eqref{e316}. Thus we
have proven the equivalence $\eqref{e314} \Leftrightarrow \eqref{e316}$
and since the map $\Psi^{\A}_{(\G\tp\A),\G^*}$ is invertible also the
invertibility of the transformation \eqref{L}. In fact, \eqref{T}
is equivalent to 
\begin{equation*}
\begin{picture}(15,7)
\put(0,1.5){\bildL}
\put(1,1){\makebox(0,0){$\G$}}
\put(3,1){\makebox(0,0){$\A$}}
\put(5,1){\makebox(0,0){$\G^*$}}
\put(3,6){\makebox(0,0){$\A$}}
\put(0.5,1.75){(}
\put(3.25,1.75){)}
\put(7.5,3.5){\makebox(0,0){$=$}}
\put(10,0.5){\makebox(0,0){$\G$}}
\put(12,0.5){\makebox(0,0){$\A$}}
\put(14,0.5){\makebox(0,0){$\G^*$}}
\put(10,1){\Dmatrix}
\put(14,1){\line(0,1){2}}
\put(9,3){\boxdi}
\put(12,6){\okreis}
\put(10,6.5){\makebox(0,0){$\A$}}
\end{picture}
\end{equation*}
We are left to show that $\eqref{e315} \Leftrightarrow \eqref{e317}$
(under the conditions \eqref{e314} and \eqref{e316}, respectively). To
this end we use that  Eq.\ \eqref{e317} is
graphically expressed as
{\rm  
  \begin{equation}
      \label{graf31c}
    \begin{picture}(18,16)
    \put(1,3.5){\makebox(0,0){$\G$}}
      \put(3,3.5){\makebox(0,0){$\G$}}
    \put(5,3.5){\makebox(0,0){$\A$}}
    \put(1,4){\line(0,1){2}}
   \put(3,4){\Dmatrix}
   \put(0,6){\boxdii}
    \put(5,9){\line(0,1){2}}
    \put(1,9){\Dmatrix}
    \put(5,11.5){\makebox(0,0){$\G$}}
      \put(3,11.5){\makebox(0,0){$\G$}} 
    \put(1,11.5){\makebox(0,0){$\A$}}
  \put(8,7.5){\makebox(0,0){$=$}}
  \put(10,1){\boxdii}
\put(11,0.5){\makebox(0,0){$\G$}}
      \put(13,0.5){\makebox(0,0){$\G$}}
    \put(15,0.5){\makebox(0,0){$\A$}}
  \put(10,3.5){\bildid[$\id_{\G\tp\G}$]}
  \put(15,3.5){\line(0,1){3}}
  \put(12,6.5){\Dmatrixgross}
  \put(12,9.5){\line(0,1){3}}
  \put(13,9.5){\bildidu[$\id_{\G\tp\G}$]}
   \put(11,12){\boxdii}
\put(16,15.5){\makebox(0,0){$\G$}}
      \put(14,15.5){\makebox(0,0){$\G$}}
    \put(12,15.5){\makebox(0,0){$\A$}}
    \end{picture}   
  \end{equation}}
Thus the following Lemma proves the equivalence
$\eqref{e315} \Leftrightarrow \eqref{e317}$ and therefore completes
the proof of Proposition \ref{prop31.1}
\end{proof} 
\begin{lem}
  \label{lem31.2}
  The graphical identities \eqref{graf31a} for $L^{\A}_{(\A\tp\G)\tp
    \G^*}$ and \eqref{graf31c} 
  for $T_{\G\A}$ are equivalent. 
\end{lem}
\begin{proof}
Let us prove \eqref{graf31a}
$\Rightarrow$ \eqref{graf31c}: 
  Using the definition \eqref{graf31b} we get for the l.h.s.\ of
  \eqref{graf31c} 
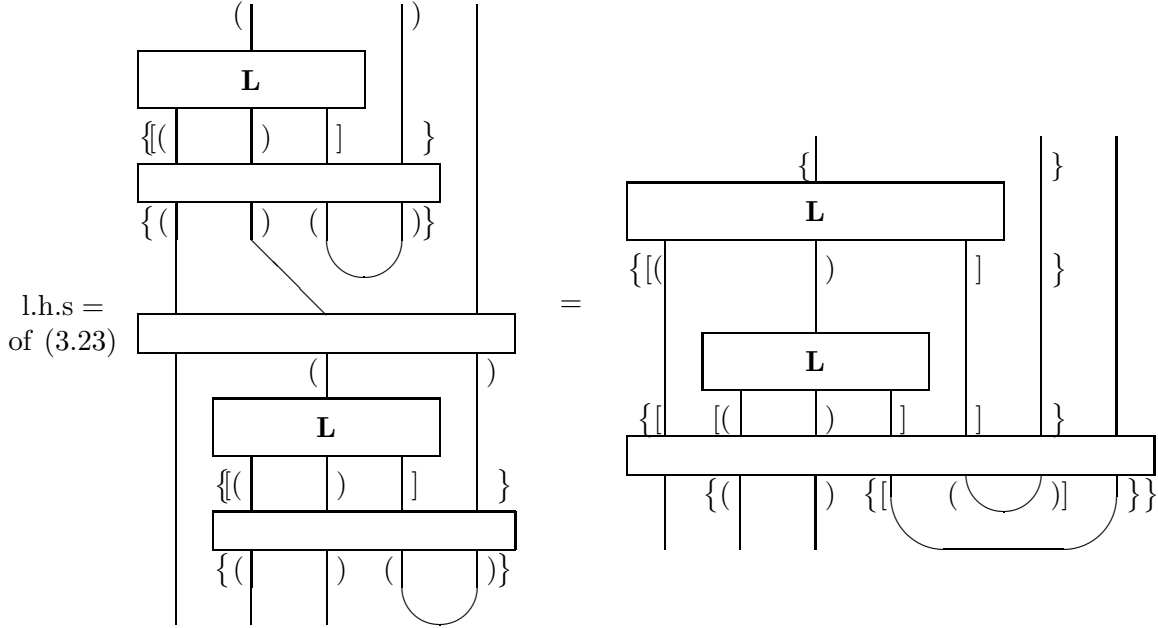
  \begin{equation*}
    \begin{picture}(34,16.5)
  \put(1,8.5){\makebox(0,0){$\text{l.h.s} =$}}
  \put(1,7.5){\makebox(0,0){of \eqref{graf31c}}}
  \multiput(4,0)(2,0){3}{\line(0,1){1}}    
  \put(5,1.25){\big\{}
\put(12.5,1.25){\big\}}
  \put(5,3.5){\big\{}
\put(12.5,3.5){\big\}}
  \put(4,1){\line(0,1){6.25}}
\put(4,8.25){\line(0,1){2}}
  \put(5,1){\boxviii}
  \put(10,1){\ukreis}
  \put(5,3.25){\bildL}
  \put(7.5,6.5){(}
\put(12.25,6.5){)}
\put(3,7.25){\framebox(10,1){}}
  \put(8,8.25){\line(-1,1){2}}
  \put(3,10.25){\boxviii}
  \put(3,10.5){\big\{}
\put(10.5,10.5){\big\}}
  \put(3,12.75){\big\{}
\put(10.5,12.75){\big\}}
  \put(8,10.25){\ukreis}
  \put(3,12.5){\bildL}
  \put(5.5,16){(}
\put(10.25,16){)}
  \put(10,12.5){\line(0,1){4}}
  \put(12,4){\line(0,1){3.25}}
  \put(12,8.25){\line(0,1){8.25}}
\put(14.5,8.5){\makebox(0,0){$=$}}
\put(16,4){\framebox(14,1){}}
\multiput(17,2)(2,0){3}{\line(0,1){2}}
\put(24.5,3.25){(}  \put(27.25,3.25){)]}  
 \put(29.25,3.25){\big\}\big\}}    
\put(22.25,3.25){\big\{[}
    \put(18,3.25){\big\{(} \put(21.25,3.25){)}      
\put(23,4){\gukreis}
\put(25,4){\ukreis}
 \put(18,5){\bildL}
\put(16.25,5.25){\big\{[}
 \put(18.25,5.25){[(} \put(21.25,5.25){)} \put(23.25,5.25){]}
 \put(25.25,5.25){]} \put(27.25,5.25){\big\}}
 \put(16,9){\bildgrossL}
  \put(20.5,12){\{} \put(27.25,12){\}}
 \put(16,9.25){\big\{[(} \put(21.25,9.25){)} \put(25.25,9.25){]}
\put(27.25,9.25){\big\}}
 \multiput(17,5)(8,0){2}{\line(0,1){4}}
 \multiput(27,5)(2,0){2}{\line(0,1){8}}
    \end{picture}
  \end{equation*}
Here we have used a pull through rule to push the lower $d$--fork up
and then we have combined all rebracketing morphisms in one
box.  Now
plugging in Eq.\ \eqref{graf31a}, splitting the rebracketing
morphism at the bottom into four factors and  pushing two of them up,
 one obtains \\ \\
 l.h.s. of  \eqref{graf31c} = 
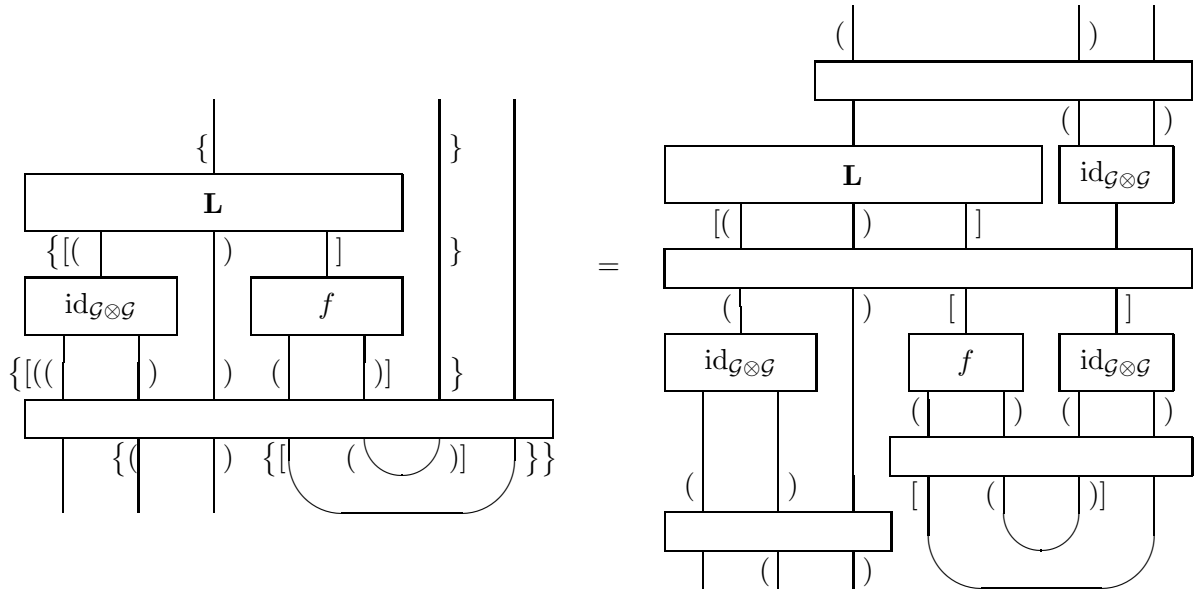
\begin{equation*}
  \begin{picture}(34,16)
\put(0,4){\framebox(14,1){}}
\multiput(1,2)(2,0){3}{\line(0,1){2}}
\put(8.5,3.25){(}  \put(11.25,3.25){)]}  \put(13.25,3.25){\big\}\big\}}  
\put(6.25,3.25){\big\{[}
   \put(2.25,3.25){\big\{(} \put(5.25,3.25){)}      
\put(7,4){\gukreis}
\put(9,4){\ukreis}
\put(-0.5,5.5){\big\{[((} \put(3.25,5.5){)}
\put(5.25,5.5){)}
 \put(6.5,5.5){(} \put(9.25,5.5){)]}
  \put(11.25,5.5){\big\}}
 \put(0,6){\bildid[$\id_{\G\tp\G}$]}
 \put(6,6){\bildid[$f$]}
 \put(5,6){\line(0,1){3}}
 \put(0.5,8.75){\big\{[(} \put(5.25,8.75){)}
 \put(8.25,8.75){]} \put(11.25,8.75){\}}
 \put(0,8.25){\bildgrossLb}
  \put(4.5,11.5){\{} \put(11.25,11.5){\}}
 \put(5,12.25){\line(0,1){0.75}}
 \multiput(11,5)(2,0){2}{\line(0,1){8}}
 \multiput(1,5)(2,0){5}{\line(0,1){1}}
  \put(15.5,8.5){\makebox(0,0){$=$}}
\multiput(18,3)(2,0){3}{\line(0,1){2}}
 \put(17,0){\boxdii}
\put(23,2){\boxvvii}
\put(25,6.75){\line(0,1){1.25}}
    \put(29,6.75){\line(0,1){1.25}}
    \put(23.5,5.25){\framebox(3,1.5){$f$}}
    \put(27.5,5.25){\framebox(3,1.5){$\id_{\G\tp \G}$}}
    \put(24,4){\line(0,1){1.25}}
    \put(26,4){\line(0,1){1.25}}
    \put(28,4){\line(0,1){1.25}}
    \put(30,4){\line(0,1){1.25}}
\put(17,4.5){\bildid[$\id_{\G\tp\G}$]}
\multiput(19,7.5)(6,0){2}{\line(0,1){0.5}}
\put(24,2){\gukreis}
\put(26,2){\ukreis}
\put(22,5){\line(0,1){3}}
\put(18.5,7.25){(} \put(22.25,7.25){)} 
\put(24.5,7.25){[}  \put(29.25,7.25){]}
\put(17,8){\framebox(14,1){}}
\put(18.25,9.5){[(} \put(22.25,9.5){)} 
\put(25.25,9.5){]} 
    \put(29,9){\line(0,1){1.25}}
 \put(27.5,10.25){\framebox(3,1.5){$\id_{\G\tp \G}$}}
 \put(28,11.75){\line(0,1){1.25}}
    \put(30,11.75){\line(0,1){1.25}}
\put(17,9){\bildgrossLb}
\put(27.5,12.25){(} \put(30.25,12.25){)}
\put(21,13){\framebox(10,1){}}
\put(21.5,14.5){(} \put(28.25,14.5){)}
\put(22,14){\line(0,1){1.5}}
\multiput(28,14)(2,0){2}{\line(0,1){1.5}}
  \end{picture}
\end{equation*}
Using the identity \eqref{grafdelta} the last picture equals the
r.h.s.\ of \eqref{graf31c}. Hence we have shown \eqref{graf31a}
$\Rightarrow$ \eqref{graf31c}. The implication \eqref{graf31c}
$\Rightarrow$ \eqref{graf31a} is shown similarly by bending the two upper
$\G$--legs in \eqref{graf31c} down again. Thus we have proved Lemma
\ref{lem31.2} and therefore also Proposition \ref{prop31.1}.
\end{proof}

%
\subsection{Left and right diagonal crossed products}
In this  subsection we sketch how the quantum double $\D(\G)$ may
equivalently be modeled on $\G\tp\dG$ instead of $\dG\tp\G$. (In fact
this is true for any diagonal crossed product as has been shown
in [HN,Thm.\ 10.2].) 
With the notation as in Proposition \ref{prop311} and with $\Om_R \in
\G^{\tp^5}$ given by  
$ \Om_R  := (h^{-1})^{21}\cdot (\Si\tp\Si\tp\idG^{\tp^3})(\Phi)$ the
right diagonal crossed product $\G\reli\dG$ is defined to be the
vector space $\G\tp\dG$ with multiplication rule 
\begin{equation}
   \label{e318}
  (a\reli \vi)( b \reli\psi) := 
   \Big[ a\, \big(\vi_{(1)}\re b \li \hat{S}^{-1}(\vi_{(3)})\big)\,\Om_R^3 
       \Big] \reli  \Big[(\Om_R^2  \arr\vi_{(2)}\arl
       \Om_R^4)(\Om_R^1\arr\psi\arl\Om_R^5)  \Big]. 
\end{equation}
This makes $\G\reli\dG$ an associative algebra with unit
$\e\tp\hat{\e}$, containing $\G \equiv \G\reli \hat{\e}$ as a unital
subalgebra. To see that the two algebras $\dG\reli \G$ and $\G\reli\dG$
are isomorphic let us begin with stating the
analogue of Lemma \ref{lem3.3}: Let $\gamma : \G\pfeil \A$ be
a unital algebra map into some target algebra $\A$. Then the relation  
\begin{equation}
  \label{e313'}
 \gamma_R
 (a\reli \vi) = \gamma(a) \cdot (\vi\tp\id)(\RR)
\end{equation}
provides a one--to--one correspondence between unital algebra
morphisms $\gamma_R : \G\reli \dG \pfeil \A$ extending $\gamma$
and elements $\RR \in \G\tp\A$ satisfying $(\ep\tp\id)(\RR) = \eA$ and
\begin{align}
  \label{e314'}
\RR\, [\eG\tp\gamma(a)] & = [ a_{(1)}\tp\gamma( a_{(0)})]\, \RR\,
[\Si( a_{(-1)}) \tp\eA], \quad \forall a \in \G \\
  \label{e315'}
\RR[13]\,\RR[23] & = [\Om_R^4 \tp \Om^5_R\tp\gamma(\Om^3_R)]\, (\cop\tp\id)(\RR)
\, [\Om^2_R \tp\Om_R^1\tp \eA)]
\end{align}
We call such elements {\it normal coherent right diagonal
  $\delta$--implementers} [HN].
With this definition one gets
\begin{lem}
  \label{prop31.3}
Let $\gamma : \G \pfeil \A$ be some unital algebra map. 
Then the relation
\begin{equation}
  \label{e320}
\RR : = [\Phi^5  \Si(\Phi^4 \beta) \tp \eA]\, \LL \, 
         [\Phi^2 \Si(\Phi^1 \beta) \tp \gamma(\Phi^3)] 
\end{equation}
defines a one--to--one correspondence between unital algebra maps 
$\gamma_L : \dG\reli \G \pfeil \A$  and
unital algebra maps $\gamma_R: \G\reli \dG \pfeil \A$
extending $\gamma$, as defined in \eqref{e313} and \eqref{e313'},
respectively. 
\end{lem}
\begin{proof}
  We will sketch the proof, using graphical methods. For more details
  see [HN, Prop.\ 10.5]. Defining the map $R^{(\G^*\tp\A)\tp\G}_\A:
  \A\pfeil (\G^*\tp\A)\tp\G$ by 
  \begin{equation*}
   R^{(\G^*\tp\A)\tp\G}_\A (A) : = \sum e^\mu \tp \big[
   \RR[21] \cdot (A\tp e_\mu)\big], \quad A\in  \A,
  \end{equation*}
property \eqref{e314'} of $\RR$ is equivalent to
$R^{(\G^*\tp\A)\tp\G}_\A$ being an intertwiner of
$\G$--modules. Depicting this intertwiner as  a {u--fork ($\equiv$ up
  fork) graph} 
\begin{equation*}
\begin{picture}(12,6)
\put(2,3){\makebox(0,0){$R^{(\G^*\tp\A)\tp\G}_\A  : =$}}
\put(5,1){\bildR}
\put(6,5.5){\makebox(0,0){$\G^*$}}
\put(8,5.5){\makebox(0,0){$\A$}}
\put(10,5.5){\makebox(0,0){$\G$}}
\put(8,0.5){\makebox(0,0){$\A$}}
\put(5.5,4.25){(}
\put(8.25,4.25){)}
\put(11.5,3){\makebox(0,0){,}}
\end{picture}
\end{equation*}
the relation \eqref{e320} may graphically be expressed as 
\begin{equation}
  \label{e321}
  \begin{picture}(21,8)
  \put(0,2){\bildR}
\put(1,6.5){\makebox(0,0){$\G^*$}}
\put(3,6.5){\makebox(0,0){$\A$}}
\put(5,6.5){\makebox(0,0){$\G$}}
\put(3,1.5){\makebox(0,0){$\A$}}
\put(0.5,5.25){(}
\put(3.25,5.25){)}
\put(7.5,4){\makebox(0,0){:=}}
\put(9.25,1.25){[(} \put(12.25,1.25){)} \put(14.25,1.25){]}
\put(15.5,1.25){(} \put(18.25,1.25){)}
\put(9,1){\bildfuenf}
\put(9.25,3.5){\{}  \put(11.5,3.5){[(} \put(14.25,3.5){)} 
\put(16.25,3.5){]\}}  
\put(10,1){\ukreis}
\put(16,1){\ukreis}
\put(14,0.5){\makebox(0,0){$\A$}}
\put(11,3){\bildL}
\multiput(10,4)(8,0){2}{\line(0,1){3}}
\put(10,7.5){\makebox(0,0){$\G^*$}}
\put(14,7.5){\makebox(0,0){$\A$}}
\put(18,7.5){\makebox(0,0){$\G$}}
  \end{picture}
\end{equation}
Since the r.h.s. of \eqref{e321} defines a $\G$--module intertwiner 
if and only if $\LL$ satisfies \eqref{e314}, the element $\RR$ defined
by \eqref{e320},\eqref{e321} satisfies \eqref{e314'} if and only if
$\LL$ satisfies 
\eqref{e314}. The equivalence of the coherence conditions \eqref{e315}
and \eqref{e315'} is 
shown by first expressing \eqref{e315'} as the graphical identity
\eqref{graf31a}, then plugging in the definition
\eqref{e321} and using a pull through rule to collect all rebracketing
morphisms at the bottom of the graph and finally using the identities
\eqref{grafdelta}. We leave it to the reader to draw the corresponding
pictures. 

We are left to show that relation
\eqref{e320} may be inverted. This follows by a ``two--sided version''
of Lemma \ref{umbieg}. The reader is invited to check that the
inverse is given by 
\begin{equation}
\label{e322}
\LL: =  [\Si(\alpha \bar{\Phi}^5)\bar{\Phi}^4 \tp \gamma(\bar{\Phi}^3)]
\, \RR \, [\Si(\alpha \bar{\Phi}^2)\bar{\Phi}^1 \tp \eA],
\end{equation}
with 
$\Phi^{-1} =: \bar{\Phi}^1 \tp \bar{\Phi}^2 \tp \bar{\Phi}^3 \tp
\bar{\Phi}^4 \tp \bar{\Phi}^5$.
Eq.\ \eqref{e322} may be expressed graphically as the up--side--down
mirror image 
of picture \eqref{e321} with $\LL$ and $\RR$ as well as $\G$ and
$\G^*$ exchanged. 
\end{proof}
Putting $\A = \dG\reli \G$ and $\gamma_L = \id$ or $\A = \G\reli\dG$
and $\gamma_R = \id$, respectively, Lemma~\ref{prop31.3} implies 
\begin{cor}
 \label{cor3.8}
  Let $(\delta,\Phi)$ be the two--sided coaction of $\G$ on $\G$ given
  in \eqref{Phi} and
  define $V\equiv V^1 \tp V^2 \tp V^3,\,\, W\equiv W^1 \tp W^2 \tp
  W^3\, \in \G^{\tp^3}$ by 
  \begin{equation*}
  V:= S(\bar{\Phi}^1)\alpha\bar{\Phi}^2 \tp \bar{\Phi}^3 \tp
   \Si(\alpha \bar{\Phi}^5) \bar{\Phi}^4,\quad \,
    W:= \Phi^2 \Si(\Phi^1 \beta) \tp \Phi^3 \tp \Phi^4 \beta S(\Phi^5) 
  \end{equation*}
Then the map 
\begin{equation*}
  \dG\reli \G \ni (\vi\reli a) \mapsto \Big(V^2 \reli \big(\Si(V^1) \arr
  \vi\arl V^3\big)\Big) \cdot (a\reli \hat{\e}) \, \in \G\reli\dG
\end{equation*}
provides an algebra isomorphism with inverse given by 
\begin{equation*}
  \G\reli \dG \ni (a\reli \vi) \mapsto (\hat{\e} \reli a)\cdot
   \Big(\big(W^1 \arr
  \vi\arl \Si(W^3)\big) \reli W^2\Big)\, \in \dG\reli \G. 
\end{equation*}
\end{cor}
Corollary \ref{cor3.8} has been proven for general diagonal crossed
products in [HN, Thm.\ 10.2.iii] using the notation $V \equiv
q_\delta$ and $W \equiv p_\delta$. We also remark that one may equally
well use the two--sided $\G$--coaction $\delta' : = (\id\tp\cop)\circ \cop$
with reassociator $\Phi' : = [\e\tp
(\id\tp\cop\tp\id)(\phi^{-1})][\e\tp\e\tp\phi^{-1}][(\id\tp\id\tp
\delta')(\phi)] $ to construct another to versions of quantum doubles
$\dG\reli_{\delta'}\G$ and $\G\reli_{\delta'}\dG$. Since the two--sided
coactions $(\delta,\Phi)$ and $(\delta',\Phi')$ are twist equivalent
[HN, Prop.8.4], these constructions are also isomorphic to the
previous ones, i.e. all four diagonal crossed products define
equivalent extensions of $\G$ [HN, Prop. 10.6].


\subsection{The quasitriangular quasi--Hopf structure}
In [HN] we have shown that $\D(\G)$ is a quasi--bialgebra. As one
might expect, $\D(\G)$ is 
even a quasitriangular quasi-Hopf algebra. This is the content of the
next theorem. 
\begin{thm}
\label{thm31a}
Let $\D(\G)$ be the associative algebra defined in Theorem \ref{thm31}.
   Then $(\D(\G),\cop_D,\ep_D,\phi_D,S_D,\alpha_D,\beta_D, R_D)$ is a
   quasitriangular quasi-Hopf algebra, where 
 \begin{align}
    \label{eq37}
     \phi_D :=& (i_D\tp i_D\tp i_D)(\phi),\\
    \label{eq38}
    R_D :=& (i_D \tp\id)(\DD) = \sum_\mu i_D(e_\mu) \tp D(e^\mu)
  \end{align}
   and where the structural maps are given by 
   \begin{align}
     \label{eq34}
     \copd(i_D(a)) &:= (i_D\tp i_D)(\cop(a)),\,\,\quad \forall a\in\G \\
     (i_D\tp\copd)(\DD) &:= (\phi_D^{-1})^{231} \DD[13] \phi_D^{213} \DD[12]
         \phi_D^{-1},  \label{eq35} \\
    \label{eq36}
    \ep_D (i_D(a)) &:= \ep(a),\,\,\,  \forall a \in \G, \\
    \label{eq36b}
     (\id\tp\ep_D)(\DD) &:= (\ep\tp\id)(\DD) \equiv \e_{\D(\G)}.
  \end{align} 
 Furthermore the antipode $S_D$ is defined by 
 \begin{align}
   \label{eq39}
  & S_D(i_D(a)) : = i_D (S(a)), \,\quad \forall a\in \G \\
    \label{eq310}
  &  (S\tp S_D)(\DD) := 
  (\id\tp i_D)(f^{op})\,\DD(\id\tp i_D)\,(f^{-1}) 
 \end{align}
 where $f\in \G\tp\G$ is the twist defined in \eqref{eq118}. The
 elements $\alpha_D,\beta_D$ are given by
 \begin{equation}
   \label{eq39a}
   \alpha_D := i_D(\alpha),\,\,\,\,
   \beta_D := i_D(\beta).
 \end{equation}
\end{thm}
\begin{proof}
  To simplify the notation we will frequently suppress the
  embedding $i_D$, if no confusion is possible, i.e.\ we write
  $\alpha\equiv i_D(\alpha)=\alpha_D, \,\,(\id\tp\id\tp
  i_D)(\phi)\equiv \phi$ etc. 
  To show that (\ref{eq34}),(\ref{eq35}) define an algebra morphism
  $\copd \, : 
  \D(\G) \pfeil \D(\G)\tp \D(\G)$, it is sufficient to check the
  consistency with the defining relations (\ref{eq31}) -
  (\ref{eq33}). 
  Consistency with \eqref{eq31} is obvious because of \eqref{eq34}.
  Let us go on with (\ref{eq32}): For the r.h.s.\ we get 
  \begin{equation*}
     (\id\tp\copd)\big( \cop^{op}(a)\, \DD\big) =
     [(\cop\tp\id)(\cop(a))]^{231} \cdot 
     (\phi^{-1})^{231} \DD[13] \phi^{213} \DD[12] \phi^{-1}
  \end{equation*}
 The l.h.s.\ yields 
  \begin{equation*}
    (\id\tp\copd)\big(\DD\,\cop(a)\big) =
    (\phi^{-1})^{231} \DD[13] 
    \phi^{213} \DD[12] \phi^{-1} [(\id\tp\cop)(\cop(a))],
   \end{equation*}
  which, using (\ref{eq32}) and the property $\phi (\cop\tp\id)(\cop) =
  (\id\tp\cop)(\cop) \phi$ to shift the factor $(\id\tp\cop)(\cop(a))$ 
  to the left, equals the r.h.s..\\
  Consistency with the relation (\ref{eq33}) may be checked in a
  longer but analogues calculation, where one also has to use the
  pentagon equation for $\phi$ several times, as in the
  proof of [HN, Lem.\ 11.2]. Hence $\cop_D$ is an
  algebra map. To show that $\cop_D$ is quasi-coassociative we compute
  by a similar calculation 
  \begin{equation*}
    (\e\tp\phi_D) \cdot (\id\tp\cop_D\tp\id)\circ (\id\tp\cop_D)(\DD) =
    (\id\tp\id\tp\cop_D)\circ (\id\tp\cop_D) (\DD)\cdot (\e\tp\phi_D),
  \end{equation*}
  by using again the pentagon equation for $\phi$ and the covariance
  property (\ref{eq32}).

  The property of $\ep_D$
  being a counit for $\cop_D$ follows directly from the fact that 
  $(\id\tp\ep\tp\id)(\phi) = \e\tp\e$. Hence
  $(\D(\G),\cop_D,\ep_D,\phi_D)$ is a quasi--bialgebra, see [HN,
  Thm. 11.3] for more details.\\
  
  To show quasitriangularity we first note that
  the element $R_D =(i_D\tp\id)(\DD)$ fulfills (\ref{eq19}) and
  (\ref{eq110}) so 
  to say  by definition because of (\ref{eq33}) and (\ref{eq35}).
  The invertibility of $R_D$ is equivalent to the invertibility of
  $\DD$ which will be proved in Lemma \ref{leminverse} (i)  below.
  We are left to show that $R_D$ intertwines $\cop_D$ and
  $\cop_D^{op}$, i.e. 
  \begin{align}
    \label{eq311}
    \cop_D^{op}(i_D(a)) \cdot R_D &= R_D \cdot \cop_D(i_D(a)),\, \quad
     \forall a\in \G\\
    \label{eq311a}
    (\id\tp\cop^D_{op})(\DD)) \cdot R_D^{23} &= R_D^{23} \cdot
    (\id\tp\cop_D)(\DD).
  \end{align}
Now Eq.\ (\ref{eq311}) follows from (\ref{eq32}). Hence we also get 
in $\D(\G)^{\tp^3}$
\begin{equation}
     \label{eq312}
     R_D^{12} \cdot (\cop_D\tp\id)(R_D)  = (\cop_D^{op}\tp\id)(R_D)
     \cdot R_D^{12},
   \end{equation}
which together with (\ref{eq19}) implies the 
  quasi-Yang Baxter equation    
  \begin{equation}
   \label{eq312a}
    \Fii[321]{D}\,R_D^{12}\,\Fi[312]{D}\,R_D^{13}\,\Fii[132]{D}\, R_D^{23} =
     R_D^{23}\, \Fii[231]{D}\,R_D^{13}\,\Fi[213]{D}\,R_D^{12}\,\phi^{-1}_{D}.
  \end{equation}
Using the Definition (\ref{eq35}), Eq.\ (\ref{eq312a}) is further
equivalent to 
\begin{equation*}
    (i_D \tp\cop_D^{op})(\DD) \cdot R_D^{23} = R_D^{23} \cdot
    (i_D \tp\cop_D)(\DD)
\end{equation*}
which also proves (\ref{eq311a}). Hence $R_D$
   is quasitriangular. \\ 
  
  In order to prove that the definition of $S_D$ in
  (\ref{eq39}),(\ref{eq310}) may be extended anti-multiplicatively to the
  entire algebra $\D(\G)$, we have to show that this continuation is
  consistent with the defining relations (\ref{eq32}),(\ref{eq33}).
 This amounts to showing
  \begin{equation}
    \label{eq313}
   (S\tp S_D)(\DD) \cdot (S\tp S_D)(\cop^{op}(a)) =
       (S\tp S_D)(\cop(a)) \cdot (S\tp S_D)(\DD), \quad \text{and} 
   \end{equation}
   \begin{multline} 
    \label{eq314}
    (S\tp S\tp S_D)\big( (\cop\tp\id)(\DD)\big) = 
     (S\tp S\tp S_D)(\phi) \cdot (S\tp S\tp S_D)(\DD[23]) \\ 
    \cdot
      (S\tp S\tp S_D)(\Fii[132]{}) \cdot (S\tp S\tp S_D)(\DD[13])\cdot
     (S\tp S\tp S_D)(\phi^{312}).
    \end{multline}
 Since by definition $(S\tp S_D)(\DD) = f^{op}\DD f^{-1}$ \footnote
 {where we have 
   again suppressed the embedding $\id\tp i_D$ of $f$}, equation (\ref{eq313})
 follows directly from (\ref{eq32}) and the fact, that by (\ref{eq120}) $f$ has the
 property $f \cdot 
 \cop(S(a)) = (S\tp S) (\cop_{op}(a)) \cdot f$. For the proof of
 (\ref{eq314}) let us recall, that $\cop^f : = f\cop(\cdot) f^{-1}$
 defines a twist equivalent quasi-coassociative coproduct on $\G$ with
 twisted reassociator 
 $\phi_f$ defined in (\ref{eq114}) satisfying 
  $\phi_f = (S\tp S\tp S)(\phi^{321})$ (see (\ref{eq121})).
 Thus we get for the l.h.s.\ of (\ref{eq314}) (with $ \DD_f : = 
  f^{op} \DD f^{-1}$)
 \begin{align*}
   (S\tp S\tp S_D)\big( (\cop\tp\id)(\DD) \big) 
    & = (\cop^f_{op}\tp\id)\big((S\tp S_D)(\DD)\big) \\
    & = (\cop^f_{op}\tp\id)(\DD_f)  \\
    & = \Fi[321]{f} \, \DD[23]_f \, \Fii[231]{f} \DD[13]_f \Fi[213]{f},
 \end{align*}
  where the last equality is exactly the transformation property of
  a quasitriangular 
  R-matrix under a twist [Dr2] and may be proven analogously using
  (\ref{eq32}). By (\ref{eq121}) 
  this equals the r.h.s.\ of (\ref{eq314}). Hence $S_D$ defines an
  anti-algebra morphism on $\D(\G)$.  
  
  We are left to show that the map $S_D$ fulfills the antipode axioms
  given in (\ref{eq16}) and (\ref{eq17}). Axiom (\ref{eq17}) is
  clearly fulfilled since we have $S_D \circ i_D = i_D \circ S$ and 
  $\alpha_D = i_D(\alpha)$, $\beta_D = i_D(\beta)$, $\phi_D = (i_D\tp
  i_D \tp i_D)(\phi)$. Noting that $\cop_D(i_D(a)) = (i_D\tp
  i_D)(\cop(a)) \,,a \in \G$, the validity of axiom \eqref{eq16}
  follows from its validity in $\G$ and Lemma \ref{leminverse} (ii)
  below.
\end{proof}
\newpage
  \begin{lem}$\,$
    \label{leminverse}
  \begin{itemize}
   \item[(i)] The universal flip operator $\DD \,\in \G\tp \D(\G)$ is
     invertible where the inverse is given by 
  \begin{equation}
   \label{Di}
   \DD[-1] = 
   [X^j\beta S(P^iY^j)\tp\e ] \cdot [(S\tp\id)(q^{op}_\rho\,\DD)]
   \cdot [(R^i\tp Q^i)
   \,\cop_{op}(Z^j)], 
  \end{equation} 
where $q_\rho \in \G\tp\G$ has been defined in \eqref{defpq}.
   \item[(ii)]
   Let $\mu_D$ denote the multiplication map $\mu_D: \D(\G)\tp \D(\G)
   \pfeil \D(\G)$, then
     \begin{align}
       \label{eqS1}
      (\id \tp\mu_D)\circ(\id\tp S_D\tp\id)\Big( (\id\tp\cop_D)(\DD)
      \cdot (\e_\G\tp\e_\G\tp\alpha_D)\Big) &= \e_\G \tp\alpha_D \\
      \label{eqS2}
        (\id \tp\mu_D)\circ(\id\tp\id\tp S_D)\Big( (\id\tp\cop_D)(\DD)
      \cdot (\e_\G\tp\beta_D\tp\e_\G)\Big) &= \e_\G \tp\beta_D 
     \end{align}
   \end{itemize}
  \end{lem}
  \begin{proof}
  We will use the graphical methods adopted in Sections 2.2/2.3. To this
  end  let us view $\G$ and $\D \equiv \D(\G)$ as left $\G$--modules.
   Then, due to
  \eqref{eq32}, $\B_{\G \D} : = \tau_{\G \D}\circ \DD$ defines
  an intertwiner   $\B_{\G\D}: \, \G\tp\D \pfeil \D\tp\G$ which will
  be depicted as 
  \begin{equation*}
    \begin{picture}(6,4)
      \put(2,2){\makebox(0,0){$\B_{\G\D}=: $}}
      \put(4,0.5){\makebox(0,0){$\G$}}
      \put(6,0.5){\makebox(0,0){$\D$}}
      \put(4,3.5){\makebox(0,0){$\D$}}
      \put(6,3.5){\makebox(0,0){$\G$}}
      \put(4,1){\Dmatrix}
    \end{picture}
  \end{equation*}
(In fact this is the intertwiner $T_{\G\A}$ defined in 
\eqref{graf31b} for the special case $\A = \D$).
  For the left modules ${}^*\G$ and ${}^*\D$ the corresponding
  intertwiners $\B_{{}^*\G\D}$, $\B_{\G{}^*\D}$, $\B_{{}^*\G{}^*\D}$ 
are defined with the help of the map $S$ and/or
  $S_D$. Graphically they are represented by the same picture, except
  that the colours of the legs are replaced by ${}^*\G$ and(or)
  ${}^*\D$, respectively. The reason for distinguishing the $\D$--line
  from the $\G$--line lies in the fact that unlike in (\ref{graf3})
  $\B_{\G\D}$ is {\bf not} given in terms of a quasitriangular
  R--matrix in $\G\tp\G$, which is why we write $\B_{\G\D}$ in place
  of $B_{\G\D}$. Correspondingly, the identities derived in
  Section 2.3 are {\bf not} automatically valid for $\B_{\G\D}$. We now
  show, which of them still hold. First, since $S$ is an antipode for
  $\cop$,  Eq. (\ref{eq33}) together with $(\ep\tp\id)(\DD) = \e_\D$ implies 
  the equality (compare with (\ref{graf9})) 
  \begin{equation}
  \label{grafS1}
  \begin{picture}(15,12)
  \put(1,0.5){\makebox(0,0){${}^*\G$}}
  \put(3,0.5){\makebox(0,0){$\G$}}
  \put(5,0.5){\makebox(0,0){$\D$}}
  \put(3,1){\Dmatrix}
  \put(1,1){\line(0,1){2}}
  \put(0,3){\boxdii}
  \put(1,6){\Dmatrix}
  \put(5,6){\line(0,1){2}}
  \put(0,8){\boxdi}
  \put(3,11){\okreis}
  \put(1,11.5){\makebox(0,0){$\D$}}
  \put(8,5){\makebox(0,0){$=$}}
  \put(9,4){\boxdii}
  \put(10,7){\okreis}
  \put(14,7){\line(0,1){1}}
    \multiput(10,3)(2,0){3}{\line(0,1){1}}
  \put(10,2.5){\makebox(0,0){${}^*\G$}}
  \put(12,2.5){\makebox(0,0){$\G$}}
  \put(14,2.5){\makebox(0,0){$\D$}}
  \end{picture}
\end{equation}
 and a step by step repetition of the
prove of \eqref{graf9c} yields
  \begin{equation}
    \label{grafS2}
     \begin{picture}(18,12)
      \put(1,6){\makebox(0,0){$(\B_{\G\D})^{-1} =: $}}
      \put(4,4.5){\makebox(0,0){$\D$}}
        \put(6,4.5){\makebox(0,0){$\G$}}
       \put(4,5){\iDmatrix}  
        \put(6,7.5){\makebox(0,0){$\D$}}
        \put(4,7.5){\makebox(0,0){$\G$}}
        \put(7.5,6){\makebox(0,0){=}}
        \put(14,0.5){\makebox(0,0){$\D$}}
        \put(16,0.5){\makebox(0,0){$\G$}}
        \put(9,1){\bilddD}
        \put(12,11.5){\makebox(0,0){$\D$}}
        \put(10,11.5){\makebox(0,0){$\G$}}        
    \end{picture}
   \end{equation}
This means that algebraically we get the analogue of the first
identity in Eq. (\ref{graf7a}) which yields \eqref{Di}.
Thus we have proven part (i) \\

To prove (ii) let us translate the two claims \eqref{eqS1} and \eqref{eqS2}
into the graphical language as
\begin{equation}
  \label{grafS3}
  \begin{picture}(16,12)(-3,0)
  \put(-3,5){\makebox(0,0){\eqref{eqS1} $\Leftrightarrow$}}
  \put(1,0.5){\makebox(0,0){$\G$}}
  \put(5,0.5){\makebox(0,0){$\D$}}
  \put(3,0.5){\makebox(0,0){${}^*\D$}}
  \put(1,1){\Dmatrix}
  \put(5,1){\line(0,1){2}}
  \put(0,3){\boxdi}
  \put(3,6){\Dmatrix}
  \put(1,6){\line(0,1){2}}
  \put(0,8){\boxdii}
  \put(1,11){\okreis}
  \put(5,11.5){\makebox(0,0){$\G$}}
  \put(8,5){\makebox(0,0){$=$}}
 \put(10,4){\boxdi}
  \put(13,7){\okreis}
  \put(11,7){\line(0,1){1}}
   \multiput(11,3)(2,0){3}{\line(0,1){1}} 
  \put(11,2.5){\makebox(0,0){$\G$}}
  \put(13,2.5){\makebox(0,0){${}^*\D$}}
  \put(15,2.5){\makebox(0,0){$\D$}}
  \end{picture}
\end{equation}
and 
\begin{equation}
  \label{grafS4}
  \begin{picture}(16,12)(-3,0)
\put(-3,5){\makebox(0,0){\eqref{eqS2} $\Leftrightarrow$}}
  \put(1,0.5){\makebox(0,0){$\G$}}
  \put(0,1){\boxdii}
  \put(3,1){\ukreis} 
  \put(1,4){\Dmatrix}
  \put(5,4){\line(0,1){2}}
  \put(0,6){\boxdi}
  \put(3,9){\Dmatrix}
  \put(1,9){\line(0,1){2}}
  \put(3,11.5){\makebox(0,0){${}^*\D$}}
  \put(1,11.5){\makebox(0,0){$\D$}} 
   \put(5,11.5){\makebox(0,0){$\G$}}
  \put(8,5){\makebox(0,0){$=$}}
 \put(10,4){\boxdi}
  \put(11,4){\ukreis}
  \put(15,3){\line(0,1){1}}
   \multiput(11,7)(2,0){3}{\line(0,1){1}} 
  \put(15,8.5){\makebox(0,0){$\G$}}
  \put(11,8.5){\makebox(0,0){$\D$}}
  \put(13,8.5){\makebox(0,0){${}^*\D$}}
  \end{picture}
\end{equation}
Note that as opposed to \eqref{grafS1} the identities \eqref{grafS3} and
\eqref{grafS4} are not automatically satisfied, since $S_D$ is not yet
proved to be an antipode for $\cop_D$.
To prove \eqref{grafS3} and \eqref{grafS4} we now proceed backwards
along the proof of Lemma \ref{lem02}, i.e. we use
Lemma \ref{umbieg} to show that either of these two identities is
equivalent to 
\begin{equation}
  \label{grafS5}
  \begin{picture}(13,12)
   \put(1,0.5){\makebox(0,0){${}^*\D$}}
  \put(3,0.5){\makebox(0,0){$\G$}}
  \put(0,1){\bildcD}
   \put(5,11.5){\makebox(0,0){$\G$}}
  \put(7,11.5){\makebox(0,0){${}^*\D$}}
  \put(9.5,6){\makebox(0,0){$=$}}
  \put(11,5){\iDmatrix}
   \put(11,4.5){\makebox(0,0){${}^*\D$}}
  \put(13,4.5){\makebox(0,0){$\G$}}
  \put(18.5,6){\makebox(0,0){$:=(\B_{\G{}^*\D})^{-1}$}}
  \end{picture}
\end{equation}
More precisely \eqref{grafS3} is equivalent to \eqref{grafS5} just as
\eqref{graf10} is equivalent to the second equation in \eqref{graf7},
and ``rotating'' this proof by $180^\circ$ in the drawing plane we
also get $\eqref{grafS4} \Leftrightarrow \eqref{grafS5}$. Thus we are
left with proving \eqref{grafS5}. To this end we remark, that
\eqref{grafS1} equally holds if we replace $\D$ by ${}^*\D$, and
therefore \eqref{grafS2} also holds with $\D$ replaced by
${}^*\D$. Hence \eqref{grafS5} follows from \eqref{grafS2} provided we
can show  
\begin{equation}
  \label{grafS6}
  \begin{picture}(24,12)
  \put(1,0.5){\makebox(0,0){${}^*\D$}}
  \put(3,0.5){\makebox(0,0){$\G$}}
  \put(0,1){\bildcD}
   \put(5,11.5){\makebox(0,0){$\G$}}
  \put(7,11.5){\makebox(0,0){${}^*\D$}}
  \put(9.5,6){\makebox(0,0){$=$}}
        \put(16,0.5){\makebox(0,0){${}^*\D$}}
        \put(18,0.5){\makebox(0,0){$\G$}}
        \put(11,1){\bilddD}
        \put(14,11.5){\makebox(0,0){${}^*\D$}}
        \put(12,11.5){\makebox(0,0){$\G$}}         
\end{picture}
\end{equation}
By Lemma \ref{corantipode} this further equivalent to 
\begin{equation}
    \begin{picture}(22,8)(0,-1)
       \put(1,-0.5){\makebox(0,0){${}^*\G$}}
       \put(3,-0.5){\makebox(0,0){${}^*\D$}}
       \put(5,-0.5){\makebox(0,0){$\G$}}
       \put(7,-0.5){\makebox(0,0){$\D$}}
       \put(1,0){\Dmatrix}
       \multiput(5,0)(2,0){2}{\line(0,2){2}}
       \put(0,2){\bilde}
       \put(10,4){=}
       \put(13,-0.5){\makebox(0,0){${}^*\G$}}
       \put(15,-0.5){\makebox(0,0){${}^*D$}}
       \put(17,-0.5){\makebox(0,0){$\G$}}
       \put(19,-0.5){\makebox(0,0){$\D$}}
       \put(17,0){\Dmatrix}
       \multiput(13,0)(2,0){2}{\line(0,2){2}}
       \put(12,2){\bilde}
    \end{picture} 
    \label{grafS7}
  \end{equation}
Using \eqref{graf6} and \eqref{graf6al}, Eq.\ \eqref{grafS7} is
algebraically equivalent to
\begin{equation*}
  \gamma^{op}\, \DD = (S\tp S_D)(\DD)\,\gamma,
\end{equation*}
which finally holds by  \eqref{eq310},\eqref{eq32}  
and \eqref{eq119b}. This concludes the proof of Lemma \ref{leminverse}
(ii) and therefore of Theorem \ref{thm31a}.
 \end{proof}
%
%
%
%

Clearly, if $\G$ is a Hopf algebra and $\phi = \e\tp\e\tp\e$, one
recovers the well-known
definitions of $\cop_D, S_D$ and $R_D$ in Drinfelds quantum double
\begin{align*}
  \cop_D(i_D(g)) &= (i_D\tp i_D)(\cop(g)) \\
  \cop_D(D(\vi)) &= (D\tp D)(\hat{\cop}^{op}(\vi)) \\
  S_D(i_D(g)) &= i_D(S(g)) \\
  S_D(D(\vi)) &= D(S^{-1}(\vi))\\
  R_D &= (\hat{\e} \tp e_\mu ) \tp (e^\mu \tp \e),
\end{align*}
where $D(\vi) : = (\vi\tp\id)(\DD),\,\vi \in \dG$.
As in the Hopf algebra case, one may take the construction of the
quasitriangular R-Matrix in $\D(\G)$ as the starting point and
formulate Theorem 
\ref{thm31} together with Theorem \ref{thm31a} differently:
\begin{cor}
  \label{cor31}
   Let $\G$ be a finite dimensional quasi-Hopf algebra with invertible
   antipode. Then there exists a unique 
   quasi-Hopf algebra $\D(\G)$ such that 
   \begin{enumerate}
   \item[(i)] $\D(\G) = \dG\tp\G$ as a vector space,
   \item[(ii)] the canonical embedding  $i_D : \G \hookrightarrow 
    \edG\tp\G \subset \D(\G)$ is
     a unital injective homomorphism of quasi-Hopf algebras,
   \item[(iii)]  
    Let $\DD \in \G\tp \D(\G)$ be given by Eq.\ \eqref{D1},
   then  the element 
      $R_D : = (i_D \tp \id)(\DD) \in \D(\G)\tp\D(\G)$
  is quasitriangular.    
   \end{enumerate}
This quasi--Hopf algebra structure is given by the definitions in
  Theorem \ref{thm31} and \ref{thm31a}.
\end{cor}
\begin{proof}    
The property (ii) implies (\ref{eq31}), (\ref{eq34}), \eqref{eq37},
\eqref{eq36} and \eqref{eq39a}, yielding also $f_D = (i_D \tp i_D)(f)$.
  The quasitriangularity of $R_D$ implies (\ref{eq32}), (\ref{eq33}),
  (\ref{eq35}) and \eqref{eq36b} and
  according to (\ref{cor01}) 
  $(S_D\tp S_D)(R_D) = f_D^{op} \,R_D f_D^{-1}$. Hence the antipode
  is uniquely fixed to 
  be the one defined in Theorem \ref{thm31a}.
\end{proof}
We are now in the position to prove Theorem A. 
\begin{proof}[Proof of Theorem A]
  First note that Corollary \ref{cor31} already proves the existence
  parts (i),(ii) of Theorem A by putting $D(\vi) : =
  (\vi\tp\id)(\DD)$. The fact that $\mu: \dG\tp\G \pfeil \D(\G)$
  provides a linear isomorphism follows from the last
  statement in Theorem \ref{thm31}. Moreover, if $\tilde{\D} \supset
  \G$ is another Hopf algebra extension and if $\tilde{D} : \dG \pfeil
  \tilde{\D}$ is a linear map such that $\tilde{\D}$ is algebraically
  generated by $\G$ and $\tilde{D}(\dG)$ and $R_{\tilde{D}}: = e_\mu
  \tp \tilde{D}(e^\mu) \in \tilde{\D}\tp \tilde{\D}$ is
  quasitriangular, then $\nu : \D(\G) \pfeil \tilde{\D}$ 
  \begin{equation}
    \label{eq:3.58'}
    \nu (\vi\tp a) : = (\id\tp \vi_{(1)})(q_\rho)\,
    \tilde{D}(\vi_{(2)}) \, a
  \end{equation} 
is a uniquely and well defined algebra map satisfying $\nu \circ D =
\tilde{D}$ by Prop.~\ref{prop31.1}. In fact, the quasitriangularity of
$R_{\tilde{D}}$ implies that $\nu$ is even a quasi--bialgebra
homomorphism. Thus $\D(\G)$ also solves the
universality property (iii) of Theorem A. In particular
the extension $\D(\G)\supset \G$ is unique up to equivalence.
\end{proof}


\subsection{The category Rep $\!\D(\G)$}  
We will now give a representation theoretical interpretation of the
quantum double $\D(\G)$ by describing its representation category in terms of
the representation category of the underlying quasi-Hopf algebra
$\G$. In this way we will show that $\D(\G)$ is a concrete
realization of the quantum double as defined by Majid in [M2] with the
help of a Tannaka-Krein-like reconstruction theorem. 
We denote the monoidal category of finite dimensional unital
representations of $\D(\G)$  and of $\G$ by $\rep\D(\G)$ and
$\rep\G$, respectively.  The next
proposition states a necessary and sufficient condition, under which a
representation of $\G$ extends to a representation of $\D(\G)$:
\begin{prop}$\,$\\
  \label{thm32}
  1.)   The objects of $\rep\D(\G)$ are in one to one correspondence with
   pairs $\{(\pi_V,V),\DD_V\}$, where $(\pi_V,V)$ is a finite
   dimensional representation of $\G$ and where $\DD_V \in \G\tp
   \End_\Co(V)$ is a normal coherent $\cop$-flip, i.e.  
   \begin{enumerate}
   \item[(i)] $(\ep\tp\id)(\DD_V) = \id_V$
   \item[(ii)] $\DD_V \cdot (\id\tp\pi_V)(\cop(a)) =
               (\id\tp\pi_V)(\cop^{op}(a)) \cdot \DD_V, \,\quad
               \forall a\in\G$
   \item[(iii)] $\phi_V^{312}\,\DD[13]_V\,(\phi_V^{-1})^{132}\,
     \DD[23]_V\,\phi_V = (\cop\tp\id)(\DD_V)$, where $\phi_V :=
     (\id\tp\id\tp\pi_V)(\phi)$. 
   \end{enumerate}
2.) Let $\{(\pi_V,V),\DD_V\}$ and $\{(\pi_W,W),\DD_W\}$ be as above,
then 
$$\makebox{\rm Hom}_{\D(\G)} = \big\{ t \in \makebox{\rm Hom}_\G(V,W) \mid
(\id\tp t)(\DD_V) = \DD_W  \big\}$$
\end{prop}
\begin{proof}
  We define the extended representation $\pi^D_V$ on the generators of
  $\D(\G)$ by
  \begin{align}
    \label{eq321}
     \pi_V^D (i_D(g)) : &= \pi_V(g), \quad g \in \G \\
     \pi_V^D (D(\vi)) : &= (\vi\tp\id_{End_V})(\DD_V), \quad \vi \in\dG
  \end{align}
    Condition (i) implies that $\pi_V^D$ is unital whereas conditions
    (ii),(iii) just reflect the defining relations (\ref{eq32})
    and (\ref{eq33}) of $\D(\G)$, which ensures, that $\pi_V^D$ is a
    well defined algebra morphism. On the other hand, given a
    representation $(\pi_V^D,V)$ of $\D(\G)$, we define  
    \begin{equation*}
     \DD_V := (\id_\G\tp\pi_V^D)(\DD)
    \end{equation*}
    which clearly satisfies conditions (i) - (iii). This proves part~1. Part~2. follows trivially.
\end{proof}
To get the relation with Majid's formalism [M2] we now write
 $a \cdot  v : = \pi_V(a) v,\, a \in \G, v \in V$ and define
$\beta_V:\, V\pfeil \G\tp V; \quad
   v\mapsto   v^\eb \tp v^\zb : = \DD_V \big(\eG\tp v \big)$. 
With this notation we get the following Corollary:
\begin{cor}
The conditions (i)-(iii) of Proposition \ref{thm32} are equivalent to
the following three conditions for $\beta_V$ (as before denoting
$P^i\tp Q^i \tp R^i = \phi^{-1}$): 
\begin{itemize}
  \item[(i')] 
$ (\ep\tp\id_V)\circ \beta_V = \id_V $
 \item[(ii')]
$      (a_{(2)}\cdot   v)^\eb  a_{(1)} \tp (a_{(2)}\cdot v)^\zb 
   =  a_{(2)}v^\eb \tp a_{(1)} \cdot v^\zb, \quad \forall v \in V $
  \item[(iii')]
$  R^i v^\eb  \tp (Q^i \cdot v^\zb)^\eb \,
        P^i \tp (Q^i \cdot v^\zb)^\zb  =    $ \\
\hfill $\Fii[321]{}\cdot \Big[ {(R^i\cdot v)^\eb}_{(2)}\,Q^i 
  \tp {(R^i \cdot v)^\eb}_{(1)} \, P^i \tp (R^i\cdot v)^\zb
  \Big],\quad \forall v \in V$
\end{itemize}
\end{cor}
\begin{proof}
  The equivalences (i) $\Leftrightarrow$ (i') and (ii)
$\Leftrightarrow$ (ii') are obvious. The equivalence (iii)
  $\Leftrightarrow$ (iii') follows by multiplying (iii) with
  $\Fii[312]{V}$ from the left and with $\phi^{-1}_V$ from the right
  and permuting the first two tensor factors.  
\end{proof}
The conditions stated in the above Corollary agree with those formulated
in [M2, Prop.2.2] by taking $\G^{cop}\equiv
(\G,\cop^{op},\Fii[321]{})$ 
instead of $(\G, \cop, \phi)$ as the underlying quasi-bialgebra. This
means that we have identified the category $\rep\D(\G)$ with what is
called the double category of modules over $\G$ in [M2].


\section{Doubles of weak quasi--Hopf algebras}
Allowing the coproduct $\cop$ to be non--unital (i.e.\ $\cop(\e)\neq
\e\tp\e$) leads to the definition of {\it weak quasi--Hopf algebras} as
introduced by G.~Mack and V.~Schomerus in [MS].
In this Section we sketch how the construction of the quantum double $\D(\G)$
generalizes to this case. As it will turn out, there are only minor
adjustments to be made. The reason for this lies in the fact, that we
have used mostly graphical identities (i.e.\ identities in $\rep \G$)
to derive and describe our results. But since $\rep \G$ is a rigid
monoidal category also in the case of  {\it weak} quasi--Hopf algebras
$\G$, all graphical identities in Section 2.2 and 2.3 stay
valid. Thus the only adjustements required refer to
those points, where we have 
translated graphical identities into algebraic ones. \\

Following [MS] we define a {\it weak quasi--Hopf algebra}
$(\G,\e,\cop,\ep,\phi)$ to be an
  associative algebra $\G$ with unit $\e$, a non--unital algebra
  map $\cop : \G \pfeil \G\tp\G$, an algebra map $\ep: \G \pfeil
  \mathbb{C}$ and an element
  $\phi \in \G\tp\G\tp\G$ satisfying \eqref{eq11}-\eqref{eq13},
  whereas \eqref{eq14} is replaced by 
  \begin{equation}
    \label{3.1}
    (\id\tp\ep\tp\id)(\phi) = \cop(\e) 
  \end{equation}
and where in place of invertibility $\phi$ is supposed to have a
{\it quasi--inverse} $\Fib{} \equiv \phi^{-1}$ with respect to the
intertwining property \eqref{eq11}. By this we mean that $\Fib{}$
satisfies  $\phi\Fib{}\phi =
\phi$, $\Fib{}\phi\Fib{} = \Fib{}$ as well as 
\begin{equation}
  \label{3.2}
    \phi \, \Fib{} = (\id\tp\cop)(\cop(\e)),\quad 
   \Fib{}\,\phi = (\cop\tp\id)(\cop(\e))
\end{equation}
which implies the further identities
\begin{align}
  \label{3.4}
  (\id\tp\cop)(\cop(a)) &= \phi \, (\cop\tp\id)(\cop(a)) \, \Fib{},
  \quad \forall a \in \G \\
   \label{3.4'}
  \phi = \phi \, (\cop\tp\id)(\cop(\e)), &\quad \Fib{} = \Fib{}\,
  (\id\tp\cop)(\cop(\e))\\
  \label{3.3}
  (\id\tp\ep\tp\id)(\Fib{}) & = \cop(\e)
\end{align}
More generally we call an element $t \in \A$ an intertwiner between
two (possibly non--unital) algebra maps $\alpha,\beta : \G\pfeil \A$,
if 
\begin{equation*}
  t \, \alpha(a) = \beta(a)\, t,\, \forall a \in \G\quad \text{and} \quad
  t \, \alpha(\e) \equiv \beta(\e)\, t = t
\end{equation*}
In this case by a quasi--inverse of $t$ (with respect to this
intertwiner property) we mean the unique (if existing) element
$\bar{t}\equiv t^{-1} \in \A$ satisfying $\bar{t} t = \alpha(\e)$, $t
\bar{t} = \beta(\e)$ and $\bar{t} t \bar{t} = \bar{t}$. Note that this
implies
\begin{equation*}
  \bar{t}\, \beta(a) = \alpha(a)\, \bar{t}, \quad 
   \bar{t}\, \beta(\e) \equiv \alpha(\e)\, \bar{t} = \bar{t}
\end{equation*}
and therefore $t$ is also the quasi--inverse of $t^{-1}$.

A weak quasi--bialgebra is called {\it weak quasi--Hopf algebra}, if there
exists  a unital algebra antimorphism $S: \G\pfeil \G$ and elements
$\alpha,\beta\in 
\G$ satisfying \eqref{eq16} and \eqref{eq17}. We will also always
suppose that $S$ is invertible.

Furthermore, $\G$ is said to be quasitriangular if there exists an
element $R\in \G\tp\G$ satisfying \eqref{eq18}--\eqref{eq110} and
possessing a quasi--inverse $\bar{R} \equiv R^{-1}$ with respect to the
intertwining property \eqref{eq18}. 

With these substitutions Theorem A generalizes as follows
\vspace{0.25cm}
  \begin{sloppypar} \noindent{\bf Theorem B }{\it 
Let $(\G,\cop,\phi)$ be a finite dimensional weak quasi--Hopf algebra
with invertible antipode $S$. 
Assume $\D(\G) \supset \G$ to be a  weak quasi--Hopf
algebra extension  satisfying (i)-(iii) of Theorem~A.
Then $\D(\G)$ exists uniquely up to equivalence and the 
linear map $\mu: \dG \tp \G \pfeil \D(\G)$ 
  \begin{equation*}
    \mu(\vi\tp a) : = (\id\tp \vi_{(1)})(q_\rho)\, D(\vi_{(2)})
  \end{equation*}
is surjective with $\makebox{\rm Ker}\, \mu = \makebox{\rm Ker}\,P$,
where $P:\dG\tp \G \pfeil \dG\tp \G$ is the linear
projection
\begin{equation*}
  P(\vi\tp a) : = \vi_{(2)} \tp\Big(\big(\hat{S}^{-1}(\vi_{(1)})
\arr \eG\big)   \arl \vi_{(3)}\Big)\, a =: \vi \reli a
\end{equation*}}
\end{sloppypar} \vspace{0.25cm}  

To adapt our previous strategy to weak quasi--Hopf algebras we first
recall that
due to the coproduct being non--unital the definition of the
tensor product functor in $\rep \G$ has to be slightly modified. First note
that the element $\cop(\e)$ (as well as higher coproducts of $\e$)
is idempotent and commutes with all elements in $\cop(\G)$. Thus,
given two representations $(V,\pi_V), (W,\pi_W)$, the operator
$(\pi_V\tp\pi_W)(\cop(\e))$ is a  projector,
whose image is precisely the $\G$--invariant subspace of $V\tp W$ on
which the tensor product representation operates non trivial. Thus one
is led to define the tensor product $\bo$ of two representations of $\G$
by setting
\begin{equation}
\label{bo}
  V\bo W : = (\pi_V \tp \pi_W)(\cop(\e))\, (V\tp W), \quad
  \pi_V\bo \pi_W : = (\pi_V \tp \pi_W)\circ \cop |_{ V\bo W}
\end{equation}
One readily verifies that with these definitions $\phi_{UVW}$ -
restricted to the subspace $(U\bo V)\bo W$ - furnish a natural family of
isomorphisms defining an associativity constraint for the tensor
product functor
$\bo$, where the tensor product of morphisms is
defined by restricting the ``usual'' tensor product map to the
truncated subspace.

With these adjustments, the graphical calculus described in Section
2.2 and 2.3 carries over to the present case. The collection of
colored upper (or  lower) legs
represent the (truncated) tensor product of $\G$--modules associated with the
individual legs. One just has to take care when translating the pictures
into algebraic identities. For example the graph
\begin{equation*}
  \begin{picture}(10,4)
 \put(0,1){\line(0,1){3}}
  \put(2,1){\line(0,1){3}}
\put(0,0.5){\makebox(0,0){$\G$}}
\put(2,0.5){\makebox(0,0){$\G$}}
\put(3.5,2){\makebox(0,0){$\equiv$}}
 \put(6,1){\line(0,1){3}}
\put(6,0.5){\makebox(0,0){$\G\bo\G$}}
  \end{picture}
\end{equation*}
is a pictorial representation of $\cop(\e)$ and not of $\e\tp\e$! Thus
the graph \eqref{graf5e} is equivalent to the algebraic identity 
$  [S(p^1_\la)\tp \e]\, q_\la \, \cop(p^2_\la) = \cop(\e)$
in place of the first equation of \eqref{eqvw}, etc. In this way all graphical
identities of Section 2 stay valid as 
well as Theorem \ref{SS}, where now $R^{-1}$ is meant to be the
quasi--inverse of $R$. \\

The definition of the diagonal crossed product in Proposition
\ref{prop311} yields an associative algebra which in general is
not unital, but $(\edG\tp \eG)$ is still a right unit and in
particular idempotent [HN, Thm.\ 14.2]. This may be cured by taking
the right ideal 
generated by $(\edG\tp\eG)$. Thus, let $P: \dG\tp \G \pfeil \dG\tp\G$
be the linear projection given by left multiplication with
$\edG\tp\eG$ with respect to the algebra structure \eqref{e312}, i.e. 
\begin{equation}
  \label{eq41}
  P(\vi\tp a) : = \vi_{(2)} \tp
  \big(\hat{S}^{-1}(\vi_{(1)})\re \e \li\vi_{(3)}\big) a  = 
  \e_{(-1)} \arr \vi\arl \Si(\e_{(1)}) \tp \e_{(0)}\, a.
\end{equation}
As in [HN, Sect.\ 14] we introduce the notation 
\begin{equation}
  \label{eq42}
  \vi\reli a : = P(\vi\tp a) \in \dG\tp \G
\end{equation}
and define the quantum double $\D(\G)$ as the subalgebra 
\begin{equation}
  \label{eq43}
  \D(\G) : = \dG\reli \G \equiv P(\dG\tp \G)
\end{equation}
Then $\edG \reli \eG \equiv \edG \tp \eG$ is the unit of $\D(\G)$ and
in terms of the notation \eqref{eq42} the multiplication in $\D(\G)$ is
still given by \eqref{e312}. In particular 
\begin{equation*}
  i_D: \G \ni a \mapsto \edG \reli a \equiv \edG \tp a \in \D(\G)
\end{equation*}
still provides a unital algebra inclusion. Interpreting Eq.\
\eqref{e313} also via \eqref{eq42}, Corollary \ref{lem3.3} likewise
extends to the present scenario. However note that now the definition
\eqref{e314} for left diagonal $\delta$--implementers $\LL \in \G\tp
\A$ also implies the nontrivial relation
\begin{equation}
\label{Lneu}
 \LL \equiv [\eG\tp\eA]\,\LL = [\Si(\e_{(1)}) \tp \eA]\,\LL\,[\e_{(-1)}\tp
  \tau(\e_{(0)})] ,
\end{equation}
This leads to a slight
modification of Proposition \ref{prop31.1} where one has to add the
requirement that $\cop$--flip operators $\TT$ fulfill also 
\begin{equation*}
  \TT\, (\id\tp\tau)(\cop(\e)) \equiv (\id\tp\tau)(\cop^{op}(\e))\,\TT
  = \TT
\end{equation*}
which follows directly from \eqref{Lneu} or by
multiplicating both sides of \eqref{L} from the right with
$(\id\tp\tau)(\cop(\e))$. 

Taking this additional identity into account, Theorem \ref{thm31} now reads
\begin{thm}
  \label{thm31w}
Using the notation \eqref{eq42} we 
 define the element $\DD \in \G\tp \D(\G)$ by 
 \begin{equation*}
   \DD: = \sum_\mu \Si(p^2_\rho)\, e_\mu \, p^1_{\rho (1)}\tp
     (e^\mu\reli p^1_{\rho(2)})
 \end{equation*}
Then the multiplication \eqref{e312} is the unique algebra
structure on $\D(\G)$ satisfying
\eqref{eq31}--\eqref{eq33} together with 
\begin{equation}
  \label{eq31w}
  \DD\, (\id\tp i_D)(\cop(\e)) \equiv (\id\tp i_D)(\cop^{op}(\e))\,
  \DD = \DD
\end{equation}
Moreover the identity \eqref{eq:3.15} also stays valid.
\end{thm}
The quasitriangular quasi--Hopf structure is now defined
precisely as in Theorem \ref{thm31a} and is proven analogously, where
in \eqref{eq310} $f^{-1}$ is the quasi--inverse of $f$. 
Correspondingly $\DD[-1]$ given by \eqref{Di} 
becomes the quasi--inverse of $\DD$ with respect to the $\cop$--flip
property \eqref{eq32}. Thus we arrive at a 
\begin{proof}[Proof of Theorem B] The existence parts (i),(ii) of
  Theorem B follow by putting as before $D(\vi) = (\vi\tp
  \id)(\DD)$. The universality (and therefore uniqueness) property
  (iii) follows analogously as in the proof of Theorem~A, Eq.\
  \eqref{eq:3.58'}. Here one just has to note that by Prop.~\ref{prop31.1}
  \begin{equation*}
    \LLt : = q^{op}_\rho \, \DDt \in \G\tp \tilde{\D}
  \end{equation*}
  is a normal coherent $\delta$--implementer, where $\DDt : = \sum e_\mu \tp
  \tilde{D}(e^\mu)$. Hence $\LLt$ satisfies \eqref{Lneu} and therefore
  the map $\nu: \dG\tp \G \pfeil \tilde{\D}$ 
  \begin{equation*}
    \nu(\vi \tp a) : = (\vi \tp \id)(\LLt) \, a
  \end{equation*}
satisfies $\nu \circ P = \nu$, where $P$ is the projection
\eqref{eq41}. Since as in Corollary \ref{lem3.3} the relations
\eqref{e314},\eqref{e315} guarantee that $\nu$ is an algebra map with
respect to the multiplication \eqref{e312} on $\dG\tp\G$, it passes
down to a well defined algebra map $f : \D(\G) \pfeil
\tilde{\D}$, $f(\vi\reli a) : = \nu(\vi\tp a)$, thus proving
(iii). Since $\DDt$ 
(as an element in $\tilde{\D}\tp \tilde{\D}$) is also required to be
quasitriangular, $f$ is even a quasi--bialgebra
homomorphism. In the case $\tilde{\D} = \D(\G)$ we have $\nu =
\mu$ and $\makebox{Ker}\,P =\makebox{Ker}\,\mu$ by definition, proving
also the second part of Theorem~B.
\end{proof}

%

%
%
%
%
\begin{appendix}
\section{The twisted double of a finite group}
As an application we now use Theorem \ref{thm31} and Theorem
\ref{thm31a}  to recover the
``twisted'' quantum double $\D^\om(G)$ of [DPR] where 
$G$ is a finite group and $\om : \,G\times G\times G \rightarrow U(1)$
is a normalized 3-cocycle. By definition this means $\om(g,h,k) = 1 $
whenever at least one 
of the three arguments is equal to the unit $e$ of G and
\begin{equation*}
  \om(g,x,y)\om(gx,y,z)^{-1}\om(g,xy,z)\om(g,x,yz)^{-1}\om(x,y,z)=1,\,\quad \forall g,x,y,z \in {\rm G}.
\end{equation*}
The Hopf algebra $\G:= Fun(G)$ of functions on $G$ may then also be viewed
as a quasi-Hopf algebra with its standard coproduct, counit and
antipode but with reassociator given by 
\begin{equation}
  \label{eq316}
  \phi : = \sum_{g,h,k \in G} \om(g,h,k) \cdot (\de_g \tp\de_h \tp \de_k),
\end{equation}
where $\de_g (x) := \de_{g,x}$. The identities (\ref{eq12}) and (\ref{eq14})
for $\phi$ are equivalent to $\om$ being a normalized 3-cocycle. Also
note that choosing $\alpha = \eG$ the antipode axioms now require
$\beta = \sum_g \om(g^{-1},g,g^{-1}) \delta_g$. In this special
example our quantum double $\D(\G) \equiv \dG \reli \G$ allows for
another identification with the linear space $\dG\tp \G$.
\begin{lem}
  \label{lemA1}
Let $\G$ be as above and define $\sigma: \dG\tp\G \pfeil \D(\G)$ by 
$\sigma(\vi\tp a) : = D(\vi)\,a,\, \vi \in \dG,a\in\G$. Then $\sigma$
is a linear bijection.
\end{lem}
\begin{proof}
  Since $(\G,\cop,\ep,S)$ is also an ordinary Hopf algebra, the
  relation \eqref{eq32} is equivalent to (suppressing the symbol
  $i_D$)
  \begin{equation}
    \label{eq:A1}
    a\,D(\vi) = D\big(a_{(1)} \arr \vi \arl \Si(a_{(3)})\big)\,
    a_{(2)},\quad \forall a \in \G,\vi \in \dG.
  \end{equation}
Using \eqref{eq:3.15} this implies
\begin{equation*}
  \vi\reli a \equiv (\id\tp \vi_{(1)})(q_\rho)\, D(\vi_{(2)})\, a =
  D\big( q^1_{\rho (1)} \arr \vi \arl (q^2_\rho\, \Si(q^1_{\rho
  (3)}))\big) \, q^1_{\rho (2)}\, a,
\end{equation*}
which lies in the image of $\sigma$. Hence, $\sigma$ is surjective and
therefore also injective.
\end{proof}
We note that in general the map $\sigma$ need not be surjective (nor
injective). Due to Lemma \ref{lemA1} we may now identify $\D(\G)$ with
the new algebraic structure on $\dG\tp\G$ induced by $\sigma^{-1}$. We
call this algebra $\dG\tp_D \G$. Putting $a \equiv \hat{\e} \tp_D a,\,
a \in \G$ and $\DD : = e_\mu \tp (e^\mu \tp_D \e) \in \G\tp (\dG\tp_D
\G)$ it is described by the relations \eqref{eq:A1},(\ref{eq33})  and the
requirement of $\G \equiv \hat{\e} \tp_D \G$ being a unital
subalgebra. To compute these multiplication rules we now use that the
group elements $g \in G$ provide a basis in $\dG$ with dual basis
$\delta_g \in \G$. Hence a basis of $\dG\tp_D \G$ is given by 
$\{h\tp_D\de_g\}_{h,g \in \G}$.
In this basis the generating matrix $\DD$ is given by
\begin{equation}
  \label{eq317}
   \DD = \sum_{k\in G} \de_k \tp (k\tp_D\eG), \quad \eG = \sum_{h\in G}\de_h.
\end{equation}
Let us know compute the multiplication laws according to the
definitions in Theorem \ref{thm31}. To begin with, we have 
\begin{equation*}
  (h \tp \eG )(e \tp \de_g)  = (h\tp \de_g) \quad \text{and} \quad
  (g \tp \eG )(h \tp \eG) = (gh \tp \eG).
\end{equation*}
Taking $(x\tp\id)$ of both sides of (\ref{eq32}), where $x\in G$, and 
using $\cop(\de_g) = \sum_{k\in G} \de_k \tp \de_{k^{-1}g}$ we get
\begin{equation*}
  (x\tp\eG)(e\tp\de_{x^{-1}g}) = (e\tp\de_{gx^{-1}})(x\tp\eG),
\end{equation*}
or equivalently
\begin{equation}
  \label{eq318}
  (e\tp\de_g)(x\tp\eG) = (x\tp\de_{x^{-1}gx}).
\end{equation}
Finally, pairing equation (\ref{eq33}) with $x\tp y \in \dG\tp\dG$ in
the two auxiliary spaces, the l.h.s. yields
\begin{multline*}
 \sum_{s,r,t \in G}\om(s,x,y) (\e\tp\de_s)(x\tp\eG) \cdot \om(x,r,y)^{-1}
   (e\tp\de_r)(y\tp\eG) \cdot \om(x,y,t)(\e\tp\de_t) \\
  = (x\tp\eG)(y\tp\eG) \cdot
    \sum_{s,r,t \in G} \frac{\om(s,x,y)\om(x,y,t)}{\om(x,r,y)} 
    (e\tp\de_{(xy)^{-1}sxy}\de_{y^{-1}ry}\de_t) \\
  = \sum_{t\in G} (x\tp\eG)(y\tp\eG)(\e\tp\de_t)\,
    \frac{\om(xyt(xy)^{-1},x,y)\om(x,y,t)}{\om(x,yty^{-1},y)},
\end{multline*}
where we have used (\ref{eq318}) in the first equality.
The right hand side of (\ref{eq33}) gives $(xy \tp \eG)$ so that we end up with
\begin{equation}
  \label{eq319}
  (x\tp\eG)(y\tp\eG) = \sum_{t\in G}
        \frac{\om(x,yty^{-1},y)}{\om(xyt(xy)^{-1},x,y)\om(x,y,t)}    
        (xy\tp\de_t).
\end{equation}
Similarly the coproduct is  computed as $\cop_D(e\tp\de_g) = \sum_{k\in
  G}(e\tp\de_k) \tp (e\tp\de_{k^{-1}g})$  and 
\begin{equation}
  \label{eq320}
  \cop_D (x\tp\eG) = \sum_{r,s \in G}  
  \frac{\om(xrx^{-1},x,s)} {w(x,r,s)\om(xrx^{-1},xsx^{-1},x)} \big( 
      (x\tp\de_r)\tp (x\tp\de_s)\big).  
\end{equation}
The above construction agrees with the definition of $\D^\om(G)$ given in
[DPR] up to 
the convention, that they have build $\D(\G)$ on $\G\tp\dG$ instead of
$\dG\tp\G$. 
%
%
%
%
%
\section{The Monodromy Algebra}
The definition of monodromy algebras (see e.g. [AFFS]) associated with
quasitriangular 
Hopf algebras may now easily be generalized to
the case of quasi-Hopf algebras. This has already been done in
[AGS]. We will give an explicit proof that the defining relations
of [AGS] indeed define an associative algebra structure on
$\dG\tp\G$, which in fact is isomorphic to our quantum
double $\D(\G)$. For ordinary Hopf algebras this has recently been
shown in [N1]. 

Let $\G$ be a finite dimensional quasi-Hopf algebra
with quasitriangular R-matrix $R \in \G\tp\G$. Following [N1] we define
the monodromy matrix $\MM \in \G\tp \D(\G)$ to be
\begin{equation*}
  \label{eq325}
   \MM : = (\id\tp i_D)(R^{op}) \, \DD. 
\end{equation*}
Defining also $\hat{R} \in \G\tp\G\tp \D(\G)$ by 
\begin{equation*}
  \label{eq326}
  \hat{R} := \Fi[213]{}  \, R^{12} \, 
                \phi^{-1},
\end{equation*}
we get the following Lemma:
\begin{lem}
  \label{lem32}
   The monodromy matrix $\MM$ obeys the following three conditions
   (dropping the symbol $i_D$):  
  \begin{align}
    \label{eq328}
       (\ep\tp\id)(\MM) &= \e_{\D(\G)},\\
    \label{eq329}
     \cop(a) \, \MM &= \MM \, \cop(a), \quad a\in \G \\
    \label{eq330}
      \MM[13]\,\hat{R}\,\MM[23] &= \hat{R} \, \phi \,(\cop\tp\id)(\MM)\,
     \phi^{-1} 
  \end{align}
\end{lem}
\begin{proof}
  We will freely suppress the embedding $i_D$. 
  Since the R-Matrix has the property $(\id\tp\ep)(R) = \e$, equation
  (\ref{eq328}) follows from $(\ep\tp\id)(\DD) =  \e_{\D(\G)}$. The
  identity \eqref{eq329} is implied by \eqref{eq32} and the
  intertwiner property of the R--Matrix.
  Let us now compute the l.h.s. of (\ref{eq330}):
  \begin{align*}
   \MM[13] \, \hat{R} \,\MM[23]& =
       R^{31}\, \DD[13] \phi^{213}\, R^{12}\,
       \phi^{-1}\,R^{32}\,\DD[23] \\
    & = R^{31}\, \DD[13]\, [(\cop\tp\id)(R)\phi^{-1}]^{132} \,\DD[23]
    \\
    & = \big[(R\tp\e) \cdot (\cop\tp\id)(R)\big]^{312}
         \,\DD[13]\,\Fii[132]{}\,\DD[23]
  \end{align*}
  where we have used the quasitriangularity of $R$ in the second
  line and property (\ref{eq32}) of $\DD$ in the third line. The r.h.s. of
  (\ref{eq330}) yields
  \begin{align*}
    \hat{R} \, \phi \,(\cop\tp\id)(\MM)\,
     \phi^{-1} &= \phi^{213}\,R^{12}\, (\cop\tp\id)(R^{op})\,
    \phi^{312}\,\DD[13]\, \Fii[132]{}\DD[23] \\
     & = \big[\phi^{321}\, (\e\tp R) \,
       (\id\tp\cop)(R)\,\phi\big]^{312} \, 
      \DD[13]\, \Fii[132]{}\DD[23],
  \end{align*}
  where we have used the definitions of $\MM$ and $\hat{R}$ and
  (\ref{eq33}). Now, the quasitriangularity of $R$ implies
  \begin{equation*}
  (R\tp\e) \, (\cop\tp\id)(R) = \phi^{321}\, (\e\tp R) \, (\id\tp\cop)(R) \,
   \phi 
  \end{equation*}
  which finally proves (\ref{eq330}).
\end{proof}
Note that the relations \eqref{eq328} - \eqref{eq330} are precisely
the defining relations postulated by [AGS] to describe the algebra
generated by the entries of a monodromy matrix around a closed loop
together with the quantum group of gauge transformations sitting at
the initial ($\equiv$ end) point of the loop. Thus we define similarly
as in [N1]
\begin{df}
  The {\bf gauged monodromy algebra} $M_R(\G) \supset \G$ is the algebra
 extension  generated by $\G$ and 
  elements $M(\vi),\, \vi\in\dG$ with defining
  relations given by \eqref{eq328} - \eqref{eq330}, where
  $M(\vi) \equiv (\vi\tp\id)(\MM)$.
\end{df}
Lemma \ref{lem32} then implies the immediate
\begin{cor} 
  Let $(\G,R)$ be a finite dimensional quasitriangular quasi--Hopf
  algebra. Then
  the monodromy algebra $M_R(\G)$ and the quantum double $\D(\G)$ are
  equivalent extensions of $\G$, where the isomorphism is given on the
  generators by 
  \begin{equation*}
 M(\vi) \leftrightarrow (\vi\tp\id)(R^{op}\,\DD)    
  \end{equation*}

\end{cor}
\end{appendix}


\end{document}